\pdfoutput=1
\RequirePackage{ifpdf}
\ifpdf 
\documentclass[pdftex]{sigma}
\else
\documentclass{sigma}
\fi

\usepackage{bbm}
\usepackage[all]{xy}
\usepackage{tikz, tkz-euclide}
\usetikzlibrary{calc,arrows,shapes,positioning}
\usetikzlibrary{decorations.markings}
\tikzstyle arrowstyle=[scale=1]
\tikzstyle directed=[postaction={decorate,decoration={markings,
 mark=at position .65 with {\arrow[arrowstyle]{stealth}}}}]
\tikzstyle reverse directed=[postaction={decorate,decoration={markings,
 mark=at position .65 with {\arrowreversed[arrowstyle]{stealth};}}}]

\definecolor{MichiganBlue}{HTML}{00274C}
\definecolor{MichiganYellow}{HTML}{FFCB05}
\definecolor{NicePurple}{RGB}{75,56,76} 
\definecolor{NiceRed}{RGB}{230,37,52}
\definecolor{MidnightBlue}{rgb}{0.1, 0.1, 0.44}

\DeclareMathOperator{\iiI}{I}

\DeclareMathOperator{\K}{K}
\DeclareMathOperator{\HH}{H}
\DeclareMathOperator{\Bes}{Bessel}
\newcommand{\iia}{\iiI_\alpha}
\newcommand{\Ka}{\K_\alpha}
\newcommand{\Ht}{\HH_\alpha^{(2)}}
\newcommand{\Ho}{\HH_\alpha^{(1)}}

\renewcommand{\vec}{\mathbf}

\newcommand{\lvertiii}{\left\vert\kern-0.25ex\left\vert\kern-0.25ex\left\vert}
\newcommand{\rvertiii}{\right\vert\kern-0.25ex\right\vert\kern-0.25ex\right\vert}
\newcommand{\dbar}{\overline{\partial}}
\newcommand{\siggg}{{\boldsymbol{\sigma}_3}}
\newcommand{\exop}{\sharp}
\newcommand{\eexop}{{\natural}}
\newcommand{\one}{\mathbbm{1}}

\numberwithin{equation}{section}

\newtheorem{Theorem}{Theorem}[section]
\newtheorem*{Theorem*}{Theorem}

\newtheorem{Lemma}[Theorem]{Lemma}
\newtheorem{Proposition}[Theorem]{Proposition}

\theoremstyle{definition}

\newtheorem{Remark}[Theorem]{Remark}

\let\Re=\undefined\DeclareMathOperator{\Re}{Re}
\let\Im=\undefined\DeclareMathOperator{\Im}{Im}

\DeclareMathOperator{\sign}{sign}

\newcommand*\dd{\mathop{}\!\mathrm{d}}
\newcommand{\ii}{\ensuremath{\mathrm{i}}}
\newcommand{\ee}{\ensuremath{\mathrm{e}}}

\newcommand{\OO}{\mathrm{O}}
\newcommand{\btau}{{\boldsymbol{\tau}}}
\newcommand{\ra}{\rightarrow}

\begin{document}

\allowdisplaybreaks

\newcommand{\arXivNumber}{2412.18656}

\renewcommand{\PaperNumber}{006}

\FirstPageHeading

\ShortArticleName{On the Asymptotics of Orthogonal Polynomials with Non-Analytic Weights}

\ArticleName{On the Asymptotics of Orthogonal Polynomials \\ on Multiple Intervals with Non-Analytic Weights}

\Author{Thomas TROGDON}

\AuthorNameForHeading{T.~Trogdon}

\Address{Department of Applied Mathematics, University of Washington, Seattle, WA, USA}
\Email{\mail{trogdon@uw.edu}}

\ArticleDates{Received April 05, 2025, in final form January 06, 2026; Published online January 28, 2026}

\Abstract{We consider the asymptotics of orthogonal polynomials for measures that are differentiable, but not necessarily analytic, multiplicative perturbations of Jacobi-like measures supported on disjoint intervals. We analyze the Fokas--Its--Kitaev Riemann--Hilbert problem using the Deift--Zhou method of nonlinear steepest descent and its $\dbar$ extension due to Miller and McLaughlin. Our results extend that of Yattselev in the case of Chebyshev-like measures with error bounds that give similar rates while allowing less regular perturbations. For the general Jacobi-like case, we present, what appears to be the first result for asymptotics when the perturbation of the measure is only assumed to be differentiable with bounded second derivative.}


\Classification{42C05; 33C47}

\section{Introduction}

Consider a Borel measure $\mu$ on $\mathbb R$ of the form
\begin{align}\label{eq:mu}
 \mu(\dd x) = \sum_{j=1}^{g+1} h_j(x) \mathbbm{1}_{[a_j,b_j]}(x) (b_j - x)^{\alpha_j} (x - a_j)^{\beta_j} \dd x + \sum_{j = 1}^P {r}_j \delta_{{c}_j} = \rho(x)\!\dd x + \sum_{j = 1}^P {r}_j \delta_{c_j},
\end{align}
where $a_j < b_j < a_{j+1}$ and $\{c_1,\dots,c_P\}$ is a subset of the complement of $\bigcup_j [a_j,b_j]$. We assume that $\alpha_j > -1$, $\beta_j > -1$, ${r}_j > 0$ for each $j$ and $h_j(x) > 0$ on $[a_j,b_j]$. We consider the problem of determining the strong asymptotics of the (monic) orthogonal polynomials, denoted by $\pi_n(;\mu)$, with respect to $\mu$ as $n \to \infty$ when $h_j$ is not assumed to be analytic. Polynomials orthogonal over disjoint intervals were originally considered by Akhiezer \cite{Akhiezer1960} and Akhiezer and Tomchuck \cite{Akhiezer1961} for a special class of weight functions. These authors described the orthogonal polynomials and their associated weighted Cauchy integrals, two of the main objects of study here. See \cite{Geronimo1988} for a method to study a different special class using polynomial mappings.

We are interested in obtaining asymptotic formulae for many reasons. One is in the use of the asymptotics in approximation theory. Specifically, the approximation of functions across multiple disjoint intervals has been found to be fruitful in computing matrix functions \cite{Ballew2023a}. Within random matrix theory there are other applications. First, the local asymptotics near the edge(s) and bulk(s) imply universality for certain invariant ensembles \cite{DeiftOrthogonalPolynomials,KuijlaarsInterval} and bounds on the growth of the polynomials can imply a certain stability of the recurrence coefficients that is useful in statistical estimation problems \cite{Ding2021,Ding2019}.

The so-called strong asymptotics of orthogonal polynomials was developed in the book of Szeg\H{o} \cite{Szego1939} for the classical families of orthogonal polynomials where one has a related differential equation and integral formula. The asymptotics are termed ``strong'' because they provide precise, pointwise, leading-order behavior. Two slightly more modern references are \cite{Levin2006,VanAssche1987}. The majority of these results focus on orthogonal polynomials on a single finite interval, on the entire real axis or on a semi-axis.

The Riemann--Hilbert approach to determining strong asymptotics was pioneered in \cite{DeiftOrthogonalPolynomials,DeiftWeights1,Kriecherbauer1999} for weights on the real axis and in \cite{KuijlaarsInterval} for the interval $[-1,1]$. Crucially, these works used that weight function $\rho$ in $\mu(\dd x) = \rho(x)\dd x$ had an appropriate analytic continuation to a~neighborhood of the support of $\mu$. It is worth noting here that the approach in \cite{Levin2006}, for example, does not require analyticity. {For a special choice of weight function, explicit formulae exist in terms of Riemann theta functions~\cite{Chen2008}.}

The Riemann--Hilbert approach was extended to non-analytic weights on the unit circle and the real line in \cite{McLaughlin2006,McLaughlin2008}. Then, in \cite{Baratchart2010}, the approach was used to obtain asymptotics for orthogonal polynomials on $[-1,1]$ for smooth, but not analytic, perturbations of the Jacobi weight $(1-x)^\beta (1 + x)^\alpha$. Continuing, in \cite{Yattselev2023} the asymptotics were refined in the case where~${\alpha,\beta \in \{-1/2,1/2\}}$. This refinement is important because if $\rho(x) = h(x) (1-x)^\beta (1 + x)^\alpha$, $\alpha,\beta \in \{-1/2,1/2\}$ and $h$ is analytic then one knows the error terms in the asymptotic expansions are exponentially small with respect to the degree of the polynomial \cite{Kuijlaars2003}. The work of \cite{Yattselev2023} uses a Riemann--Hilbert analysis to show that as the smoothness of $h$ increases, the order of the error term decreases accordingly. As noted in \cite{Yattselev2023}, this can also be seen to be a consequence of theorems concerning orthogonal polynomials on the unit circle, using the canonical (Joukowski) mapping of the unit circle to the unit interval $[-1,1]$.

For orthogonal polynomials for weight functions supported on multiple disjoint intervals of the form \eqref{eq:mu} for $\alpha_j,\beta_j \in \{-1/2, 1/2\}$, i.e., the Chebyshev-like case, the relevant works are first \cite{Yattselev2015} and then \cite{Ding2021}, where the quantities in~\cite{Yattselev2015} were made more explicit and used for perturbation theory. We are unaware of results for multiple intervals for general $\alpha_j$, $\beta_j$, i.e., the general Jacobi-like case, even for analytic perturbations. Thus, specifically, in this work we are interested in extending these results to the case where (i) $h_j$ in~\eqref{eq:mu} is not analytic and~(ii) considering general $\alpha_j, \beta_j > -1$. The goal is not specifically to obtain new formulae for the asymptotics in such situations, since the formulae for analytic perturbations $h_j$ broadly hold for non-analytic perturbations. Rather, one is interested in the size of the error terms in the asymptotics. For this reason, we largely leave our main results, Theorems~\ref{t:jacobi} and~\ref{t:cheb}, in an abstract form. If more specifics about the details of the asymptotics is desired, we refer the reader to~\cite{Ding2021} for asymptotics away from the edges and~\cite{KuijlaarsInterval} for local behavior near the edge.

The notation $C^{k}(I)$ is used to denote the space of $k$ times continuously differentiable functions~$f$ on a set $I$ with norm
\smash{$
 \|f\|_{C^{k}(I)} := \max_{0 \leq \ell \leq k} \max_{x \in I} \big|f^{(\ell)}(x)\big|$}.
And let $C^{k,\alpha}(I)$ denote the space of $k$ times continuously differentiable functions $f$ on a set $I$ such that $f^{(k)}$ is $\alpha$-H\"older continuous, with norm
\begin{align*}
 \|f\|_{C^{k, \alpha}(I)} := \max_{0 \leq \ell \leq k} \sup_{x \in I} \big|f^{(\ell)}(x)\big| + \sup_{y \neq x} \frac{\big|f^{(k)}(x) - f^{(k)}(y)\big|}{|x-y|^\alpha}.
\end{align*}
For general $\alpha_j$, $\beta_j$, we assume that $h_j \in C^{k,1}([a_j,b_j])$, $k \geq {1}$ {and $k > 2\max_j \max \{|\alpha_j|,|\beta_j|\} - \frac 1 2$} (see Theorem~\ref{t:jacobi}). We also utilize the Sobolev spaces $W^{k,p}([a_j,b_j])$ consisting of functions that have $k$ weak derivatives that are all in $L^p([a_j,b_j])$ with respect to standard Lebesgue measure.
 In the special case, where $\alpha_j, \beta_j \in \{1/2, -1/2\}$, we allow $h_j \in W^{k,q}([a_j,b_j])$, $k \geq 1$, $q > 4$. Below, we reserve the index~$k$ to refer to this assumed smoothness class.

As noted above, \cite{Yattselev2023} considered the case of a single interval with $\alpha_1, \beta_1 \in \{1/2, -1/2\}$. It was shown that the error term is, in particular, $\OO\bigl( \frac{\log n}{n^{k + \alpha}}\bigr)$ if $h_1\in C^{k,\alpha}([a_1,b_1])$ and $k \geq 3$. Here we obtain similar error bounds\footnote{In a previous version of this manuscript, an error bound without the $\log n$ factor was stated. While we believe this bound to still be valid, there was a flaw in the argument.} while allowing $g > 0$, $P > 0$, and only requiring $k \geq 1$, see Theorem~\ref{t:cheb}. There is a strong reason to believe this is still less than optimal, see Figures~\ref{fig:3_2} and~\ref{fig:2}. These plots seem to indicate that this bound is suboptimal by a factor of $1/n$ in the situation where the density has, in a sense, a single point of non-analyticity.

In the same vein as \cite{Yattselev2023}, the method described in this paper works by constructing close-to-analytic extensions of each function $h_j$ using polynomial approximation. Then, again borrowing directly from \cite{Yattselev2023}, we reuse an argument based on the Bernstein--Walsh inequality. Here it requires a bit more care in its implementation using a comparison, see Proposition~\ref{p:ggrowth}.

For the case of general $\alpha_j$, $\beta_j$, our approach differs moderately from the most closely related work in \cite{Baratchart2010}. In \cite{Baratchart2010}, the authors construct an approximately analytic extension of a Szeg\H{o}-like function that appears in the exponent. This has the benefit that one can directly use the exact classical parametrices to solve local problems. It appears to have the drawback that as the assumed smoothness is increased, one does not recover the order of the error term that is found in the case of analytic $h_j$. The method proposed here directly constructs a nearly analytic extension of $h_j$ and one obtains error terms on the order of $\OO\bigl(n^{-1 + \epsilon}\bigr)$ for any $\epsilon >0$ by supposing that $h_j$ is smooth enough. The drawback of our approach is that there is a small interval on the real axis where we do not obtain uniform estimates, determined by the $\delta$ parameter in Theorem~\ref{t:jacobi}. Importantly, $\delta \gg 1/n^2$ is always possible and the so-called Bessel asymptotics for the orthogonal polynomials near $a_j$, $b_j$ can be determined using the method here. If uniform asymptotics are desired using this method, we conjecture that a Sobolev improvement of our $L^2$ bounds is possible, with possibly greater smoothness assumptions on $h_j$.

The outline of the paper is as follows. In the remainder of this section, we introduce the Fokas--Its--Kitaev Riemann--Hilbert problem and fix more notation. In Section~\ref{s:dbar}, we present the~$\dbar$ deformations for the general Jacobi-type case $\alpha_j, \beta_j > -1$. This leads to Theorem~\ref{t:jacobi}. In Section~\ref{s:ext}, we first present an improvement for the Chebyshev-like case $\alpha_j, \beta_j \in \{-1/2,1/2\}$ in Theorem~\ref{t:cheb} and then discuss the addition of a finite number of point masses. Section~\ref{s:rec} discusses the implied asymptotics for recurrence coefficients and uses this to present numerical experiments demonstrating that Theorem~\ref{t:cheb} is likely suboptimal. We also include five appendices that include technical developments.

On a notational note, with possible subscripts, $C$, $C'$, $C''$, $c$, $c'$, $D$, $D'$ will be used to denote generic constants that may vary from line to line. Capital bold symbols $\vec A, \vec B, \dots$ will be used to denote matrices and matrix-valued functions. Lower-case bold characters are used to denote vectors. For concreteness, $\|\vec A\|$ denotes the Frobenius (Hilbert--Schmidt) matrix norm.

\subsection{The asymptotics of the recurrence coefficients}
We briefly review the implications of the asymptotic estimates obtained here by reviewing the formulae given in \cite{Ding2021b}. For error terms $E_1$, $E_2$, $E_3$, the recurrence coefficients for the orthonormal polynomials with respect to a measure $\mu$ \eqref{eq:mu} satisfy
\begin{gather*}
 b_n(\mu)^2 = \frac{1}{\mathfrak c^2} \frac{\displaystyle \frac{\Theta_2(\infty;\vec d_2;(n+1 - P) \vec{\Delta} + \boldsymbol{\zeta})}{\Theta_1(\infty;\vec d_2;(n+1-P) \vec{\Delta} + \boldsymbol{\zeta})} + E_1}{\displaystyle \frac{\Theta_2(\infty;\vec d_2;(n-P)\vec{\Delta} + \boldsymbol{\zeta})}{\Theta_1(\infty;\vec d_2;(n-P)\vec{\Delta} + \boldsymbol{\zeta})} + E_2},\\
a_n(\mu) = \frac{\Theta_1^{(1)}(\vec d_2; (n-P)\vec{\Delta} + \boldsymbol{\zeta})}{\Theta_1(\infty;\vec d_2;(n-P)\vec{\Delta} + \boldsymbol{\zeta})} - \frac{\Theta_1^{(1)}(\vec d_2; (n+1-P)\vec{\Delta} + \boldsymbol{\zeta})}{\Theta_1(\infty;\vec d_2;(n+1-P) \vec{\Delta} + \boldsymbol{\zeta})} + \mathfrak g_1 + E_3.
\end{gather*}
For the definition of $\Theta_j$ see \eqref{eq_vectortheta}, and the superscript $^{(1)}$ denotes the $\OO(1/z)$ term in the expansion of the function at $\infty$. The vectors $\vec d_j$, $\vec \Delta$, and $\boldsymbol{\zeta}$ are determined in \eqref{sec:global_para} and the constant $\mathfrak g_1$ is the $\OO(1/z)$ term in the expansion of $\mathfrak g$ at $\infty$ (see \eqref{eq:g}). The size of the error terms $E_1$, $E_2$, $E_3$ is determined by the size of the error terms, with respect to $n$, obtained in Theorem~\ref{t:jacobi} or Theorem~\ref{t:cheb}.

\subsection{The Fokas--Its--Kitaev Riemann--Hilbert problem}

In \cite{FokasOP}, the authors found a characterization of orthogonal polynomials in terms of a matrix Riemann--Hilbert (RH) problem. We now review such a formulation. Define the Cauchy transforms of the monic orthogonal polynomials\footnote{The monic orthogonal polynomials are defined by the requirement that $\pi_n(x;\mu) = x^n (1 + o(1))$, $n \to \infty$, and~${\int \pi_n(x;\mu) \pi_m(x;\mu) \mu( \dd x) \propto \delta_{n,m}}$.}
\[
 c_n(z;\mu) = \frac{1}{2 \pi \ii} \int_{\mathbb R} \frac{\pi_n(\lambda;\mu)}{\lambda - z} \mu(\dd \lambda),
\]
and the matrix-valued function
 \[
 \vec Y_n(z;\mu) = \begin{bmatrix} \pi_n(z;\mu) & c_n(z;\mu) \\
 \gamma_{n-1}(\mu) \pi_{n-1}(z;\mu) & \gamma_{n-1}(\mu) c_{n-1}(z;\mu) \end{bmatrix}, \qquad z \not \in \operatorname{supp}(\mu),
 \]
 where we use the notation
\smash{$
 \gamma_{n}(\mu) = - 2 \pi \ii \|\pi_n(\cdot;\mu)\|_{L^2(\mu)}^{-2}$}.
 It then follows that (see \cite{FokasOP} or \cite{Kuijlaars2003})
 \begin{align}\label{eq:def_Y_jump}
 \vec Y_n^{+}(z;\mu)&= \vec Y_n^{-}(z;\mu) \begin{bmatrix} 1 & \rho(z) \\ 0 & 1 \end{bmatrix}, \qquad \vec Y_n^{\pm}(z;\mu):=\lim_{\epsilon \to 0^+} \vec Y_n(z \pm \ii \epsilon;\mu),
 \end{align}
 at all points $z \in \mathbb R$ where $\mu$ has a continuous density $\rho$. From the discrete contributions to $\mu$,
 \begin{align} \label{eq:def_Y_res}
 \mathrm{Res}_{z = c_j} \vec Y_n(z;\mu) & = \lim_{z \to c_j} \vec Y_n(z;\mu) \begin{bmatrix} 0 & \frac{{r}_j}{2 \pi \ii} \\
 0 & 0 \end{bmatrix}, \qquad j = 1,2,\dots,P.
 \end{align}
We will initially suppose that $P = 0$ and discuss the requisite modifications later.

 Additionally,
 \begin{align}\label{e:def_Y_inf}
 \vec Y_n(z;\mu) \begin{bmatrix} z^{-n} & 0 \\ 0 & z^n \end{bmatrix} &= \vec I + \OO(1/z), \qquad z \to \infty.
 \end{align}

 We need to impose singularity conditions at the endpoints, following \cite{Kuijlaars2003}, we require the entry-wise asymptotics
 \begin{gather}
 \vec Y_n(z;\mu) = \begin{cases} \OO \begin{bmatrix} 1 & 1 + |z - b_j|^{\alpha_j} \\ 1 & 1 + |z - b_j|^{\alpha_j} \end{bmatrix}, & \alpha_j \neq 0,\vspace{1mm}\\
 \OO \begin{bmatrix} 1 & \log |z - b_j| \\ 1 & \log |z - b_j| \end{bmatrix}, & \alpha_j = 0,
 \end{cases}\qquad z \to b_j,\nonumber\\
 \vec Y_n(z;\mu) = \begin{cases} \OO \begin{bmatrix} 1 & 1 + |z - a_j|^{\beta_j} \\ 1 & 1 + |z - a_j|^{\beta_j} \end{bmatrix}, & \beta_j \neq 0,\vspace{1mm}\\
 \OO \begin{bmatrix} 1 & \log |z - a_j| \\ 1 & \log |z - a_j| \end{bmatrix}, & \beta_j = 0,
 \end{cases}\qquad z \to a_j.\label{eq:endpoint}
 \end{gather}
We make an important note here. First, it follows that $\det \vec Y_n(z;\mu) = 1$. Second, if we construct another function, $\tilde{\vec Y}_n(z;\mu)$ that also satisfies \eqref{eq:def_Y_jump}, \eqref{eq:def_Y_res}, \eqref{e:def_Y_inf} and \eqref{eq:endpoint} then the ratio~${
 \tilde{\vec Y}_n(z;\mu){\vec Y}_n(z;\mu)^{-1}}
$
will be analytic except for possibly $a_j$, $b_j$ where isolated singularities could persist. Indeed, we check, as $z \to b_j$, if $\alpha_j \neq 0$, for example,
\begin{align*}
 \tilde{\vec Y}_n(z;\mu){\vec Y}_n(z;\mu)^{-1} &= \OO \begin{bmatrix} 1 & 1 + |z - b_j|^{\alpha_j} \\ 1 & 1 + |z - b_j|^{\alpha_j} \end{bmatrix} \begin{bmatrix} 1+ |z - b_j|^{\alpha_j} & 1 + |z - b_j|^{\alpha_j} \\ 1 & 1 \end{bmatrix} \\
 & = \OO \begin{bmatrix} 1 + |z - b_j|^{\alpha_j}& 1 + |z - b_j|^{\alpha_j} \\ 1 + |z - b_j|^{\alpha_j} & 1 + |z - b_j|^{\alpha_j} \end{bmatrix}.
\end{align*}
From the condition $\alpha_j > -1$, we see that this isolated singularity must be removable. Liouville's theorem gives $ \tilde{\vec Y}_n(z;\mu){\vec Y}_n(z;\mu)^{-1} = \vec I$.

\begin{Remark} First, in all our calculations, $z = x + \ii y$ where $x, y \in \mathbb R$. We will treat a~function~$f(z)$ of a complex variable also a function two variables $f(x,y)$ and will abuse notation here.
 We typically write $f(x,y)$ when $f$ may not be an analytic function of $z = x + \ii y $.
\end{Remark}

\begin{Remark} We typically use $f^\pm(z)$ to denote the boundary values of a function $f$ from the left ($+$) or right ($-$) side of an oriented curve. In some cases, it may be convenient to use subscripts to denote the boundary values, e.g., $f'_\pm(z)$.
\end{Remark}

\section{Dbar deformation}\label{s:dbar}

We begin the deformation procedure by first recalling the definition and properties of the so-called $g$-function, see \eqref{eq:g}. {Historically speaking, this construction is found in \cite{Akhiezer1961}}. This function is also related to the exterior Green's function with pole at $\infty$ \cite{Walsh1935}. Often this function is determined during the process of deformation, but we need it at the outset. Throughout what follows, we now suppose that $\epsilon' > \epsilon > 0$ are sufficiently small so that both $\varphi(z;a_j)$, $\varphi(z;b_j)$ from~\eqref{eq:varphi} are conformal on the balls $B_{\epsilon'}(a_j)$, $B_{\epsilon'}(b_j)$, respectively,\footnote{Here we use the notation $B_\epsilon(c) = \{z \in \mathbb C\mid |z - c| < \epsilon\}$.} for all $j$. We assume that all the sets $B_{\epsilon'}(a_j)$, $B_{\epsilon'}(b_j)$, $j = 1,2,\dots,g+1$ are at a distance at least, say, $ 10\epsilon'$ from one another. Next, define
$
 \mathcal A_j = \varphi(B_\epsilon(a_j);a_j)$, $ \mathcal B_j = \varphi(B_\epsilon(b_j);b_j)$.
We also now suppose that $\epsilon > 0$ is sufficiently small so that these sets are convex.

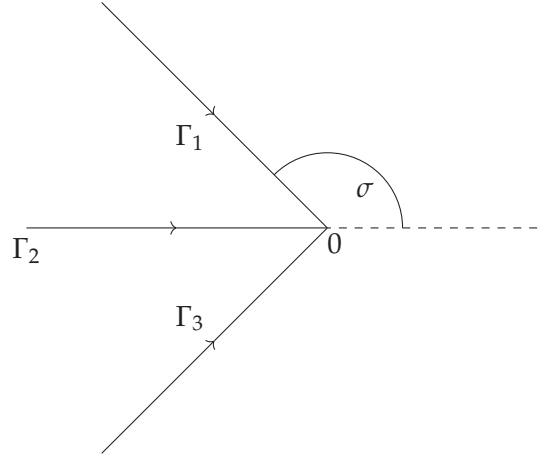
\begin{figure}[h!]
\centering\scalebox{.68}{
\begin{tikzpicture}

\coordinate (O) at (0,0);
\coordinate (A) at (-3,3);
\coordinate (Ah) at (-1.5,1.5);
\coordinate (B) at (3,0);
\coordinate (Bh) at (1.5,0);
\coordinate (C) at (-3,-3);
\coordinate (Ch) at (-1.5,-1.5);
\coordinate (D) at (-4,0);
\coordinate (Dh) at (-2,0);

\draw[-] (O) -- (Ah) node[anchor=north east] {$\Gamma_1$};
\draw[->] (A) -- (Ah);
\draw[-,dashed] (B) -- (O) ;
\draw[-] (O) -- (Ch) node[anchor=south east] {$\Gamma_3$};
\draw[->] (C) -- (Ch);
\draw[-] (O) -- (Dh);
\draw[<-] (Dh) -- (D) node[anchor=north] {$\Gamma_2$} ;

\draw (1,0) arc[start angle=0,end angle=135,radius=1cm];
\node at (0.5,0.5) {$\sigma$};

\node at (0.1,-0.2) {$0$};

\end{tikzpicture}}
\caption{Model contours. Here we use $\sigma = 2\pi/3$. The contours are symmetric about the real axis. \label{fig:model}}
\end{figure}

When $h_j$ is assumed analytic, one arrives at a local RH problem that needs to be solved explicitly near $a_j$, $b_j$. After employing the local conformal mappings, one finds a model problem with piecewise constant jump matrices. The solution of this problem is built out of Bessel and Hankel functions \cite{KuijlaarsInterval}. As will be made clearer below, the jump contours for the model problem are given in Figure~\ref{fig:model}. The contours will then need to be mapped using the inverse mappings. So, define
\smash{$
 \Sigma_{a_j,\ell} = \varphi^{-1}\bigl(\overline{\mathcal A_j} \cap \Gamma_\ell; a_j\bigr)$}, \smash{$ \Sigma_{b_j,\ell} = \varphi^{-1}\bigl(\overline{\mathcal B_j} \cap \Gamma_\ell; b_j\bigr)$},
for $j = 1,2,\dots, g+1$, $\ell = 1,2,3$.

\begin{figure}[!ht]
\centering

\scalebox{.67}{%
\begin{tikzpicture}

\coordinate (O) at (0,0);
\coordinate (A) at (3,3);
\coordinate (Ah) at (1.5,1.5);
\coordinate (B) at (3,0);
\coordinate (Bh) at (1.5,0);
\coordinate (C) at (3,-3);
\coordinate (Ch) at (1.5,-1.5);
\coordinate (D) at (4,0);
\coordinate (Dh) at (2,0);
\coordinate (E) at (4,2);
\coordinate (Eh) at (2,1);
\coordinate (F) at (4,-2);
\coordinate (Fh) at (2,-1);

\draw[-] (O) to[out=60,in=220] (Ah) ;
\draw[<-] (A) node[anchor=south east] {$\Sigma_{3,a_j}$} to[out=190,in=40] (Ah) ;
\draw[-,dashed] (B) -- (O) ;
\draw[-] (O) to[out=-60,in=140] (Ch) ;
\draw[<-] (C) node[anchor=south west] {$\Sigma_{1,a_j}$} to[out=170,in=-40] (Ch);
\draw[-] (O) -- (Dh);
\draw[->] (Dh) -- (D) node[anchor=north] {$\Sigma_{2,a_j}$} ;
\draw[fill=none](0,0) circle (3.5);
\node[anchor=north east] at (0,0) {$a_j$};

\draw[-,dashed] (O) -- (3.6/1.15,2.4/1.5);
\draw (1,0) arc[start angle=0,end angle=28,radius=1];
\node[anchor=south west] at (1,0) {\small $\theta_2$};

\draw[-,dashed] (O) -- (3.6/1.1,1.2/1.1);
\draw (2,0) arc[start angle=0,end angle=17.5,radius=2];
\node[anchor=south west] at (2,0) {\small $\theta_1$};

\end{tikzpicture}}
\scalebox{.67}{%
\begin{tikzpicture}

\coordinate (O) at (0,0);
\coordinate (A) at (-3,3);
\coordinate (Ah) at (-1.5,1.5);
\coordinate (B) at (-3,0);
\coordinate (Bh) at (-1.5,0);
\coordinate (C) at (-3,-3);
\coordinate (Ch) at (-1.5,-1.5);
\coordinate (D) at (-4,0);
\coordinate (Dh) at (-2,0);
\coordinate (E) at (-4,2);
\coordinate (Eh) at (-2,1);
\coordinate (F) at (-4,-2);
\coordinate (Fh) at (-2,-1);

\draw[-] (O) to[out=120,in=-40] (Ah) ;
\draw[<-] (A) node[anchor=east] {$\Sigma_{1,b_j}$} to[out=-10,in=140] (Ah) ;
\draw[-,dashed] (B) -- (O) ;
\draw[-] (O) to[out=240,in=40] (Ch) ;
\draw[<-] (C) node[anchor=east] {$\Sigma_{3,b_j}$} to[out=10,in=220] (Ch);
\draw[-] (O) -- (Dh);
\draw[->] (Dh) -- (D) node[anchor=south] {$\Sigma_{2,b_j}$} ;
\draw[fill=none](0,0) circle (3.5);
\node[anchor=south west] at (0,0) {$b_j$};

\end{tikzpicture}}
\caption{Mapping the contours $\Gamma_j$ using the local conformal mappings, with valid choices for $\theta_1$, $\theta_2$ displayed.\label{fig:mapped-1}}
\end{figure}
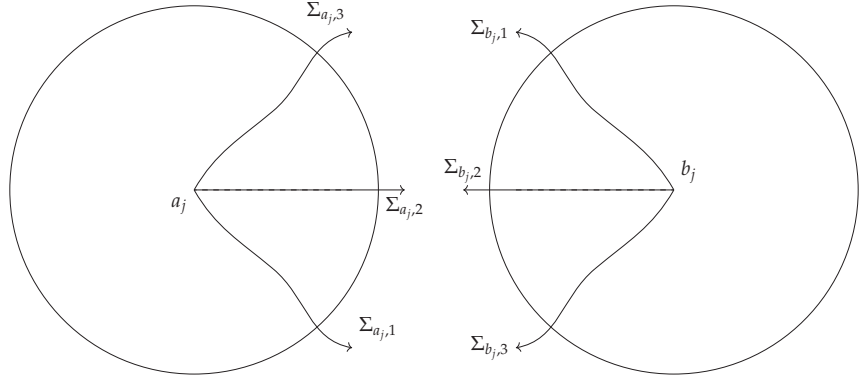

The contours in Figure~\ref{fig:mapped-1} will become jump contours in a~deformed hybrid RH-$\dbar$ problem.

For $z \in \mathbb C$, $z = r \ee^{\ii \theta}, \theta \in [0, 2\pi)$, define $(z)_\ra^\alpha := r^\alpha \ee^{\ii \alpha \theta}$. Similarly, for $z \in \mathbb C$, $z = r \ee^{\ii \theta}$, $\theta \in [-\pi, \pi)$, define $(z)^\alpha := r^\alpha \ee^{\ii \alpha \theta}$, {$\sqrt{z} = z^{1/2}$}.
Consider the function
$ w_j(x) = (b_j - x)^{\alpha_j} (x - a_j)^{\beta_j}$, $ a_j \leq x \leq b_j$,
and
\[
 \omega_j(z) = \ee^{- \ii \alpha_j \pi} (z - b_j)_\ra^{\alpha_j} (z - a_j)^{\beta_j}, \qquad z \in \mathbb C \setminus ((-\infty, a_j] \cup [b_j,\infty)).
\]
Then for $a_j < x < b_j$,
$
 \lim_{\epsilon \to 0} \omega_j(x + \ii \epsilon) = \ee^{- \ii \alpha_j \pi}(x - b_j)_\ra^{\alpha_j} (x - a_j)^{\beta_j} = w_j(x)$,
and therefore it is the analytic continuation of $w_j(x)$ to the upper- and lower-half planes.

In Appendix~\ref{a:extend}, we develop a \smash{$^\exop$} operator that extends a function off the real axis, while satisfying several convenient conditions. For $f_j(x) = 1/h_j(x)$, consider \smash{$f_j^{\exop}(x,y)$} with $\kappa = k$, $\tau = {\epsilon/2}$. The choice of $\theta_1$, $\theta_2$ that appear in the definition of \smash{$f_j^\exop$} will be discussed below.

\begin{figure}[t]
\centering

\scalebox{.63}{\begin{tikzpicture}

\coordinate (O) at (-8,0);
\coordinate (A) at (-5,3);
\coordinate (Ah) at (-6.5,1.5);
\coordinate (B) at (-5,0);
\coordinate (Bh) at (-7.5,0);
\coordinate (C) at (-5,-3);
\coordinate (Ch) at (-6.5,-1.5);
\coordinate (D) at (-4,0);
\coordinate (Dh) at (-6,0);
\coordinate (E) at (-4,2);
\coordinate (Eh) at (-6,1);
\coordinate (F) at (-4,-2);
\coordinate (Fh) at (-6,-1);

\draw[->] (-5.67,2.615) to (0,2.615);
\draw[-] (5.67,2.615) to (0,2.615);
\draw[-] (5.67,-2.615) to (0,-2.615);
\draw[->] (-5.67,-2.615) to (0,-2.615);

\draw[-] (O) to[out=60,in=220] (Ah) ;
\draw[<-] (-5.67,2.615) to[out=230,in=40] (Ah) ;
\node[anchor=south east] at (Ah) {$\Sigma_{3,a_j}$} ;
\draw[-,dashed] (B) -- (O) ;
\draw[-] (O) to[out=-60,in=140] (Ch) ;
\draw[<-] (-5.67,-2.615) to[out=130,in=-40] (Ch);
\node[anchor=north east] at (Ch) {$\Sigma_{1,a_j}$};
\draw[-] (O) -- (Dh);
\draw[>-] (Dh) node[anchor=north] {$\Sigma_{2,a_j}$} -- (D) ;
\draw[fill=none,dashed](-8,0) circle (3.5);
\node[anchor=north east] at (-8,0) {$a_j$};

\coordinate (O) at (8,0);
\coordinate (A) at (5,3);
\coordinate (Ah) at (6.5,1.5);
\coordinate (B) at (5,0);
\coordinate (Bh) at (6.5,0);
\coordinate (C) at (5,-3);
\coordinate (Ch) at (6.5,-1.5);
\coordinate (D) at (4,0);
\coordinate (Dh) at (6,0);
\coordinate (E) at (4,2);
\coordinate (Eh) at (6,1);
\coordinate (F) at (4,-2);
\coordinate (Fh) at (6,-1);

\draw[-] (O) to[out=120,in=-40] node[anchor=south west] {$\Sigma_{1,b_j}$} (Ah) ;
\draw[<-] (5.67,2.615) to[out=-50,in=140] (Ah) ;
\draw[-,dashed] (B) -- (O) ;
\draw[-] (O) to[out=240,in=40] node[anchor=north west] {$\Sigma_{3,b_j}$} (Ch) ;
\draw[<-] (5.67,-2.615) to[out=50,in=220] (Ch);
\draw[-] (O) -- (Dh);
\draw[<-] (Dh) node[anchor=south] {$\Sigma_{2,b_j}$} -- (D) ;
\draw[fill=none,dashed](8,0) circle (3.5);
\node[anchor=south west] at (8,0) {$b_j$};

\node at (0,2) {$\Sigma_{j,+}$};
\node at (0,-2) {$\Sigma_{j,-}$};

\node at (1,1) {$\Omega_{j,+}$};
\node at (-1,-1) {$\Omega_{j,-}$};

\draw[-] (-8,0) to (8,0) ;

\end{tikzpicture}}
\caption{The addition of contours to form $\Sigma_{j,\pm}$ and enclose the regions $\Omega_{j,\pm}$.
 The definition of $\Sigma_{j,\pm}$ includes the mapped contours near $a_j$, $b_j$. \label{fig:full}}
\end{figure}
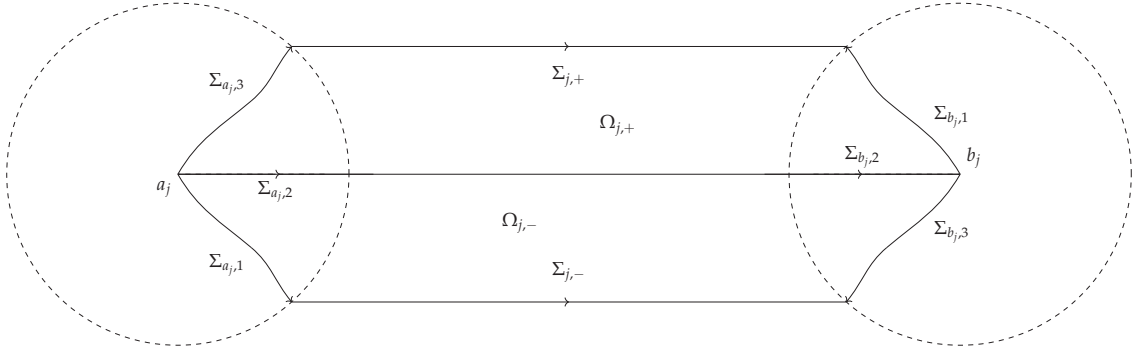

Before we state some properties of \smash{$f^\exop_j$}, we define a region where \smash{$f^\exop_j$} will be used. The contours~${\Sigma_{1,a_j}}$ and $\Sigma_{3,a_j}$ intersect the circle $\{|z - a_j| = \epsilon \}$ at points off the real axis, and to the right of the line $a_j + \ii \mathbb R$. The same is true of the intersection of the contours $\Sigma_{1,b_j}$ and $\Sigma_{3,b_j}$ and the circle $\{|z - b_j| = \epsilon \}$ to the left of $b_j + \ii \mathbb R$. So, a straight line, lying within the lower-half plane, can be used to connect
$
 \Sigma_{1,a_j} \cap \{|z - a_j| = \epsilon \} $ and $ \Sigma_{3,b_j} \cap \{|z - b_j| = \epsilon \}$,
giving a contour~$\Sigma_{j,-}$, {using left-to-right orientation}. And a straight line, lying within the upper-half plane, can be used to connect
$
 \Sigma_{3,a_j} \cap \{|z - a_j| = \epsilon \} $ and $ \Sigma_{1,b_j} \cap \{|z - b_j| = \epsilon \}$,
giving a contour $\Sigma_{j,+}$, {using left-to-right orientation}. This is depicted in Figure~\ref{fig:full}. The contours $\Sigma_{j,\pm}$ now enclose open regions $\Omega_{j,\pm}$, also depicted in Figure~\ref{fig:full}.

Define the $\dbar$ derivative
$
 \dbar = \frac 1 2 ( \partial_x + \ii \partial_y)$.
We have the following properties of \smash{$f_j^\exop$} (again, see Appendix~\ref{a:extend}) so long as we choose $\theta_2$ sufficiently small so that $\beta(x-a_j,y)$ and $\beta(x - b_j,y)$ in the definition of the \smash{$^\exop$} operator are identically equal to 1 on $\Sigma_{a_j,\ell}$, $\Sigma_{b_j,\ell}$, $\ell = 1, 3$. The extension~\smash{$f_j^\exop$} of $f_j$ then satisfies:
\begin{enumerate}\itemsep=0pt
 \item[$(1)$] for $a_j \leq x \leq b_j$, $f_j^\exop(x,0) = f_j(x)$,
 \item[$(2)$] for $c = a_j, b_j$, $|z - c| < \epsilon/2$, and $z \in \Sigma_{c,1} \cup \Sigma_{c,3}$, we have
 \[
 f_j^\exop(x,y) = f_j(z;c) = \sum_{\ell=0}^k \frac{f_j^{(\ell)}(c)}{\ell!} (z - c)^\ell,
\]
 \item[$(3)$] \smash{$f_j^\exop$} extends to be continuous on \smash{$\overline{\Omega_{j,+} \cup \Omega_{j,-}}$}, and
 \item[$(4)$] for $(x,y) \in {\Omega_{j,+} \cup \Omega_{j,-}}$, \smash{$\big|\dbar f_j^\exop(x,y)\big| \leq C \|f_j\|_{C^{k,1}([a_j,b_j)]}|y|^k$}, for some $C > 0$.
\end{enumerate}
{We further suppose that $\epsilon$ is sufficiently small so that $f_j(z;c)$ does not vanish for $|z-c|\leq \epsilon/2$.}

\subsection[The first dbar deformation of the FIK RH problem]{The first $\boldsymbol{ \dbar}$ deformation of the FIK RH problem}

To begin to deform the RH problem, we notice the factorization
\begin{align*}
 \begin{bmatrix} 1 & \rho(x) \\ 0 & 1 \end{bmatrix} = \begin{bmatrix} 1 & 0 \\ 1/\rho(x) & 1 \end{bmatrix} \begin{bmatrix} 0 & \rho(x) \\ -1/\rho(x) & 0 \end{bmatrix} \begin{bmatrix} 1 & 0 \\ 1/\rho(x) & 1 \end{bmatrix}.
\end{align*}
We will use this factorization locally near each interval $[a_j,b_j]$. Define
\begin{align*}
 \vec Z_n(x,y) = \vec Y_n(z ;\mu)\begin{cases} \begin{bmatrix} 1 & 0 \\ -f_j^\exop(x,y)/\omega_j(z) & 1 \end{bmatrix} ,& z \in \Omega_{j,+},\vspace{1.5pt}\\ \vspace{1.5pt}
 \begin{bmatrix} 1 & 0 \\ f_j^\exop(x,y)/\omega_j(z) & 1 \end{bmatrix},& z \in \Omega_{j,-},\\
 \vec I ,& \text{otherwise}.\end{cases}
\end{align*}
We compute
\begin{align*}
\dbar \vec Z_n(x,y) = \begin{cases}
\vec Y_n(z; \mu) \begin{bmatrix} 0 & 0 \\ -\dbar f_j^\exop(x,y)/\omega_j(z) & 0 \end{bmatrix} ,& z \in \Omega_{j,+},\vspace{1.5pt}\\ \vspace{1.5pt}
 \vec Y_n(z; \mu) \begin{bmatrix} 0 & 0 \\ \dbar f_j^\exop(x,y)/\omega_j(z) & 0 \end{bmatrix},& z \in \Omega_{j,-},\\
 \vec 0 ,& z \in \mathbb C \setminus \bigl( \overline \Omega_{j,+} \cup \overline \Omega_{j,-}\bigr).
\end{cases}
\end{align*}
This can actually be rewritten as
\begin{align*}
\dbar \vec Z_n(x,y) = \begin{cases}
\vec Z_n(z; \mu) \begin{bmatrix} 0 & 0 \\ -\dbar f_j^\exop(x,y)/\omega_j(z) & 0 \end{bmatrix}, & z \in \Omega_{j,+},\vspace{1.5pt}\\ \vspace{1.5pt}
 \vec Z_n(z; \mu) \begin{bmatrix} 0 & 0 \\ \dbar f_j^\exop(x,y)/\omega_j(z) & 0 \end{bmatrix},& z \in \Omega_{j,-},\\
 \vec 0 ,& z \in \mathbb C \setminus \bigl( \overline \Omega_{j,+} \cup \overline \Omega_{j,-}\bigr).
\end{cases}
\end{align*}
Setting
$
 \Sigma_j = \Sigma_{j,+} \cup \Sigma_{j,-}$,
as in Figure~\ref{fig:full}, we have the jump conditions
\begin{align*}
 \vec Z_n^+(x,y) = \vec Z_n^-(x,y) \begin{cases}
 \begin{bmatrix}
 1 & 0 \\ f_j^\exop(x,y)/\omega_j(z) & 1
 \end{bmatrix} ,& z \in \Sigma_j, \vspace{1.5pt}\\
 \begin{bmatrix} 0 & \rho(x) \\ -1/\rho(x) & 0 \end{bmatrix}, & z \in [a_j,b_j].
 \end{cases}
\end{align*}

Next, to handle the asymptotics at infinity, which we have not changed, we recall the definition of $\mathfrak g(z)$ in Appendix~\ref{a:gfun} and set
\smash{$
 \mathfrak c = \lim_{z \to \infty} \frac{\ee^{\mathfrak g(z)}}{z}$}.
We define
\begin{align*}
 \vec S_n(x,y) = \mathfrak c^{n \siggg} \vec Z_n(x,y) \ee^{-n \mathfrak g(z) \siggg}, \qquad \siggg = \begin{bmatrix}
 1 & 0 \\
 0 & -1
 \end{bmatrix}.
\end{align*}
We find that $\vec S_n(x,y)$ solves the normalized hybrid RH-$\dbar$ problem
\begin{gather*}
 \vec S_n^+(x,y) = \vec S_n^-(x,y) \begin{cases}
 \begin{bmatrix}
 1 & 0 \\ \ee^{-2n\mathfrak g(z)}f_j^\exop(x,y)/\omega_j(z) & 1
 \end{bmatrix} ,& z \in \Sigma_j, \vspace{1.5pt}\\
 \begin{bmatrix} 0 & \rho(x) \\ -1/\rho(x) & 0 \end{bmatrix} & z \in (a_j,b_j),\vspace{1.5pt}\\
 \ee^{- n\Delta_j \siggg}, & z \in (b_j,a_{j+1}),
 \end{cases}\\
 \dbar\vec S_n(x,y) = \vec S_n(z) \begin{cases}
\begin{bmatrix} 0 & 0 \\ - \ee^{-2n\mathfrak g(z)}\dbar f_j^\exop(x,y)/\omega_j(z) & 0 \end{bmatrix}, & z \in \Omega_{j,+},\vspace{1.5pt}\\ \vspace{1.5pt}
 \begin{bmatrix} 0 & 0 \\ \ee^{-2n\mathfrak g(z)}\dbar f_j^\exop(x,y)/\omega_j(z) & 0 \end{bmatrix},& z \in \Omega_{j,-},\\
 \vec 0 ,& z \in \mathbb C \setminus \bigl( \overline \Omega_{j,+} \cup \overline \Omega_{j,-}\bigr),
\end{cases}\\
\vec S_n(x,y) = \vec I + \OO\bigl(z^{-1}\bigr), \qquad z \to \infty.
\end{gather*}
We recall here that $\mathfrak g^+(x) + \mathfrak g^-(x) = 0$ for $x \in [a_j,b_j]$ and $\mathfrak g^+(x) - \mathfrak g_-(x) = {\Delta}_j$ for $x \in [b_{j}, a_{j+1}]$ where $\Delta_j$ is a constant.

\subsection[The second dbar deformation of the FIK RH problem]{The second $\boldsymbol{\dbar}$ deformation of the FIK RH problem}

While one can proceed to construct a parametrix for $\vec S_n(x,y)$ --- an approximate solution as~${n \to \infty}$, there is a fundamental difficulty. The outer solution will be determined, as we discuss further below, by a function that depends on the whole of $\rho$ --- the Szeg\H{o} function. The natural way to solve $\vec S_n(x,y)$ near $z = a_j, b_j$ is to use the Bessel parametrices that depend on~${f_j(z;c)}$, $c = a_j, b_j$, using its analyticity in a crucial way. The matching of this parametrix with the outer solution will not remove the jump on the real axis if $h_j$ does not have an analytic extension. The discrepancy can be estimated but it is not the most convenient way to handle the issue.

Define (using the notation in \eqref{eq:beps})
\begin{gather*}
 d_{j,-1}(x;c) = \left[\sqrt{\frac{f_j(x;c)}{f_j(x)}} - 1\right] {( 1 - b_\tau(x;c))}, \qquad c= a_j, b_j,\\
 d_{j,+}(x;c) = \left[\sqrt{\frac{f_j(x)}{f_j(x;c)}} - 1\right] {( 1 - b_\tau(x;c))}, \qquad c= a_j, b_j,
\end{gather*}
and set \smash{$d_{j,\pm}^\exop(x,y;c)$} to be the extension, as in Appendix~\ref{a:extend}, with $\tau = \delta < \epsilon/4$ and $\kappa = 0$ --- \smash{$d_{j,\pm}^\exop(x,y;c) = d_{j,\pm}(x;c)$} for $|y| < \epsilon/4$. {More explicitly, for example,
\begin{align*}
 d_{j,\pm}^\exop(x,y;a_j)={}& (1 - b_\delta(|y|;0))\bigl[\tilde b(x) d_{j,\pm}(x;a_j)+ (1- b_\delta(x;a_j) (1 - \beta(x - a_j,y)) d_{j,\pm}(x;a_j) \\
 &+ \beta(x-a,y) (1 - b_\delta(x;a_j)) d_{j,\pm}(a_j;a_j)\bigr].
\end{align*}
}
We write
\begin{align*}
 \left|\sqrt{\frac{f_j(x;c)}{f_j(x)}} - 1\right| = \left|\frac{\sqrt{f_j(x;c)}-\sqrt{f_j(x)}}{\sqrt{f_j(x)}}\right| \leq \frac{1}{2 \sqrt{\xi} \sqrt{f_j(x)}} | f_j(x;c)-f_j(x)|,
\end{align*}
for $\xi$ between $f_j(x;c)$ and $f_j(x)$. Since both of these functions are bounded below on $[0,2\delta]$, we find
\begin{gather*}
 \|d_{j,\pm}(\cdot;a_j)\|_{C^0([a_j,b_j])} \leq C \|f_j(\cdot;a_j)-f_j(\cdot) \|_{C^0([a_j,a_j+2\delta])},\\
 \|d_{j,\pm}(\cdot;b_j)\|_{C^0([a_j,b_j])} \leq C \|f_j(\cdot;a_j)-f_j(\cdot) \|_{C^0([b_j-2\delta,b_j])}.
\end{gather*}
Similarly, we find
\begin{gather*}
 \|d_{j,\pm}(\cdot;a_j)\|_{C^{1}([a_j,b_j])} \leq C \|f_j(\cdot;a_j)-f_j(\cdot) \|_{C^{0,1}([a_j,a_j+2\delta])},\\
 \|d_{j,\pm}(\cdot;b_j)\|_{C^{1}([a_j,b_j])} \leq C \|f_j(\cdot;b_j)-f_j(\cdot) \|_{C^{0,1}([b_j-2\delta,b_j])}.
\end{gather*}
We immediately have the following.
\begin{Lemma}\label{l:dest}
 \begin{align*}
 \big\|d_{j,\pm}^\exop(\cdot,y;c)\big\|_{L^\infty([a_j,b_j])} \leq C \delta^{k+1}, \qquad \big\|\dbar d_{j,\pm}^\exop(\cdot,y;c)\big\|_{L^\infty([a_j,b_j])} \leq C \delta^{k}.
 \end{align*}
\end{Lemma}

Next, define
\begin{align*}
 \vec T_n(x,y) = \vec S_n(x,y)\begin{cases}
 \begin{bmatrix} d_{j,{-}}^\exop(x,y;c) + 1 & 0 \\
 0 & d_{j,{+}}^\exop(x,y;c) + 1 \end{bmatrix}, & \begin{array}{@{}l@{}} z \in \Omega_{j,\pm}, \
 |z -c| < \epsilon,\\ c= a_j, b_j,\end{array}\\
 \vec I, & \text{otherwise}.
 \end{cases}
\end{align*}

We find that $\vec T_n(x,y)$ solves the hybrid RH-$\dbar$ problem for $c = a_j, b_j$
\begin{gather*}
 \vec T_n^+(x,y) = \vec T_n^-(x,y) \begin{cases}
 \begin{bmatrix}
 1 & 0 \\ \ee^{-2n\mathfrak g(z)}f_j^\exop(x,y)/\omega_j(z) & 1
 \end{bmatrix},& z \in \Sigma_j, \vspace{1.5pt}\\
 \begin{bmatrix} 0 & w_j(x)/f_j(x;c) \\ -f_j(x;c)/w_j(x) & 0 \end{bmatrix},& \begin{array}{@{}l@{}} |z - c| < \delta, \\ z \in (a_j,b_j),\end{array}\vspace{1.5pt}\\
 \begin{bmatrix} 0 & \rho(x) \\ -1/\rho(x) & 0 \end{bmatrix} ,& \begin{array}{@{}l@{}} |z - c| > 2\delta, \\ z \in (a_j,b_j),\end{array}\vspace{1.5pt}\\
 \begin{bmatrix} 0 & \rho(x) \frac{d_{j,-}^\exop(x,y;c) + 1}{d_{j,+}^\exop(x,y;c) + 1} \\ -\frac{d_{j,+}^\exop(x,y;c) + 1}{d_{j,-}^\exop(x,y;c) + 1} \rho(x)^{-1} & 0 \end{bmatrix} ,& \begin{array}{@{}l@{}} \delta \leq |z - c| \leq 2\delta, \\ z \in (a_j,b_j),\end{array}\\
 \ee^{- \Delta_j \siggg} ,& z \in (b_j,a_{j+1}),
 \end{cases}
\\
 \dbar\vec T_n(x,y) = \vec T_n(x,y) \vec W_n(x,y),\\
 \vec W_n(x,y) = \begin{cases}
 \begin{bmatrix} \frac{\dbar d_{j,+}^\exop(x,y;c)}{d_{j,+}^\exop(x,y;c) + 1} & 0 \\ - \ee^{-2n\mathfrak g(z)}\dbar f_j^\exop(x,y)/\omega_j(z) \frac{d_{j,+}^\exop(x,y;c) + 1}{d_{j,-}^\exop(x,y;c) + 1}\! & \frac{\dbar d_{j,-}^\exop(x,y;c)}{d_{j,-}^\exop(x,y;c) + 1} \end{bmatrix}\!,& \begin{array}{@{}l@{}} z \in \Omega_{j,+}, \\ |z - c| < \epsilon,\end{array}\vspace{1.5pt}\\ \vspace{1.5pt}
 \begin{bmatrix} \frac{\dbar d_{j,+}^\exop(x,y;c)}{d_{j,+}^\exop(x,y;c) + 1} & 0 \\ \ee^{-2n\mathfrak g(z)}\dbar f_j^\exop(x,y)/\omega_j(z)\frac{d_{j,+}^\exop(x,y;c) + 1}{d_{j,-}^\exop(x,y;c) + 1} & \frac{\dbar d_{j,-}^\exop(x,y;c)}{d_{j,-}^\exop(x,y;c) + 1}\end{bmatrix}\!,& \begin{array}{@{}l@{}} z \in \Omega_{j,-},\\ |z - c| < \epsilon, \end{array}\\
 \begin{bmatrix} 0 & 0 \\ - \ee^{-2n\mathfrak g(z)}\dbar f_j^\exop(x,y)/\omega_j(z) & 0 \end{bmatrix}\!, & \begin{array}{@{}l@{}} z \in \Omega_{j,+},\\ |z - a_j| > \epsilon, \\ |z - b_j| > \epsilon, \end{array}\vspace{1.5pt}\\ \vspace{1.5pt}
 \begin{bmatrix} 0 & 0 \\ \ee^{-2n\mathfrak g(z)}\dbar f_j^\exop(x,y)/\omega_j(z) & 0 \end{bmatrix}\!,& \begin{array}{@{}l@{}} z \in \Omega_{j,-}, \\ |z - a_j| > \epsilon, \\ |z - b_j| > \epsilon,\end{array} \\
 \vec 0,& z \in \mathbb C\! \setminus\! \bigl( \overline \Omega_{j,+} \!\cup\! \overline \Omega_{j,-}\bigr),
 \end{cases}\\
 \vec T_n(x,y) = \vec I + \OO\bigl(z^{-1}\bigr), \qquad z \to \infty.
\end{gather*}
And motivated by this, we define $\rho(z;c) := \omega_j(z)/f_j(z;c)$,
\begin{align*}
 \tilde \rho(x;\delta) := \begin{cases}
 \rho(x;c),& |x-c| < \delta, \ c = a_j, b_j,\\
 \rho(x),& |x - c| > 2\delta, \ c = a_1,b_1,\dots,a_{g+1},b_{g+1},\\
 \rho(x)\frac{d_{j,-}^\exop(x,y;c) + 1}{d_{j,+}^\exop(x,y;c) + 1},& \delta \leq |x-c| \leq 2 \delta,\quad c = a_j, b_j.
 \end{cases}
\end{align*}

\subsection{Construction of the global parametrix}\label{sec:global_para}
The value of the so-called Szeg\H{o} function $G(z)$ is that it allows one to replace the jumps $( a_j ,b_j)$ with something simpler at the cost of adding to the jumps on $( b_j, a_{j+1})$. Define\footnote{See Appendix~\ref{a:RS} for the definition of $R(z)$.}
 \[
 G(z) = -\frac{R(z)}{2 \pi \ii} \left[\sum_{j=1}^{g+1} \int_{a_j}^{b_j} \frac{\log \rho(x) }{x -z} \frac{\dd x}{R_+(x)} + \sum_{j=1}^{g} \int_{b_j}^{a_{j+1}} \frac{{\zeta}_j}{x -z} \frac{\dd x}{R(x)}\right],
 \]
 where the constants $\zeta_j$ are yet to be determined. Note that
 \begin{align*}
 G^+(z) + G^-(z) &= -\log \rho(z), \qquad z \in ( a_j, b_j),\\
 G^+(z) - G^-(z) &= -\zeta_j, \qquad z \in ( b_j, a_{j+1}).
 \end{align*}
 Since $R(z) = \OO\bigl(z^{g+1}\bigr)$, we see that $G(z) = \OO(z^{g})$. To avoid unbounded behavior of $G$ at infinity, we choose $\boldsymbol{\zeta} = (\zeta_j)_{j=1}^{g}$ so that as $z \to \infty$
$
 G(z) = \OO(1)$.
 Indeed, we find a linear system of equations
 \[
 m_\ell = -\sum_{j=1}^g \int_{a_j}^{b_j} \log \rho(x) x^{\ell-1} \frac{\dd x}{R_+(x)} - \sum_{j=1}^{g} \int_{b_j}^{a_{j+1}} \zeta_j x^{\ell-1} \frac{\dd x}{R(x)} = 0, \qquad \ell = 1,2,\dots,g.
 \]
 This system of equations is uniquely solvable for $\boldsymbol{\zeta}$ using the fact that the normalized differentials exist, and involves the same coefficient matrix that is used to determine the polynomials $Q_{g}$ in Appendix~\ref{a:RS}.

 As a first approximation, we consider the matrix $\tilde {\vec T}_n(z)$ that is obtained from $\vec T_n(x,y)$ by dropping the $\dbar$ conditions, and just using $\rho$ on the real axis
\begin{gather*}
 \tilde {\vec T}_n^+(z) = \tilde {\vec T}_n^-(z) \begin{cases}
 \begin{bmatrix}
 1 & 0 \\ \ee^{-2n\mathfrak g(z)} f_j^\exop(x,y)/\omega_j(z)& 1
 \end{bmatrix}, & z \in \Sigma_j, \vspace{1.5pt}\\
 \begin{bmatrix} 0 & \rho(x) \\ -1/\rho(x) & 0 \end{bmatrix}, & z \in (a_j,b_j),\vspace{1.5pt}\\
 \ee^{- \Delta_j \siggg}, & z \in (b_j,a_{j+1}),
 \end{cases}\\
\tilde {\vec T}_n(z) = \vec I + \OO\bigl(z^{-1}\bigr), \qquad z \to \infty.
\end{gather*}

 Then consider
\smash{$
 \vec U_n(x,y;\mu) = \ee^{\siggg G(\infty)} \tilde {\vec T}_n(z;\mu) \ee^{-\siggg G(z)}$}.
 We check the jumps of $\vec U_n$
\begin{align*}
 \vec U^+_n(z;\mu) = \begin{cases}
 \vec U^-_n(z;\mu) \begin{bmatrix} 1 & 0 \\ \ee^{-2(n\mathfrak g(z)-G(z))}f_j^\exop(x,y)/\omega_j(z) & 1 \end{bmatrix}, & z \in \Sigma_j \setminus \mathbb R, \vspace{4pt}\\
 \vec U^-_n(z;\mu) \begin{bmatrix} 0 & 1 \\ -1 & 0 \end{bmatrix}, & z \in ( a_j, b_j), \vspace{4pt}\\
 \vec U^-_n(z;\mu)) \begin{bmatrix} \ee^{-n\Delta_j-\zeta_j} & 0 \\ 0 & \ee^{n\Delta_j+\zeta_j} \end{bmatrix} ,& z \in ( b_{j}, a_{j+1}),
 \end{cases}
\end{align*}
and asymptotics
$
\vec U_n(z) = \vec I + \OO\bigl(z^{-1}\bigr)$, $ z \to \infty$.

Due to the exponential decay that $\ee^{-n \mathfrak g(z)}$ will induce, see Appendix~\ref{a:RS}, we expect the dominant contribution to the solution of this RH problem to come from the jump conditions on the real axis, at least away from the endpoints $a_j$, $b_j$.
From Appendix~\ref{a:RS}, \eqref{eq:L} specifically, we expect
$
 \vec U_n(x,y) \approx \vec L(z; n \boldsymbol{\Delta} + \boldsymbol{\zeta})$.
And so, the global parametrix to $\vec T_n$ is
\begin{align*}
 \vec G_n(z) = \ee^{-\siggg G(\infty)} \vec L(z; n \boldsymbol{\Delta} + \boldsymbol{\zeta}) \ee^{\siggg G(z)}.
\end{align*}
While this will indeed give the asymptotics of the orthogonal polynomials away from the endpoints, to complete the analysis, we need a modified version of this function. Define
\begin{gather*}
 \tilde G(z) = -\frac{R(z)}{2 \pi \ii} \left[\sum_{j=1}^{g+1} \int_{a_j}^{b_j} \frac{\log \tilde \rho(x;\delta) }{x -z} \frac{\dd x}{R_+(x)} + \sum_{j=1}^{g} \int_{b_j}^{a_{j+1}} \frac{{\zeta}_j}{x -z} \frac{\dd x}{R(x)}\right],\\
 \tilde{\vec G}_n(z) = \tilde{\vec G}_n(z;\delta)= \ee^{-\siggg G(\infty)} \vec L(z; n \boldsymbol{\Delta} + \boldsymbol{\zeta}) \ee^{\siggg \tilde G(z)}.
\end{gather*}
{Note that the constants $\zeta_j$ are the same as those in the definition of $G(z)$. This is so that the ratio $G(z)/\tilde G(z)$ will be analytic across the gaps $(a_{j+1},b_j)$.}

For a given $z \in \mathbb C \setminus[a,b]$ let $\eta(z)$ be the closest point in $[a,b]$ to $z$. Set $S(z) = \sqrt{z- a} \sqrt{z-b}$. Then
\begin{align*}
 S(z)\int_{a}^{b} \frac{f(x')}{x' -z} \frac{\dd x'}{S_+(x')} = f(\eta(z))\int_{a}^{b} \frac{1}{x' - z} \frac{S(z)}{S_+(x')} \dd x' + \int_{a}^{b} \frac{f(x')-f(\eta(z))}{x' - z} \frac{S(z)}{S_+(x')} \dd x'.
\end{align*}
The first integral can be computed explicitly and seen to be bounded on all of $\mathbb C$. If $L$ is the Lipschitz constant for $f$ on $[a,b]$ then because $|x'-\eta(z)| \leq |x' - z|$ the last integral is bounded by
\smash{$
L \int_{a}^b \frac{|S(z)|}{|S_+(x')|} \dd x'$},
which is bounded on bounded subsets of $\mathbb C$. Then using Lemma~\ref{l:dest}
\begin{align} \label{eq:diffG}
 |G(z) - \tilde G(z)| \leq C \sum_{j=1}^{g+1} \|\rho/w_j - \tilde \rho(\cdot,\delta)/w_j\|_{C^{0,1}([a_j,b_j])} \leq C' \delta^k,
\end{align}
where the constants depend on the subset.

\subsection{Using local solutions}
In this subsection, we make heavy use of the definitions and jump conditions established in Appendix~\ref{a:local}.
Consider the function
\begin{align*}
 \vec R_n(x,t) = \vec T_n(x,y) \bigl(\tilde{\vec G}_n(z) \vec G_n(z)^{-1} \vec P_n(z;b_j)\bigr)^{-1}.
\end{align*}
We claim that for $|z - b_j| < \epsilon$, that the function $\vec R_n$ is continuous for~${z \not \in (b_j - \epsilon, b_j)}$, and for~${|z - b_j| < \delta}$. Indeed, $\vec T_n(x,y)$ and $\vec P_n(z;b_j)$ have the same jump condition for $z \not\in (b_j - \epsilon, b_j)$ and $\tilde{\vec G}_n(z) \vec G_n(z)^{-1}$ is analytic for $z \not \in (b_j - \epsilon, b_j)$.

We then recall that for $b_j - \epsilon < z < b_j$ we have
\begin{align*}
 {\vec T}_n^+(x,y) &= {\vec T}_n^-(x,y)
 \begin{bmatrix} 0 & \tilde \rho(z;b_j) \\ -1/\tilde \rho(z;b_j) & 0 \end{bmatrix}.
\end{align*}
Using that for $b_j - \epsilon < z < b_j$
\begin{align*}
 \vec A_n^+(z;b_j) & = \vec A_n^-(z;b_j) \Psi^+(z)^{-\siggg} \vec E \left( \frac{\rho(z)}{\rho(z;b_j)} \right)^{\siggg} \vec E^{-1} \Psi^+(z)^\siggg,
\end{align*}
and
\begin{align*}
 \vec Q_n^+(z;b_j) & = \vec Q_n^-(z;b_j) \begin{bmatrix} 0 & \rho(z;b_j) \\ -1/ \rho(z;b_j) & 0 \end{bmatrix}, \qquad b_j - \epsilon < z < b_j,
\end{align*}
we compute, again for $b_j - {\epsilon} < z < b_j$, that
\begin{align*}
 \vec R_n^+(z) = {}&\vec R_n^-(z) \tilde{\vec G}_n^-(z) \vec G_n^-(z)^{-1}\vec A_n^-(z;b_j) \vec Q_n^-(z;b_j) \\
 & \times \left( \frac{\tilde\rho(z;\delta)}{\rho(z;b_j)} \right)^{\siggg}\vec Q_n^-(z;b_j)^{-1}\Psi^+(z)^{-\siggg} \vec E\left( \frac{\tilde\rho(z;\delta)}{\rho(z;b_j)} \right)^{-\siggg} \\
 & \times \vec E^{-1} \Psi^+(z)^\siggg\vec A^-_n(z;b_j)^{-1} \vec G_n^-(z) \tilde{\vec G}_n^-(z)^{-1}.
\end{align*}
For $b_j - \delta < z < b_j$, we have that $\tilde \rho(z;\delta) = \rho(z;b_j)$ and the claims about $\vec R_n$ follow.

Define
\begin{align*}
& \vec B_1(z)= \tilde{\vec G}_n^-(z) \vec G_n^-(z)^{-1}\vec A_n^-(z;b_j) \vec Q_n^-(z;b_j),\\
& \vec B_2(z)= \vec Q_n^-(z;b_j)^{-1}\Psi^+(z)^{-\siggg} \vec E,\\
& \vec B_3(z) = \vec E^{-1} \Psi^+(z)^\siggg\vec A^-_n(z;b_j)^{-1}\vec G_n^-(z) \tilde{\vec G}_n^-(z)^{-1}.
\end{align*}
The following is a direct consequence of Lemma~\ref{l:largez}.

\begin{Lemma}\label{l:d12}
 Let $b_j - \epsilon < z < b_j$, and
 \begin{align*}
 \frac{\tilde\rho(z;\delta)}{\rho(z;b_j)} = 1 + d_1(z), \qquad \frac{\rho(z;b_j)}{\tilde \rho(z;\delta)} = 1 + d_2(z), \qquad e(z) = \max\{|d_1(z)|, |d_2(z)|\},
 \end{align*}
 and define
 \begin{align*}
 {\varrho}_n(z;b_j) = \left\|\vec I - \vec B_1(z) \left( \frac{\tilde\rho(z;\delta)}{\rho(z;b_j)} \right)^{\siggg} \vec B_2(z) \left( \frac{\tilde\rho(z;\delta)}{\rho(z;b_j)} \right)^{-\siggg} \vec B_3(z)\right\|.
 \end{align*}
As $n^2(b_j-z) \to \infty$,
\smash{$
 \varrho_n(z;b_j) = \OO\bigl( n^{-1} |z-b_j|^{-1} e(z)^2 \bigr)$}.
\end{Lemma}

\subsection{Global approximation}

We have already defined $\vec R_n(z)$ for $|z - b_j| < \epsilon$. To complete the definition, we set
\begin{gather*}
 \vec R_n(x,y) = \vec T_n(x,y) \begin{cases} \bigl(\tilde{\vec G}_n(z) \vec G_n(z)^{-1} \vec P_n(z;b_j)\bigr)^{-1}, & |z - b_j| < \epsilon,\\
 \bigl(\tilde{\vec G}_n(z) \vec G_n(z)^{-1} \vec P_n(z;a_j)\bigr)^{-1}, & |z - a_j| < \epsilon,\\
\vec G_n(z)^{-1}, & \text{otherwise}.
 \end{cases}
\end{gather*}
We now make the restriction that $n^{-2} \ll \delta \ll \epsilon$ to find that, for some $c > 0$, $\vec R_n(x,y)$ should solve the hybrid $\dbar$-RH problem{\samepage
\begin{gather*}
 \vec R_n^+(x,y) = \vec R_n^-(x,y) \begin{cases}
 \vec I + \OO(\ee^{-c n}), & z \in \Sigma_{j,+} \cup \Sigma_{j,-},\ |z - b_j| > \epsilon, \ |z - a_j| > \epsilon,\\
 \vec I + \OO\bigl(n^{-1} + \delta^{k}\bigr) ,& |z - c| = \epsilon, \ c = a_j, b_j,\\
 \vec I + \OO\bigl(n^{-1}\delta^{-1}\bigr), & z \in (b_j- \epsilon,b_j-\delta) \cup (a_j + \delta,a_j+ \epsilon),
 \end{cases}\\
 \dbar\vec R_n(x,y)= \vec R_n(x,y) \vec X_n(x,y),\\
 \vec X_n(x,y)= \begin{cases}
 \vec G_n(z)
 \vec W_n(x,y)
 \vec G_n(z)^{-1} z \in \Omega_{j,\pm} \cap \bigl(\overline{B_\epsilon(a_j) \cup B_{\epsilon}(b_j)}\bigr)^c, \\
 \tilde{\vec G}_n(z) \vec G_n(z)^{-1}\vec P_n(z;c)
 \vec W_n(x,y)
 \vec P_n(z;c)^{-1} \vec G_n(z)\tilde{\vec G}_n(z)^{-1} \\
 \qquad \times z \in \Omega_{j,\pm} \cap B_\epsilon(c), \qquad c = a_j, b_j, \\
 \vec 0, \qquad z \in \mathbb C \setminus \bigl( \overline \Omega_{j,+} \cup \overline \Omega_{j,-}\bigr),
\end{cases}\\
\vec R_n(x,y) = \vec I + \OO\bigl(z^{-1}\bigr), \qquad z \to \infty,
\end{gather*}
where the jump condition error terms are uniform in $z$.}

The jump condition on $(b_j-\epsilon,b_j-\delta)$, for example, comes from Lemma~\ref{l:d12}, and the jump condition on $|z - c| = \epsilon$ follows from \eqref{eq:diffG} and the fact that $\vec G_n^{-1}(z) \vec P_n(z;c) = \vec I + \OO\bigl(n^{-1}\bigr)$ for~${|z-c| = \epsilon}$, see Appendix~\ref{a:local}.

To solve this problem, asymptotically, we first consider the problem with the RH component removed. We seek $\vec C_n(x,y)$ that is continuous on $\mathbb C$ and satisfies
\begin{gather*}
\dbar\vec C_n(x,y) = \vec C_n(z) \vec X_n(x,y), \qquad z \in \Omega_{j,+} \cup \Omega_{j,-},\\
\vec C_n(x,y) = \vec I + \OO\bigl(z^{-1}\bigr), \qquad z \to \infty,
\end{gather*}
where derivatives are understood to hold in a distributional sense. We set
\smash{$
 \Omega = \bigcup_{o = \pm} \bigcup_j \Omega_{j,o}$}.

\begin{Lemma}\label{l:dbar_est}
For $p > 2$, $k \geq 0$, suppose
\begin{align*}
 2p({2}|\alpha_j| - k) + p - 4 < 0, \qquad 2p({2}|\beta_j| - k) + p - 4 < 0,
\end{align*}
for all $j$. Then
\begin{gather*}
 \|\vec X_n\|_{L^p(\Omega)} = \OO\Bigl( {(1 + \one_{\eta = 0} \log n )}n^{2(\eta - k) + 1 - \frac 4 p} + n^{- \frac{((3 +k)p -1)((2k -1) p + 4)}{(2 k+5) p^2}} + \delta^k\Bigr), \\
 \eta = 2\max_j \max\{|\alpha_j|,|\beta_j|\}.
 \end{gather*}
\end{Lemma}
\begin{proof}
 We define three regions
 \begin{gather*}
 \Omega_{b_j,\mathrm{I}} = (\Omega_{j,+} \cup \Omega_{j,-}) \cap \big\{z\mid n^2|z - b_j| \leq c\big\},\\
 \Omega_{b_j,\mathrm{II}} = (\Omega_{j,+} \cup \Omega_{j,-}) \cap \{z\mid c < n^2|z - b_j| < N\},\\
 \Omega_{b_j,\mathrm{III}} = (\Omega_{j,+} \cup \Omega_{j,-}) \cap \big\{z\mid N \leq n^2|z - b_j| < n^2 \epsilon\big\},
 \end{gather*}
 where $c$ is sufficiently small, but fixed, such that the estimates in Lemma~\ref{l:smallz} hold, and $N \ll n^2$. In $\Omega_{b_j,\mathrm{I}}$ and $\Omega_{b_j,\mathrm{II}}$ we simply estimate $|\ee^{- \mathfrak g(z)}| \leq 1$. Then for $\alpha_j \neq 0$, $z \in \Omega_{{b_j},\mathrm{I}}$ we have for~${D_{\mathrm I} > 0}$, using Lemma~\ref{l:smallz},
 \begin{align*}
 \|\vec X_n(x,y)\| \leq D_{\mathrm I} \bigl( n^{2|\alpha_j| + 1} |z - b_j|^{- |\alpha_j| } |y|^k
 + \delta^k\bigr)\leq D_{\mathrm I} \bigl( n^{2|\alpha_j| + 1} |z - b_j|^{- |\alpha_j| + k} + \delta^k\bigr),
 \end{align*}
 where we used that
 \begin{align*}
 W(z;b_j)^{-\siggg} \begin{bmatrix} 0 & 0\\ 1/\omega_j(z) & 0 \end{bmatrix} W(z;b_j)^{\siggg} = \begin{bmatrix} 0 & 0\\ W(z;b_j)^2/\omega_j(z) & 0 \end{bmatrix}
 \end{align*}
 is bounded. Then rescaling and using polar coordinates\footnote{For convenience, we will routinely use the fact that for $a,b \geq 0$, $p > 0$ that $(a + b)^p \leq C_p (a^p + b^p)$.}
 \begin{align*}
 \int_{\Omega_{b_j,\mathrm{I}}} \|\vec X_n(x,y)\|^p \dd x \dd y \leq D_{\mathrm I}^p n^{2 p ({2}|\alpha_j|-k) + p - 4} \int_0^{2\pi} \int_0^c r^{- p |\alpha_j| + pk + 1} \dd r \dd \theta + O\bigl(\delta^{kp} \bigr).
 \end{align*}
 From this, we find the conditions
$
 p (k - |\alpha_j|) + 1 > -1$, $
 2 p ({2}|\alpha_j|-k) + p - 4 < 0$, {where the first is required for integrability and the second is needed to obtain a meaningful estimate. The second is more restrictive.} Now, if $\alpha_j = 0$, we have
 \begin{align*}
 \|\vec X_n(x,y)\| \leq D_{\mathrm I} \bigl(n|z - b_j|^{k} |\log n^2 | z -b_j| |^2 + \delta^k\bigr),
 \end{align*}
 giving
 \begin{align*}
 \int_{\Omega_{b_j,\mathrm{I}}} \|\vec X_n(x,y)\|^p \dd x \dd y &\leq D_{\mathrm I}^p n^{-2 p k + p - 4}
 (\log n)^2 + O\bigl( \delta^{kp}\bigr),\qquad
 -2 p k + p - 4 < 0.
 \end{align*}
 Then, again, rescaling and using polar coordinates {and Lemma~\ref{l:midz},}
 \begin{align*}
 \int_{\Omega_{b_j,\mathrm{II}}} \|\vec X_n(x,y)\|^p \dd x \dd y &\leq D_{\mathrm{II}}^p \left[ N^{3p/2}n^{-4-2kp} \int_0^N r^{kp + 1} \dd r\right] + \OO\bigl( \delta^{kp}\bigr)\\
 &= \frac{D_{\mathrm{II}}^p}{kp + 2} N^{3p/2+kp + 2}n^{-4-2kp} + \OO\bigl( \delta^{kp}\bigr).
 \end{align*}

 Now, in $\Omega_{b_j,\mathrm{III}}$, the change of variables will not be sufficient. So, we need an estimate on $\mathfrak g(z)$. For any $\epsilon > 0$, sufficiently small, there is $c' > 0$, such that for $ |y| \leq \epsilon$ we have
 \begin{align}\label{eq:gy}
 {\Re \mathfrak g(z) \leq - c' |y|,}
 \end{align}
 which, of course, implies $\big|\ee^{\mathfrak g(z)}\big| \leq \ee^{-c ' |y|}$. For $z \in \Omega_{b_j,\mathrm{III}}$,
 \begin{align*}
 \|\vec X_n(x,y)\|^p \leq D_{\mathrm{III}}^p \bigl(|y|^{pk} |x - b_j|^{-p/2}\ee^{- c' p n |y|} + \delta^{kp}\bigr).
 \end{align*}
 And then because $\Omega_{b_j,\mathrm{III}} \subset \big\{z\mid b_j + m^{-1} N^2/n^2 \leq x \leq b_j + \epsilon,\, |y| \leq \epsilon\big\}$ for some $m > 0$, we have, {using Lemma~\ref{l:largez}},
 \begin{align*}
 \int_{\Omega_{b_j,\mathrm{III}}} \|\vec X_n(x,y)\|^p \dd x \dd y \leq 2 D_{\mathrm{III}}^p \int_{b_j + m^{-1} N^2/n^2}^{b_j + \epsilon} \int_0^{\epsilon} y^{pk} \ee^{-c'pny} |x - b_j|^{-p/2} \dd x \dd y + O\bigl( \delta^{kp}\bigr).
 \end{align*}
 So, we set $x' = x - b_j$, $y' = c' n y$ giving
 \begin{align*}
 & \int_{\Omega_{b_j,\mathrm{III}}} \|\vec X_n(x,y)\|^p \dd x \dd y\\
 & \qquad {} \leq 2 D_{\mathrm{III}}^p (c' p n)^{-pk - 1} \left[\int_0^\infty y^{pk} \ee^{- y} \dd y\right] \frac{1}{1 - p/2}\phantom{\frac 1 2}x^{1 -p/2 } \bigg|_{m^{-1} N^2/n^2}^\epsilon
 + O\bigl( \delta^{kp}\bigr).
 \end{align*}
 We are left with, by possibly increasing $D_{\mathrm{III}}$,
 \begin{align*}
 \int_{\Omega_{b_j,\mathrm{III}}} \|\vec X_n(x,y)\|^p \dd x \dd y \leq D_{\mathrm{III}}^p n^{p(1-k) -3} N^{2 - p} + O\bigl( \delta^{kp}\bigr).
 \end{align*}
 To set $N$, we set
 \begin{gather*}
 n^{p(1-k) -3} N^{2 - p} = N^{3p/2+kp + 2}n^{-4-2kp},\qquad
 N^{3p/2+(k+1)p} = n^{1 + (k + 1)p},\\
 N = n^{\frac{1 + (k + 1)p}{3p/2 + (k + 1)p}}.
 \end{gather*}
 We note that this is a valid choice for $N$ because, for $p > 2$,
 \begin{align*}
 \frac{1 + (k + 1)p}{3p/2 + (k + 1)p} < 1.
 \end{align*}
 Then for $\Sigma = \Omega \setminus \bigcup_j (B_\epsilon(a_j) \cup B_\epsilon(b_j))$, {we use \eqref{eq:gy}, the formula for $\vec W_n$ and the boundedness of $\vec G_n$} to conclude
 \begin{align*}
 \int_{\Sigma} \|\vec X_n(x,y)\|^p \dd x \dd y \leq D \int_0^\epsilon y^{pk} \ee^{-np c' y} \dd y \leq D' n^{- pk - 1} .\tag*{\qed}
 \end{align*}\renewcommand{\qed}{}
\end{proof}

\begin{Proposition}\label{p:dbar}
 Suppose $p> 2$ and $p(k - {2}|\alpha_j|) + 2 > 0, p(k - {2}|\beta_j|) + 2 > 0$ for all $j$. Suppose also that $\vec V_n(x,y) \in L^\infty(\mathbb C)$ satisfies
 \begin{align}\label{eq:Un}
 \vec V_n(x,y) - \frac{1}{\pi} \int_\Omega \frac{\vec V_n(x',y')\vec X_n(x',y')}{z' - z} \dd A(z') = \frac{1}{\pi} \int_\Omega \frac{\vec X_n(x',y')}{z' - z} \dd A(z').
 \end{align}
 Then $\vec V_n(x,y) + \vec I$ is a solution of the $\dbar$ problem solved by $\vec C_n(x,y).$ Lastly, $\vec V_n$ is $\alpha$-H\"older continuous for some $0 < \alpha < 1$ and $\vec V_n|_{\Omega}$ is differentiable almost everywhere.
\end{Proposition}
\begin{proof}
 We use \cite[Theorem 4.3.10]{Astala2008} which states that, in a distributional sense, the integral operator
 \begin{align*}
 \mathcal K_\Omega u(x,y) = \frac{1}{\pi} \int_\Omega \frac{u(x',y')}{z' - z} \dd A(z')
 \end{align*}
 is inverse to $\dbar$, provided $u \in L^2(\mathbb C)$. The condition imposed on $k$ implies that $\vec X_n \in L^2(\mathbb C)$. So, if $\vec V_n \in L^\infty(\mathbb C)$, $\vec V_n \vec X_n \in L^2(\mathbb C)$ so that
$
 \dbar \vec V_n = \vec V_n\vec X_n + \vec X_n$, $
 \dbar (\vec V_n + \vec I) = (\vec V_n + \vec I)\vec X_n$.
 Next, from \cite[Theorem 4.3.13]{Astala2008}, we have that $\mathcal K_\Omega$ maps $L^p(\mathbb C)$ into the space of $\alpha$-H\"older continuous functions for $\alpha = 1 - 2/p$, $p > 2$. And we again use the relation
$
 \vec V_n = \mathcal K_\Omega ( \vec V_n \vec X_n) + \mathcal K_\Omega \vec X_n$,
 to get the desired conclusion because $\vec X_n \in L^p(\Omega)$. The claim about $\vec V_n|_{\Omega}$ follows from the elliptic regularity theorem, see \cite[Theorem 9.26]{Folland}, since $\vec V_n|_{\Omega} \in L^\infty(\Omega)$.
\end{proof}

We also see that for any $p > 2 > q, 1/p + 1/q = 1$, and $\vec U \in L^\infty(\mathbb C)$,
\begin{gather*}
 \| \vec K_n \|_{L^\infty(\mathbb C)} \leq C\| \vec U \|_{L^\infty(\mathbb C)} \|\vec X_n\|_{L^p(\Omega)} \left(\sup_{z \in \Omega} \int_{\Omega} \frac{\dd A(z')}{\pi|z' - z|^q}\right)^{1/q}, \\
 \vec K_n(x,y) = \frac{1}{\pi} \int_\Omega \frac{\vec U(x',y')\vec X_n(x',y')}{z' - z} \dd A(z').
\end{gather*}
We arrive at our main theorem concerning the existence of a solution $\vec C_n(x,y)$.
\begin{Theorem}
 Fix $k \geq 1$ and set $\eta = {2}\max_{j} \max\{|\alpha_j|,|\beta_j|\}$. Suppose $2(\eta -k) - 1 < 0$. Then for every $\gamma > 0$, and $n$ sufficiently large, there is a unique $L^\infty(\mathbb C)$ solution of \eqref{eq:Un} that satisfies
 \begin{align*}
 \|\vec V_n\|_{L^\infty(\mathbb C)} = \OO\bigl(n^{2(\eta -k) - 1 + \gamma} + \delta^{k} \bigr),
 \end{align*}
 which is also H\"older continuous on $\mathbb C$. Furthermore, $\vec V_n + \vec I$ solves the $\dbar$ problem for $\vec C_n$ satisfying the pointwise estimate
 \begin{align*}
 |\vec V_n(x,y)| = \OO\left( \frac{\|\vec V_n\|_{L^\infty(\mathbb C)}}{1 + |z|} \right), \qquad z \in \mathbb C.
 \end{align*}
\end{Theorem}
\begin{proof}
 {By choosing $p > 2$ in Lemma~\ref{l:dbar_est} sufficiently close to $2$, the theorem follows by Proposition~\ref{p:dbar}.}
\end{proof}

This implies that $\|\vec C_n(x,y)\|$ and $\| \vec C_n(x,y)^{-1}\|$ are uniformly bounded for $n$ sufficiently large. So, set
$
 \tilde {\vec R}_n(x,y) = \vec R_n(x,y) \vec C_n(x,y)^{-1}$.
It then follows that, for some $c > 0$, $\tilde{\vec R}_n(x,y)$ should solve the RH problem
\begin{align*}
 &\tilde{\vec R}_n^+(z) = \tilde{\vec R}_n^-(z) \begin{cases}
 \vec I + \OO(\ee^{-c n}), & z \in \Sigma_{j,+} \cup \Sigma_{j,-},\ |z - b_j| > \epsilon, \ |z - a_j| > \epsilon,\\
 \vec I + \OO\bigl(n^{-1} + \delta^{k}\bigr) ,& |z - c| = \epsilon, \ c = a_j, b_j,\\
 \vec I + \OO\bigl(n^{-1}\delta^{-1}\bigr) ,& z \in (b_j- \epsilon,b_j-\delta) \cup (a_j + \delta,a_j+ \epsilon),
 \end{cases}\\
&\tilde{\vec R}_n(z) = \vec I + \OO\bigl(z^{-1}\bigr), \qquad z \to \infty,
\end{align*}
where the jump condition error terms are uniform in $z$. Standard theory, see~\cite{DeiftOrthogonalPolynomials}, for example, gives estimates on $\tilde {\vec R}_n$. We find the following.

\begin{Theorem}\label{t:jacobi}
{Suppose $h_j \in C^{k,1}([a_j,b_j])$ for each $j$.} For $k\geq 1$ such that $2(\eta-k) < 1$, $\eta = {2}\max_{j} \max\{|\alpha_j|,|\beta_j|\}$ and any fixed $\gamma > 0$, sufficiently small, set
\begin{gather*}
\vec D^\eexop_{j}(z;c) =\begin{bmatrix} d_{j,+}^\exop(x,y;c) + 1 & 0 \\
 0 & d_{j,-}^\exop(x,y;c) + 1 \end{bmatrix},\\
 \vec D^\eexop(x,y) = \begin{cases}
 \vec D^\eexop_{j}(z;a_j)^{-1}\vec D^\eexop_{j}(z;b_j)^{-1},& z \in \Omega_{j,\pm},\\
 \vec I, & \text{otherwise},
 \end{cases}\\
 \vec F^\eexop(x,y) =\begin{cases}
 \begin{bmatrix} 1 & 0 \\
 \pm f_j^\exop(x,y)/\omega_j(z) & 1 \end{bmatrix}, & z \in \Omega_{j,\pm},\\
 \vec I, & \text{otherwise}.
 \end{cases}
\end{gather*}
Then for $n^{-1} \ll \delta \ll 1$,
\begin{gather*}
 \vec Y_n(z;\mu) =\mathfrak c^{-n \siggg} \vec T_n(x,y)\vec D^\eexop(x,y) \ee^{n \mathfrak g(z) \siggg} \vec F^\eexop(x,y) ,\\
 \vec T_n(x,y) = \tilde {\vec R}_n(z) (\vec V_n(x,y) + \vec I) \begin{cases} \tilde{\vec G}_n(z) \vec G_n(z)^{-1}\vec P_n(z;b_j), & |z - b_j| < \epsilon,\\
 \tilde{\vec G}_n(z) \vec G_n(z)^{-1}\vec P_n(z;a_j), & |z - a_j| < \epsilon,\\
 \vec G_n(z), & \text{otherwise},
 \end{cases}
\end{gather*}
where
\begin{gather*}
\tilde {\vec R}_n(z) = \vec I + \frac{1}{2\pi \ii} \int_{\Gamma_\delta} \frac{\vec O_n(z')}{z' - z} \dd z', \qquad \|\vec O_n\|_{L^2(\Gamma_\delta)} = \OO\bigl( n^{-1} + n^{-1}\delta^{-1} + \delta^k\bigr),\\
\vec V_n(x,y) = \OO\left( \frac{n^{-k +\gamma} + n^{2(\eta -k) - 1 + \gamma} + \delta^{k}}{1 +|z|} \right),
\end{gather*}
and therefore
\begin{gather*}
 \tilde {\vec R}_n(z) (\vec V_n(x,y) + \vec I) = \vec I + \OO \left( \frac{n^{-k + \gamma} + n^{2(\eta - k) - 1 + \gamma} + n^{-1} \delta^{-1} +\delta^k}{1 + |z|} \right),\\
 \tilde{\vec G}_n(z) \vec G_n(z)^{-1} = \vec I + \OO \bigl(\delta^k\bigr), \qquad \vec D^\eexop(x,y) = \vec I + \OO \bigl(\delta^k\bigr),
\end{gather*}
uniformly for $z$ in sets bounded away from $\Gamma_\delta$, where
\begin{align*}
 \Gamma_\delta = \bigcup_j \left( \{z\mid |z - a_j| = \epsilon \} \cup \{z\mid |z - b_j| = \epsilon \} \cup (b_j- \epsilon,b_j-\delta) \cup (a_j + \delta,a_j+ \epsilon) \right).
\end{align*}
\end{Theorem}
\begin{proof}
 Let
\smash{$
 \tilde{\vec Y}_n(z;\mu) =\mathfrak c^{-n \siggg} \vec T_n(x,y)\vec D^\eexop(x,y) \ee^{n \mathfrak g(z) \siggg} \vec F^\eexop(x,y)$}.
 All the estimates have been previously discussed and the remaining issue is the equality --- that $ \tilde{\vec Y}_n(z;\mu) = \vec Y_n(z;\mu)$. It can be shown that the jump conditions for $\tilde {\vec R}_n(z)$ satisfy the so-called product condition \cite[Definition~2.55]{TrogdonSOBook} and derivatives of the jump matrix are essentially bounded. Thus $\vec O_n \in H^1_z(\Gamma_\delta)$, \cite[Definition~2.48]{TrogdonSOBook}, an appropriate Sobolev space to ensure that~$\tilde {\vec R}_n(z)$, as defined, is uniformly $\alpha$-H\"older continuous and takes $\alpha$-H\"older continuous boundary values for some $\alpha > 0$. And~$\tilde{\vec Y}_n(z;\mu)$ satisfies the jump conditions set for the FIK RH problem, in a continuous sense (away from $a_j$, $b_j$).

 Next, we need to verify analyticity. But from classical elliptic regularity, \cite[Theorem 9.26]{Folland}, we have that $\tilde{\vec Y}_n(\cdot;\mu) \in C^\infty(\Omega')$ for every open set $\Omega'$ that does not intersect any one of the contours used in the deformation, i.e., for
 $$
 \Omega' \subset \mathbb C \setminus \Sigma, \qquad \Sigma := \bigcup_j(\Sigma_{j,-} \cup \Sigma_{j,+} \cup [a_j,b_j] \cup \{z\mid |z - a_j| = \epsilon\} \cup \{z\mid |z - b_j| = \epsilon\}).
 $$
 Furthermore, $\dbar \tilde{\vec Y}_n = 0$ in such a set $\Omega'$, implying analyticity.

 Lastly, from the singularity structure of $\vec P_{\Bes}$ it follows that $\tilde {\vec Y}_n$ has the same singularity orders as ${\vec Y}_n$ and therefore $\tilde{\vec Y}_n = {\vec Y}_n$.
\end{proof}

Effective use of this theorem then requires choosing $\delta$. If asymptotics are desired away from the endpoints $a_j$, $b_j$, one can choose $k$ so that $n^{-1}\delta^{-1} = \delta^k$, $\delta = n^{-1/(k + 1)}$. For $2(\eta -k) \leq 0$, this will give errors of $\OO\bigl( n^{- k/(k + 1)}\bigr)$. On the other hand, if one wants asymptotics near $z = a_j$, for example, care has to be taken in estimating $\tilde {\vec R}_n(z)$. We write $\Gamma_\delta$ as a disjoint union
\begin{align*}
 \Gamma_\delta = (a_j + \delta, a_j + \epsilon) \cup \Gamma_\delta^o,
\end{align*}
and find the entry-wise estimate for $z \in B_{\delta'}(z_j)$, $\delta' \ll \delta$,
\begin{align*}
 \big|\bigl(\tilde {\vec R}_n(z) - \vec I\bigr)_{ij}\big| & \leq \|\vec O_n\|_{L^2(\Gamma_\delta)} \left( \int_{a_j+\delta}^{a_j+\epsilon} \frac{\dd x}{|x - z|^2} \right)^{1/2} + \OO(\|\vec O_n\|_{L^2(\Gamma_\delta)})\\
 &= \OO\bigl( \delta^{-1/2} \|\vec O_n\|_{L^2(\Gamma_\delta)} \bigr).
\end{align*}
So, if $n^{-1} \delta^{-3/2} \ll 1$, or $\delta = n^{-2/3 + \gamma}$, and $2(\eta -k) + 1/3 < 0$, we obtain a valid error term giving uniform asymptotics on scaled neighborhoods of the endpoints $a_j$, $b_j$. Furthermore, by taking~$k$ sufficiently large, $\gamma < 2/3$ here can be chosen arbitrarily close to $2/3$.

\section{Extensions and improvements}\label{s:ext}

In this section, we discuss two topics. The first is the improvement in the estimates for the $\alpha_j, \beta_j \in \{1/2 , -1/2\}$, the perturbed Chebyshev-like case. The second is the addition of point masses to the measure $\mu$, i.e., $P > 0$.

\subsection{Chebyshev-like polynomials on multiple intervals}

So, suppose $\alpha_j, \beta_j \in \{1/2 , -1/2\}$. In the case where $f_j = 1/h_j$ is analytic in a neighborhood of~${[a_j,b_j]}$ for every $j$, local parametrices are not needed \cite{Kuijlaars2003} (see also \cite{Ding2021}). To treat the case of non-analytic $f_j$ using this fact, one needs to find an extension of $f_j$ to an entire neighborhood of $[a_j,b_j]$ to appropriately pose a hybrid $\dbar$-RH problem. This was accomplished in \cite{Yattselev2023}. While we follow this idea, we use a slightly different method of extension, since we have already developed our extension operator $^\exop$. First, we recall that for $\epsilon > 0$, $f_j \in W^{k,p}([a_j,b_j])$ has a~$W^{k,p}([a_j-2\epsilon,b_j +2\epsilon])$ extension.
 Similarly, $f_j \in C^{k,\alpha}([a_j,b_j])$ has a $C^{k,\alpha}([a_j-2\epsilon,b_j +2\epsilon])$ extension. This can be constructed simply by polynomial interpolation of the function and its first $k$ derivatives at the endpoints --- the Taylor polynomial. For what follows, we will assume~${f_j\in C^{k,\alpha}([a_j,b_j])}$, $k \geq 1$ and any $\alpha > 0$ or $f_j \in W^{k,q}([a_j,b_j])$, $k \geq 1$ and $q > 4$. Then, we find a~sequence of polynomials $(p'_{m,j})_{m \geq 1}$ of degree $m$ such that
\begin{align*}
 E_{j,m,q} := \| f_j' - p'_{m,j}\|_{L^q([a_j-2\epsilon,b_j+2\epsilon])} \to 0, \qquad q > 6.
\end{align*}
From Jackson's theorem, see \cite[Theorem 3.7.2]{atkinson}, we can select the sequence such that $E_{j,m,\infty} = \OO\bigl(m^{-k +1 - \alpha}\bigr)$ if $f \in C^{k,\alpha}([a_j,b_j])$, but more refined estimates are possible for $W^{k,q}([a_j,b_j])$ \cite{Canuto2006,Ditzian1987}. Specifically, if $f_j \in W^{k,p}([a_j,b_j])$, we can select the sequence such that \cite[(5.4.16)]{Canuto2006}
\begin{align*}
 E_{j,m,q} &\leq C m^{-k+1} \sum_{\ell = \min\{k-1,m+1\}}^{k-1} \|f_j^{(\ell+1)}\|_{L^q([a_j-2\epsilon,b_j +2\epsilon])} \\
 &\leq C' m^{-k +1} \|f_j\|_{W^{k,p}([a_j-2\epsilon,b_j +2\epsilon])}.
\end{align*}
Set, for example, \smash{$p_{m,j} = f_j(a_j) + \int_{a_j}^x p_{m,j}'(x') \dd x'$}. Then, define
$r_{m,j} = f_j - p_{m,j}$,
and its extension~\smash{$r_{m,j}^\exop(x,y)$} with $\kappa = 0$ and $\tau = \epsilon/2$. Then the extension of $f_j$ is given by
\begin{align*}
 f_j^\ddagger(x,y) = r_{m,j}^\exop(x,y) + p_{m,j}(z).
\end{align*}
Note that for $z \in \Omega_j' := \{z\mid \epsilon/2 - a_j < x < b_j + \epsilon/2,\, 0 < |y| < \epsilon/2\}$, we simply have
\begin{align*}
 f_j^\ddagger(x,y) = r_{m,j}(x) + p_{m,j}(z).
\end{align*}

In our deformations, we reuse the analogous notation from the previous sections and trust it will not cause too much confusion. Define $\Omega_{j,\pm}' := \Omega_j' \cap \mathbb C^\pm$ and
\begin{align*}
 \vec Z_n(x,y) = \vec Y_n(z ;\mu)\begin{cases} \begin{bmatrix} 1 & 0 \\ -f_j^\ddagger(x,y)/\omega_j(z) & 1 \end{bmatrix}, & z \in \Omega_{j,+}',\vspace{1.5pt}\\ \vspace{1.5pt}
 \begin{bmatrix} 1 & 0 \\ f_j^\ddagger(x,y)/\omega_j(z) & 1 \end{bmatrix},& z \in \Omega_{j,-}',\\
 \vec I ,& \text{otherwise}.\end{cases}
\end{align*}
So that
\begin{align*}
\dbar \vec Z_n(x,y) = \begin{cases}
\vec Z_n(z; \mu) \begin{bmatrix} 0 & 0 \\ -\dbar r_{m,j}^\exop(x,y)/\omega_j(z) & 0 \end{bmatrix} ,& z \in \Omega_{j,+}',\vspace{1.5pt}\\ \vspace{1.5pt}
 \vec Z_n(z; \mu) \begin{bmatrix} 0 & 0 \\ \dbar r_{m,j}^\exop(x,y)/\omega_j(z) & 0 \end{bmatrix},& z \in \Omega_{j,-}',\\
 \vec 0 ,& z \in \mathbb C \setminus \overline \Omega_j' .
\end{cases}
\end{align*}
Set
\begin{align*}
 \Sigma_{j,+}' &= \{z\mid x = a_j - \epsilon/2, \, 0 \leq y \leq \epsilon/2\} \\
 & =\cup \{z\mid y = \epsilon/2, \, a_j-\epsilon/2 \leq x \leq b_j + \epsilon/2\}
 \cup \{z\mid x = b_j + \epsilon/2, \, 0 \leq y \leq \epsilon/2\},
\end{align*}
i.e., the ``top'' boundary of the box $\Omega_j'$. We give this negative orientation (as the boundary of~$\Omega_j')$. And set $\Sigma'_{j,-} = \big\{z\mid \bar z \in \Sigma'_{j,+}\big\}$, inheriting orientation.

It follows that the matrix used in the definition of $\vec Z_n$ is continuous across the real axis for~${a_j - \epsilon/2 < x < a_j}$ and $b_j < x < b_j + \epsilon/2$ \cite{Kuijlaars2003} giving the resulting jump conditions
\begin{align*}
 \vec Z_n^+(x,y) = \vec Z_n^-(x,y) \begin{cases}
 \begin{bmatrix}
 1 & 0 \\ f_j^\ddagger(x,y)/\omega_j(z) & 1
 \end{bmatrix}, & z \in \Sigma_{j,-}' \cup \Sigma_{j,+}', \vspace{1.5pt}\\
 \begin{bmatrix} 0 & \rho(x) \\ -1/\rho(x) & 0 \end{bmatrix}, & z \in [a_j,b_j].
 \end{cases}
\end{align*}

As before, define
$
 \vec T_n(x,y) = \mathfrak c^{n \siggg} \vec Z_n(x,y) \ee^{-n \mathfrak g(z) \siggg}$.
We find that $\vec T_n(x,y)$ solves the normalized hybrid RH-$\dbar$ problem
\begin{gather*}
 \vec T_n^+(x,y) = \vec T_n^-(x,y) \begin{cases}
 \begin{bmatrix}
 1 & 0 \\ \ee^{-2n\mathfrak g(z)}f_j^\ddagger(x,y)/\omega_j(z) & 1
 \end{bmatrix}, & z \in \Sigma_{j,-}' \cup \Sigma_{j,+}', \vspace{1.5pt}\\
 \begin{bmatrix} 0 & \rho(x) \\ -1/\rho(x) & 0 \end{bmatrix}, & z \in (a_j,b_j),\vspace{1.5pt}\\
 \ee^{- \Delta_j \siggg} ,& z \in (b_j,a_{j+1}),
 \end{cases}\\
 \dbar\vec T_n(x,y) = \vec T_n(z) \begin{cases}
\begin{bmatrix} 0 & 0 \\ - \ee^{-2n\mathfrak g(z)}\dbar r_{m,j}^\exop(x,y)/\omega_j(z) & 0 \end{bmatrix}, & z \in \Omega_{j,+}',\vspace{1.5pt}\\ \vspace{1.5pt}
 \begin{bmatrix} 0 & 0 \\ \ee^{-2n\mathfrak g(z)}\dbar r_{m,j}^\exop(x,y)/\omega_j(z) & 0 \end{bmatrix},& z \in \Omega_{j,-}',\\
 \vec 0 ,& z \in \mathbb C \setminus \overline{\Omega_{j}'},
\end{cases}\\
\vec T_n(x,y) = \vec I + \OO\bigl(z^{-1}\bigr), \qquad z \to \infty.
\end{gather*}
With the notation for $\vec G_n(z)$ being precisely the same, we set
$
 \vec R_n(x,y) = \vec T_n(x,y) \vec G_n(z)^{-1}$,
and find that $\vec R_n(x,y)$ solves
\begin{gather*}
 \vec R_n^+(x,y) = \vec R_n^-(x,y)
 \vec G_n(z)\begin{bmatrix}
 1 & 0 \\ \ee^{-2n\mathfrak g(z) }f_j^\ddagger(x,y)/\omega_j(z) & 1
 \end{bmatrix} \vec G_n(z)^{-1},\qquad z \in \Sigma_{j,-}' \cup \Sigma_{j,+}',\\
 \dbar\vec R_n(x,y) = \vec R_n(z) \begin{cases}
 \vec G_n(z) \ee^{- \siggg G(z)} \begin{bmatrix} 0 & 0 \\ - \ee^{-2n\mathfrak g(z) - 2 G(z)}\dbar r_{m,j}^\exop(x,y)/\omega_j(z) & 0 \end{bmatrix}\\
 \qquad\times\ee^{\siggg G(z)}\vec G_n(z)^{-1}, \quad z \in \Omega_{j,+}',\vspace{1.5pt}\\ \vspace{1.5pt}
 \vec G_n(z) \ee^{- \siggg G(z)}\begin{bmatrix} 0 & 0 \\ \ee^{-2n\mathfrak g(z)- 2 G(z)} \dbar r_{m,j}^\exop(x,y)/\omega_j(z) & 0 \end{bmatrix}\\
 \qquad\times\ee^{\siggg G(z)}\vec G_n(z)^{-1},\qquad z \in \Omega_{j,-}',\\
 \vec 0, \qquad z \in \mathbb C \setminus \overline{\Omega_{j}'},
\end{cases}\\
\vec R_n(x,y) = \vec I + \OO\bigl(z^{-1}\bigr), \qquad z \to \infty.
\end{gather*}

As above, we then seek $\vec C_n(x,y)$ that is continuous on $\mathbb C$ and satisfies
\begin{gather*}
\dbar\vec C_n(x,y) = \vec C_n(z) \vec X_n(x,y), \qquad z \in \Omega_{j}',\\
\vec C_n(x,y) = \vec I + \OO\bigl(z^{-1}\bigr), \qquad z \to \infty,
\end{gather*}
where derivatives are understood to hold in an $L^2(\mathbb C)$ distributional sense. As above, we seek a~solution $\vec C_n = \vec V_n + \vec I$ where $\vec V_n$ satisfies
\begin{align*}
\vec V_n(x,y) - \frac 1 \pi \int_{\Omega'} \frac{\vec V_n(x',y') \vec X_n(x',y')}{z' - z} \dd A(z') = \frac 1 \pi \int_{\Omega'} \frac{\vec X_n(x',y')}{z' - z} \dd A(z'), \qquad \Omega' = \bigcup_j \Omega_j'.
\end{align*}
{Define the operator
$
 \mathcal L_n \vec U := \mathcal K_{\Omega'}( \vec U \vec X_n)$,
and the companion operator
\begin{align*}
 \underline{\mathcal K}_{\Omega'} u(z) = \frac{1}{\pi} \int_{\Omega'} \frac{u(x',y')}{|z'-z|} \dd A(z').
\end{align*}

Below, we will establish conditions under which $\|\underline{\mathcal K}_{\Omega'} (\|\vec X_n\|)\|_{L^\infty(\mathbb C)} \to 0$. This will then imply, via a standard Neumann series argument, that $\vec C_n(x,y)$ and $\vec C_n(x,y)^{-1}$ are uniformly bounded, for $n$ sufficiently large, and tend uniformly to the identity matrix. }

Next, we make a similar observation to what was made in \cite{Yattselev2023}. The Bernstein--Walsh inequality, see \cite[Theorem III.2.1]{SaffPotential}, in its simplest form, states that if a polynomial $p(x)$, of degree $n$, satisfies
\smash{$
 \max_{-1 \leq x \leq 1} |p(x)| \leq C$},
then
\begin{align*}
 |p(z)| \leq C |\phi(z)|^n, \qquad z \in \mathbb C, \qquad \phi(z) = z + \sqrt{z-1} \sqrt{z+1}
\end{align*}
Therefore, if a polynomial $p(x)$, of degree $n$, satisfies
\smash{$
 \max_{a \leq x \leq b} |p(x)| \leq C$},
then
\begin{align}\label{eq:has_M}
 |p(z)| \leq C |\phi(M(z))|^n, \qquad z \in \mathbb C, \qquad M(z) = M(z;a,b)= \frac{2}{b-a} \left( x - \frac{b + a}{2}\right).
\end{align}
Thus, for the sequence of polynomials $p_{m,j}(x)$ that converge uniformly on $[a_j,b_j]$, we have $|p_{m,j}(x)| \leq C$ for some constant, independent of $m$ (and $j$). From this, and Proposition~\ref{p:ggrowth}, provided that $\epsilon$ is sufficiently small we have that{\samepage
\begin{align}\label{eq:nm-bound}
 |p_{m,j}(z)| \big|\ee^{- 2n \mathfrak g(z)}\big| \leq C |\phi(M(z;a_j,b_j))|^{- 2 c n + m},
\end{align}
for some $c > 0$. So, we choose, for example $m = \lfloor c n \rfloor$ so that \eqref{eq:nm-bound} is $\OO( \ee^{- cn})$.}

Now, as before, set
$
 \tilde {\vec R}_n(x,y) = \vec R_n(x,y) \vec C_n(x,y)^{-1}$.
It then follows that $\tilde{\vec R}_n(x,y)$ should solve the RH problem
\begin{align*}
 &\tilde{\vec R}_n^+(z) = \tilde{\vec R}_n^-(z) ( \vec I + \OO( \ee^{- cn}) ), \qquad z \in \Sigma', \qquad \Sigma' = \bigcup_j ( \Sigma_{j,-}' \cup \Sigma_{j,+}')\\
&\tilde{\vec R}_n(z) = \vec I + \OO\bigl(z^{-1}\bigr), \qquad z \to \infty,
\end{align*}
where the jump condition error terms are uniform in $z$. Standard theory, see \cite{DeiftOrthogonalPolynomials}, for example, again gives estimates on $\tilde {\vec R}_n$, and we find the following.

\begin{Theorem}\label{t:cheb}
For $k\geq 1$, $q > 6$ and\footnote{Recall that $\epsilon$ appears in the definition of $\Sigma_{j,\pm}'$.} $\epsilon$ sufficiently small set
\begin{align*}
 \vec F^\eexop(x,y) &=\begin{cases}
 \begin{bmatrix} 1 & 0 \\
 \pm f_j^\ddagger(x,y)/\omega_j(z) & 1 \end{bmatrix}, & z \in \Omega_{j,\pm}',\\
 \vec I, & \text{otherwise},
 \end{cases}
\end{align*}
{and suppose $\|\underline{\mathcal K}_{\Omega'} (\|\vec X_n\|)\|_{L^\infty(\mathbb C)} = o(1)$, $n \to \infty$.} Then for $n$ sufficiently large,
\begin{align*}
 \vec Y_n(z;\mu) &=\mathfrak c^{-n \siggg} \vec T_n(x,y) \ee^{n \mathfrak g(z) \siggg} \vec F^\eexop(x,y) ,\\
 \vec T_n(x,y) &= \tilde {\vec R}_n(z) (\vec V_n(x,y) + \vec I) \vec G_n(z),
\end{align*}
where
\begin{align*}
\tilde {\vec R}_n(z) = \vec I + \frac{1}{2\pi \ii} \int_{\Sigma'} \frac{\vec O_n(z')}{z' - z} \dd z', \qquad \|\vec O_n\|_{L^2(\Sigma')} = \OO( \ee^{-c n}),
\end{align*}
{and
\begin{align*}
 \vec V_n(x,y) = \OO\left( \frac{\|\underline{\mathcal K}_{\Omega'} (\|\vec X_n\|)\|_{L^\infty(\mathbb C)}}{1 + |z|}\right), \qquad z \in \mathbb C, \qquad m = \lfloor c n \rfloor.
\end{align*}}
\end{Theorem}

\begin{proof}
 The only remaining detail to spell out is the estimate on $\vec V_n$. First, it follows that
 \begin{align*}
 \|\mathcal L_n\|_{L^\infty(\Omega')} \leq \|\underline{\mathcal K}_{\Omega'} (\|\vec X_n\|)\|_{L^\infty(\mathbb C)}.
 \end{align*}
 This implies that the integral equation
$
 \vec V_n - \mathcal L_n \vec V_n = \mathcal K \vec X_n
$
 is near-identity. For $n$ sufficiently large, it has a unique solution, tending to zero as $n \to \infty$, and
 \begin{align*}
 \|\vec V_n\| \leq \|\mathcal L_n \vec V_n\| + \|\underline{\mathcal K}_{\Omega'} \vec X_n\| \leq (1 + \|\vec V_n\|_{L^\infty(\Omega')}) \|\underline{\mathcal K}_{\Omega'} (\|\vec X_n\|)\|_{L^\infty(\mathbb C)}.\tag*{\qed}
 \end{align*} \renewcommand{\qed}{}
\end{proof}

To finalize the preceding theorem, we have the following estimates.

\begin{Lemma}[norm bounds]
For $p > 2$, $q > 6$, and $1/p = 1/q + 1/r$, and $m = \lfloor c n \rfloor$,
\begin{align*}
 \|\underline{\mathcal K}_{\Omega'} (\|\vec X_n\|)\|_{L^\infty(\mathbb C)} = \OO\bigl(n^{-1/r} \max_j E_{j,m,q}\bigr).
\end{align*}
Furthermore,
\begin{align*}
 \|\underline{\mathcal K}_{\Omega'} (\|\vec X_n\|)\|_{L^\infty(\mathbb C)} = \OO\bigl(n^{-1} \log n \max_j E_{j,m,\infty}\bigr).
\end{align*}
\end{Lemma}
\begin{proof}
 Importantly, we note that by Lemma~\ref{l:Gsing}, $\ee^{-2 G(z)}/\omega_j(z)$ is bounded in a neighborhood of $[a_j,b_j]$. From the (at most) quarter-root singularities of $\vec G_n(z) \ee^{- \siggg G(z)}$, we find for $z \in \Omega_j'$ that
\begin{align*}
\vec X_n(x,y) = \vec G_n(z)\ee^{- \siggg G(z)} \begin{bmatrix} 0 & 0 \\ - \ee^{-2n\mathfrak g(z) - 2 G(z)}\dbar r_{m,j}^\exop(x,y)/\omega_j(z) & 0 \end{bmatrix}\ee^{\siggg G(z)}\vec G_n(z)^{-1}
\end{align*}
satisfies
\begin{align*}
 \|\vec X_n(x,y) \| \leq C \big|\dbar r_{m,j}^\exop(x,y)\big||z-a_j|^{-1/2} |z - b_j|^{-1/2} \ee^{-2 n \Re \mathfrak g(z)}.
\end{align*}
With the aim of using \eqref{eq:exp_est}, we have
\begin{gather*}
 \|\vec X_n\|_{L^p(\Omega_j')} \leq C \big\|\dbar r_{m,j}^\exop\big\|_{L^q(\Omega_j')} \| h_{j,n}\|_{L^r(\Omega_j')}, \\ h_{j,n}(z) := {|z-a_j|^{-1/2}}|z - b_j|^{-1/2} \ee^{-2 n \Re \mathfrak g(z)},
\end{gather*}
for $1 \leq p \leq 4$, and $1/p = 1/q + 1/r$. Since we will need $p > 2$ below, and using \eqref{eq:exp_est} requires~${\gamma < 3/2}$, we choose $r < 3$, and $q > 6$, giving
\smash{$
 \|\vec X_n\|_{L^p(\Omega_j')} = O\bigl({n^{-1/r}} E_{j,m,q}\bigr)$}.
The first claim follows from the fact that $\|\underline{\mathcal K}_{\Omega'} u\|_{L^\infty}(\mathbb C) \leq C_p \|u\|_{L^p(\Omega')}$ for some $C_p > 0$, provided~${p > 2}$.

To establish the second claim, we need to establish the estimate
\smash{$
 \underline{\mathcal K}_{\Omega'} h_{j,n}(z) = \OO\bigl(n^{-1} \log n\bigr)$}.
This essentially follows from the proof of \cite[Lemma 4]{Yattselev2023} once we note that, after the change of variable $z' = 1/2(\tau + 1/\tau)$, there is a helpful factor of $\big|\tau^2 -1\big|$ that can be used to cancel the singularity in
\begin{align*}
 \big|z^2 - 1\big|^{-1/2} = \frac{2\tau}{\big|\tau^2 -1\big|}.\tag*{\qed}
\end{align*} \renewcommand{\qed}{}
\end{proof}

\subsection{The addition of a finite number of point masses}

The addition of point masses follows a simple procedure, see, for example, \cite{Ding2021}. Diagonal and triangular rational modifications are made to the original RH problem, with the second applying only in a neighborhood of each of the point masses, converting a~residue condition \eqref{eq:def_Y_res} into a~rational jump condition. The effect of this is a rational modification of the densities on $[a_j,b_j]$ for every $j$, a change in the asymptotics at infinity and the introduction of a~jump condition on a~small curve encircling the point masses. Once the $g$-function is introduced, this jump condition becomes exponentially close to the identity. We sketch this procedure here but refer the reader to \cite{Ding2021} for more detail.

Define
\begin{align*}
 \hat {\vec Y}_n(z) = \vec Y_n(z;\mu) \begin{bmatrix} \prod_{j=1}^P (z - c_j)^{-1} & 0 \\
 0 & \prod_{j=1}^P (z - c_j) \end{bmatrix}.
\end{align*}
Then let $\Xi_j$ be a small neighborhood of $c_j$. Define
\begin{align*}
 \check {\vec Y}_n(z) = \hat {\vec Y}_n(z)\begin{cases} \begin{bmatrix} 1 & 0 \\
 -\frac{\tilde {r}_j}{z - c_j} & 1\end{bmatrix}, & z \in \Xi_j,\\
 \vec I, & \text{otherwise}, \end{cases}
\end{align*}
 where $\tilde {r}_j $ is defined as
\[
\tilde {r}_j :=\frac{2 \pi \ii}{{r}_j} \prod_{k \neq j} (c_j - c_k)^{-2}.
\]
The RH problem satisfied by $\check {\vec Y}_n$ is found by modifying the RH problem satisfied by $\vec Y_n$ in three ways:
\begin{itemize}\itemsep=0pt
 \item Replace $\rho(x)$ with
\begin{align}\label{eq:newrho}
 \rho(x) \prod_{j=1}^P (x - c_j)^2.
\end{align}
\item Include the jump
\begin{align*}
 \check {\vec Y}_n^+(z) = \check {\vec Y}_n^-(z) \begin{bmatrix} 1 & 0 \\
 -\frac{\tilde {r}_j}{z - c_j} & 1\end{bmatrix}
\end{align*}
on the positively-oriented boundary of $\Xi_j$.
\item Replace the asymptotics at infinity with
\smash{$
 \check {\vec Y}_n(z) z^{-(n - P) \sigma_3} = \vec I + \OO\bigl(z^{-1}\bigr)$}, $ z \to \infty$.
\end{itemize}

Then Theorem~\ref{t:jacobi} holds if one makes the new definition
\begin{align*}
\vec F^\eexop(x,y) &= \vec H^\eexop(x,y) \begin{bmatrix} \prod_{j=1}^P (z - c_j) & 0 \\
 0 & \prod_{j=1}^P (z - c_j)^{-1} \end{bmatrix}, \\
 \vec H^\eexop(x,y) &=\begin{cases}
 \begin{bmatrix} 1 & 0 \\
 \pm f_j^\exop(x,y)/\omega_j(z) & 1 \end{bmatrix} ,& z \in \Omega_{j,\pm},\vspace{4pt}\\
 \begin{bmatrix} 1 & 0 \\
 \frac{\tilde {r}_j}{z - c_j} & 1\end{bmatrix} ,& z \in \Xi_j,\\
 \vec I, & \text{otherwise},
 \end{cases}
\end{align*}
$\rho$ is replaced with \eqref{eq:newrho} and $n$ is replaced with $n - P$. Similarly, Theorem~\ref{t:cheb} holds if one makes the new definition
\begin{align*}
\vec F^\exop(x,y) &= \vec H^\eexop(x,y) \begin{bmatrix} \prod_{j=1}^P (z - c_j) & 0 \\
 0 & \prod_{j=1}^P (z - c_j)^{-1} \end{bmatrix}, \\
 \vec H^\eexop(x,y) &=\begin{cases}
 \begin{bmatrix} 1 & 0 \\
 \pm f_j^\ddagger(x,y)/\omega_j(z) & 1 \end{bmatrix}, & z \in \Omega'_{j,\pm},\vspace{4pt}\\
 \begin{bmatrix} 1 & 0 \\
 \frac{\tilde {r}_j}{z - c_j} & 1\end{bmatrix}, & z \in \Xi_j,\\
 \vec I, & \text{otherwise},
 \end{cases}
\end{align*}
$\rho$ is replaced with \eqref{eq:newrho} and $n$ is replaced with $n - P$.

\section[Asymptotics of recurrence coefficients and estimating optimal errors]{Asymptotics of recurrence coefficients\\ and estimating optimal errors}\label{s:rec}

While asymptotics of the polynomials themselves are directly available using the above calculations, we are primarily interested in the asymptotics of the associated recurrence coefficients. The recurrence coefficients $(a_n(\mu))_{n \geq 0}$ and $(b_n(\mu))_{n \geq 0}$, $b_n(\mu) > 0$, are the coefficients in the relation
\begin{gather*}
 p_{-1}(x;\mu) := 0,\qquad
 p_{0}(x;\mu) = 1,\\
 x p_n(x;\mu) = b_n(\mu) p_{n+1}(x;\mu) + a_n(\mu) p_n(x;\mu) + b_{n-1}(\mu) p_{n-1}(x;\mu).
\end{gather*}
If we write
\begin{align*}
 \vec Y_n(z;\mu) = \bigl( \vec I + \vec Y_n^{(1)}(\mu) z^{-1} + \OO\bigl(z^{-2}\bigr)\bigr) z^{\siggg n},\qquad \vec Y_n^{(1)}(\mu) = \begin{bmatrix}
 y_{11}(n;\mu) & y_{12}(n;\mu)\\
 y_{21}(n;\mu) & y_{22}(n;\mu)
 \end{bmatrix},
\end{align*}
we have
\begin{align*}
 p_n(z;\mu) &= \ell_n(\mu) \pi_n(z;\mu) = \ell_n(\mu) \left( z^n + y_{11}(n;\mu) z^{n-1} + \cdots \right),\\
 c_n(z;\mu) &= - \frac{1}{2 \pi \ii} \ell_n(\mu)^2 z^{n} + \cdots = y_{12}(n;\mu) z^n + \cdots.
\end{align*}
Equating coefficients in the recurrence, gives
\begin{align*}
 \ell_n(\mu) = b_n(\mu) \ell_{n+1}(\mu),\qquad
 \ell_n(\mu) y_{11}(n;\mu) = b_n(\mu) \ell_{n+1}(\mu) y_{11}(n+1;\mu) + a_n(\mu) \ell_n(\mu).
\end{align*}
And therefore
\begin{align*}
 b_n(\mu) &= \frac{\ell_{n+1}(\mu)}{\ell_n(\mu)}, \qquad
 b_n(\mu)^2 = \frac{y_{12}(n+1;\mu)}{y_{12}(n;\mu)},\qquad
 a_n(\mu) = y_{11}(n;\mu) - y_{11}(n+1;\mu).
\end{align*}

In Theorem~\ref{t:jacobi} or Theorem~\ref{t:cheb}, we have as $z \to \infty$
\begin{align*}
 \vec T_n(x,y) &= \bigl(\vec I + \vec T_n^{(1)} z^{-1} + \OO\bigl(z^{-2}\bigr)\bigr) \vec G_n(z),\qquad
 \vec G_n(z) = \vec I + \vec G_n^{(1)} z^{-1} + \OO\bigl(z^{-2}\bigr)
\end{align*}
for $\vec T_n^{(1)} \to 0$. Then using that
\begin{align*}
\mathfrak g(z) = \log \mathfrak c + \log z + \mathfrak g_1 z^{-1} + \OO\bigl(z^{-2}\bigr),
\end{align*}
we have
\begin{gather*}
 \ee^{n \mathfrak g(z) \siggg} = \mathfrak c^{n\siggg} z^{n\siggg} \bigl(\vec I + n \mathfrak g_1 \siggg z^{-1} + \OO\bigl(z^{-2}\bigr) \bigr),
\\
 \vec Y_n(z;\mu) = \mathfrak c^{-n\siggg} \bigl( \vec I + \bigl[ \vec G_n^{(1)} + \vec T_n^{(1)}\bigr] z^{-1} + \OO\bigl(z^{-2}\bigr) \bigr)\bigl(\vec I + n \mathfrak g_1 \siggg z^{-1} + \OO\bigl(z^{-2}\bigr) \bigr) \mathfrak c^{n\siggg} z^{n\siggg} \\
 \phantom{ \vec Y_n(z;\mu) }{}= \bigl( \vec I + \bigl[ \mathfrak c^{-n\siggg} \bigl[ \vec G_n^{(1)} + \vec T_n^{(1)}\bigr] \mathfrak c^{n\siggg} + n \mathfrak g_1 \sigma_3 \bigr] z^{-1} + \OO\bigl(z^{-2}\bigr) \bigr) z^{n \siggg}.
\end{gather*}
This gives the expansion
\begin{align*}
 \vec Y_n^{(1)}(\mu) = \mathfrak c^{-n\siggg} \bigl[ \vec G_n^{(1)} + \vec T_n^{(1)}\bigr] \mathfrak c^{n\siggg} + n \mathfrak g_1 \sigma_3.
\end{align*}
If we set
\begin{align*}
 \vec G_n^{(1)} = \begin{bmatrix} g_{11}(n;\mu) & g_{12}(n;\mu)\\
 g_{21}(n;\mu) & g_{22}(n;\mu) \end{bmatrix},
\end{align*}
we have
\begin{align*}
 a_n(\mu) &= g_{11}(n;\mu) - g_{11}(n+1;\mu) -\mathfrak g_1 + \OO\bigl(\vec T_n^{(1)}\bigr),\\
 b_n(\mu)^2 &= \frac{g_{12}(n+1;\mu)+ \OO\bigl(\vec T_n^{(1)}\bigr)}{g_{12}(n;\mu)+ \OO\bigl(\vec T_n^{(1)}\bigr)}.
\end{align*}

\begin{figure}[tbp]
 \centering
 \includegraphics[width=0.45\linewidth]{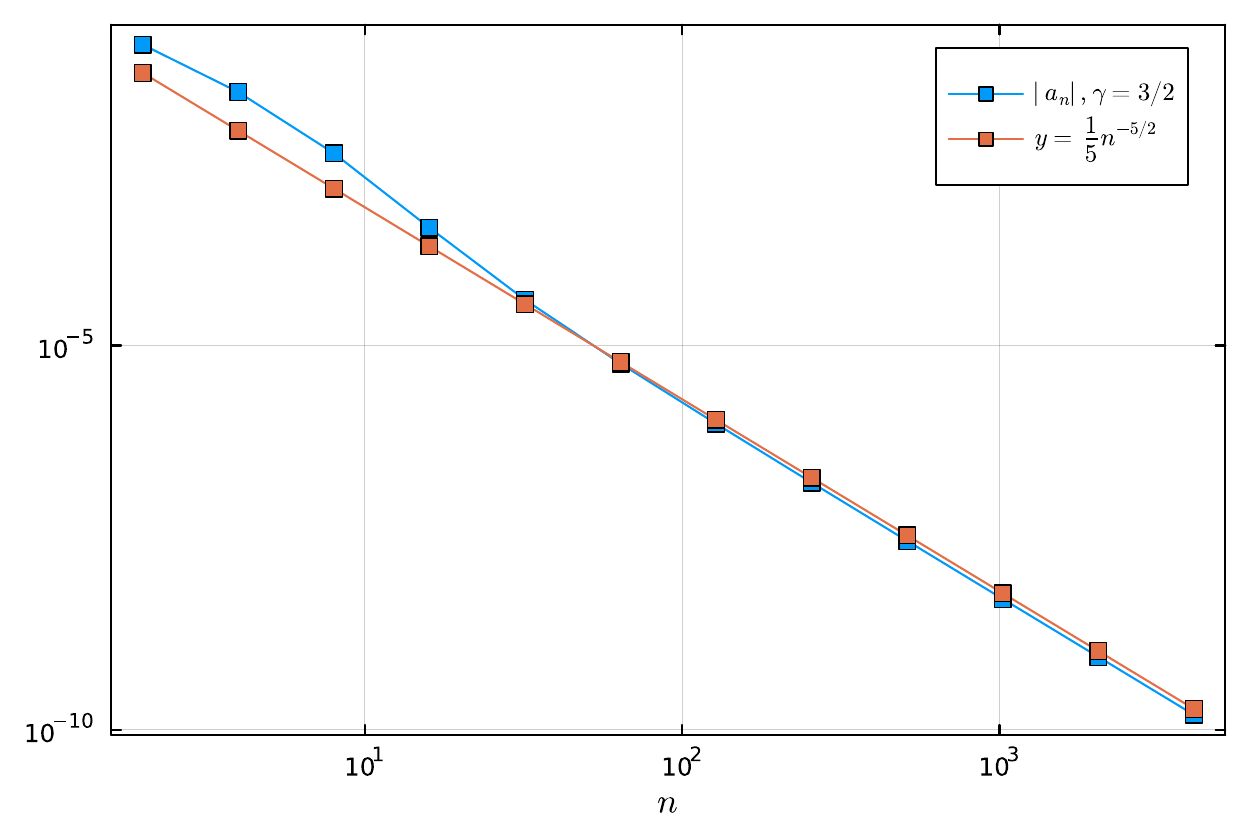} \includegraphics[width=0.45\linewidth]{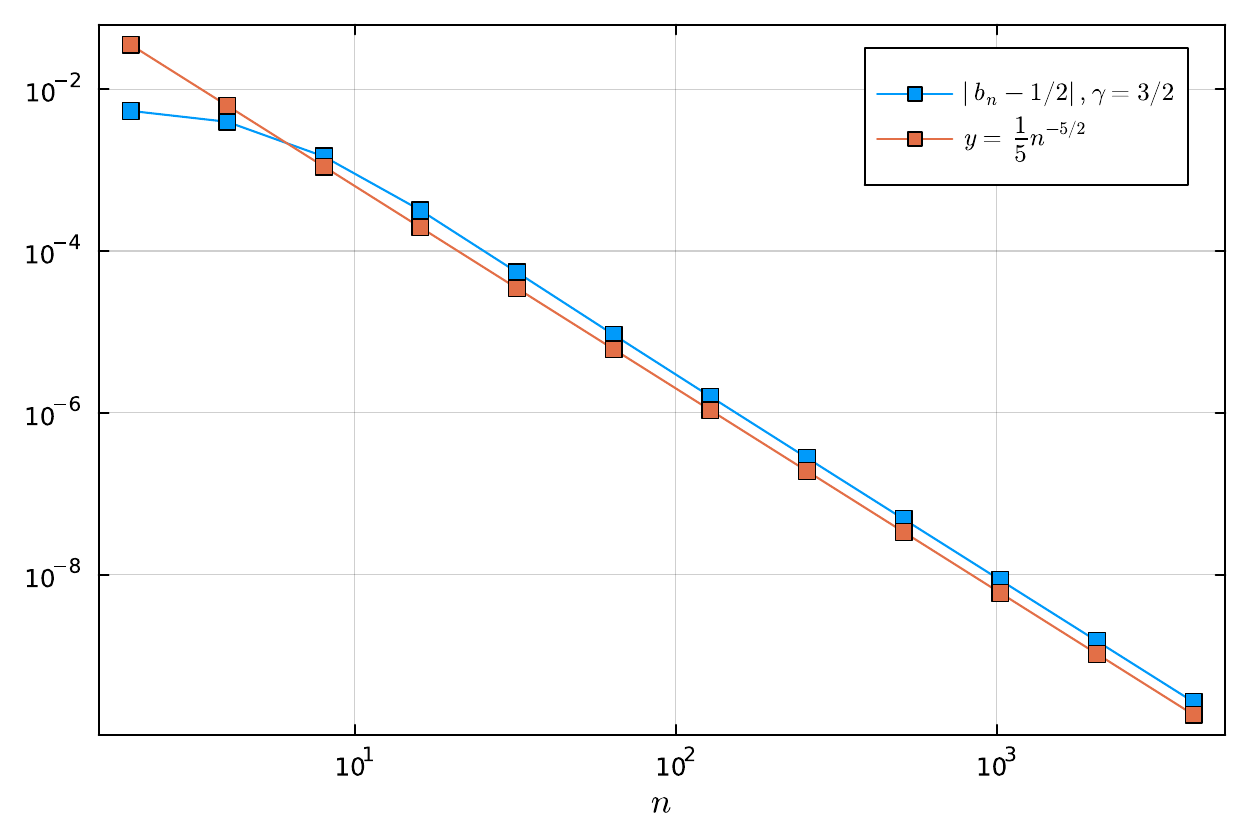}
 \caption{A comparison of the recurrence coefficients when $\gamma = 3/2$. The observed discrepancy between the recurrence coefficients and their known asymptotics appears to be $\OO\bigl(n^{-\gamma-1}\bigr)$.}
 \label{fig:3_2}
\end{figure}

\begin{figure}[tbp]
 \centering
 \includegraphics[width=0.45\linewidth]{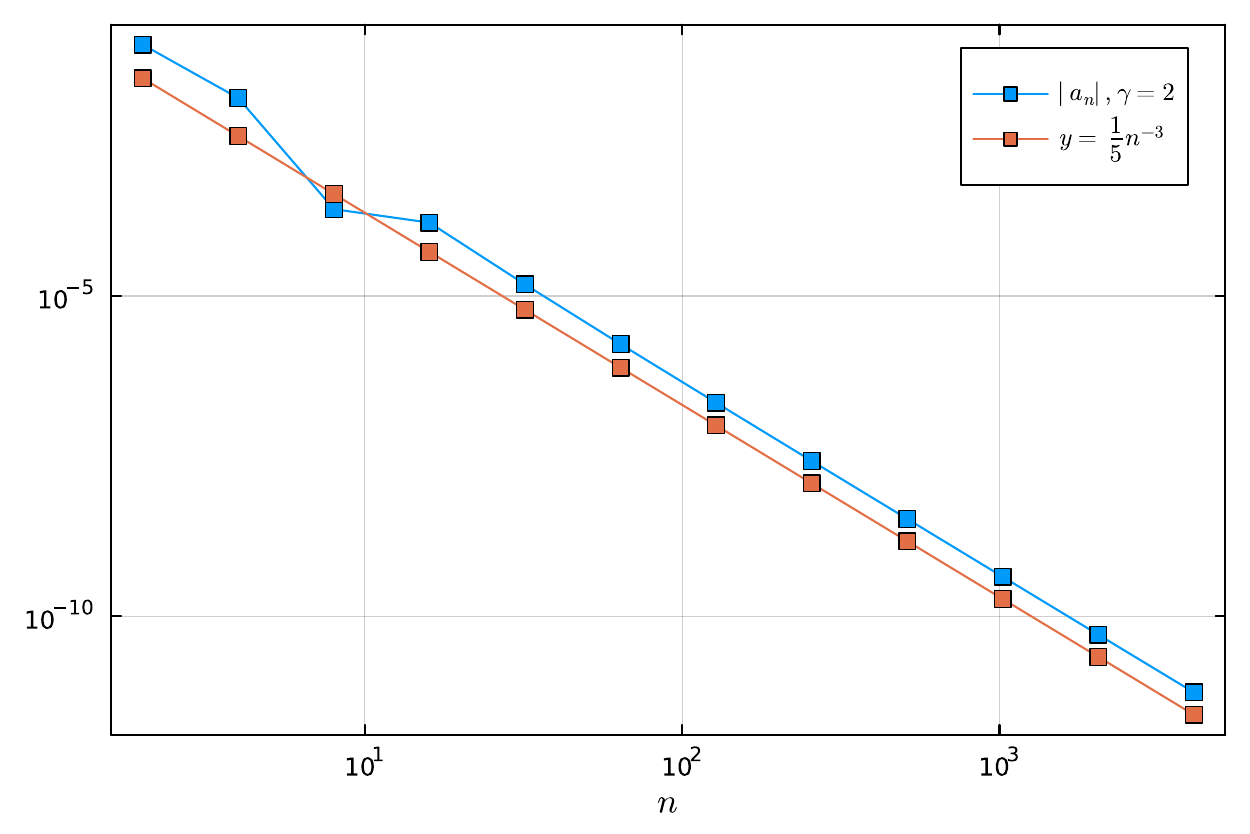} \includegraphics[width=0.45\linewidth]{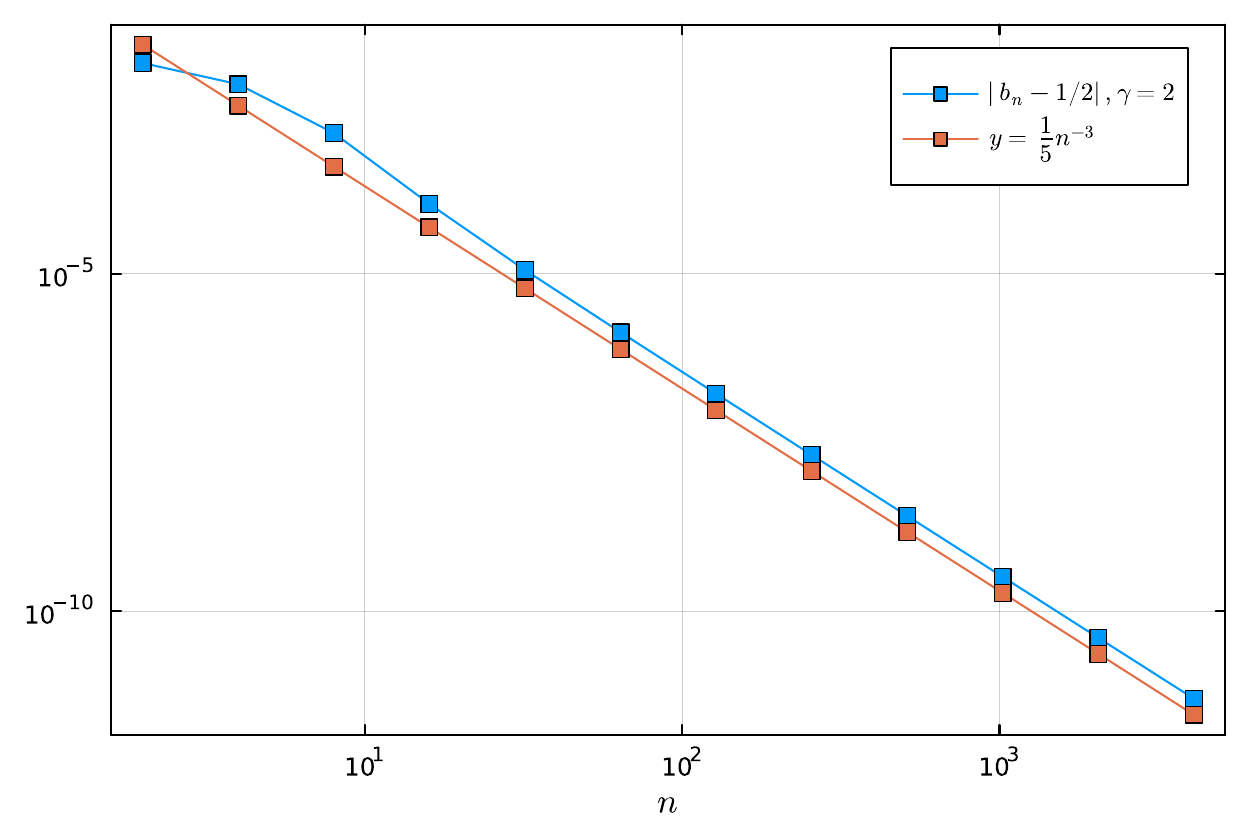}
 \caption{A comparison of the recurrence coefficients when $\gamma = 2$. The observed discrepancy between the recurrence coefficients and their known asymptotics appears to be $\OO\bigl(n^{-\gamma-1}\bigr)$.}
 \label{fig:2}
\end{figure}

To compare our error estimates to the true errors, we consider a single interval $[a_1,b_1] = [-1,1]$, $g = 0$. We use $\alpha_1 = \beta_1 = 1/2$,
\begin{align*}
 h_1(x) = \begin{cases}
 1 , & -1 \leq x \leq 0,\\
 1 + x^{\gamma}(1 + x)^{-1/2}, & 0 \leq x \leq 1.
 \end{cases}
\end{align*}
Then, if one uses the Gauss--Jacobi quadrature rules $(x_j,w_j)$, $(y_j, v_j)$ such that
\begin{gather*}
 \int_{-1}^1 p(x) \sqrt{1-x^2} \dd x = \sum_{j=1}^n p(x_j) w_j,\\
 \int_{0}^1 p(x) x^{\gamma} \sqrt{1 - x} \dd x = \sum_{j=1}^n p(y_j) v_j,
\end{gather*}
for all polynomials $p$ of degree less than $2n$, the inner product
can be approximated well because
\begin{align*}
 \int_{-1}^1 f(x) \overline{g(x)} \sqrt{1 - x^2} h_1(x) \dd x = \sum_{j=1}^n f(x_j) \overline{g(x_j)} w_j + \sum_{j=1}^n f(y_j) \overline{g(x_j)} v_j,
\end{align*}
when $f$ and $g$ are polynomials of degree less than or equal to $n-1$. This implies the recurrence coefficients can be generated exactly (up to roundoff). The ideal algorithm is the RKPW method as outlined in \cite{Gragg1984} which has $O\bigl(n^2\bigr)$ complexity. In Figures~\ref{fig:3_2} and \ref{fig:2}, we provide evidence that error term in the asymptotics of the recurrence coefficients is $O\bigl(n^{-\gamma-1}\bigr)$ whereas our method, and that of Yattselev \cite{Yattselev2023}, using $q = \infty$ predicts $O(n^{-\gamma} \log n)$. We do note that, while potentially indicative, this density is piecewise analytic. A more exotic density might be needed to see the rate predicted by Yattselev, but running computations in such a case will likely be its own entirely separate challenge.

\section{Conclusions and outlook}

The main open problem here is to obtain error bounds that reflect the true rates demonstrated in Figures~\ref{fig:3_2} and \ref{fig:2}. It is also unknown how close to optimal the rate we have obtained for the general Jacobi-type case is.

It is also of interest to see if these deformations, at least in some way, could be helpful in computing the recurrence coefficients as was done in \cite{Ballew2023,Townsend2014,Trogdon2013b}, potentially enabling computations with optimal complexity. We also anticipate using the asymptotics discussed here in perturbation theory \`a la \cite{Ding2021} for random matrices that have anomalous edge/bulk behavior.

\appendix

\section[The g-function]{The $\boldsymbol{ g}$-function}\label{a:gfun}

In this appendix, we define what is commonly called the external Green's function with pole at infinity \cite{Walsh1935}, or rather the analytic function whose real part is this Green's function, and determine some essential properties. This construction can be also be found in~\cite{Akhiezer1961}. Specifically, we find a function $\mathfrak g$ that satisfies
\begin{enumerate}\itemsep=0pt
 \item[(a)] $\mathfrak g'(z) = 1/z + \OO\bigl(1/z^2\bigr)$ as $z \to \infty$,
 \item[(b)] $\mathfrak g_+'(z), \mathfrak g_-'(z) \in \ii \mathbb R$ on $[a_j, b_j]$,
 \item[(c)] $\int_{b_j}^{a_{j+1}} {\mathfrak g'}(z) \dd z = 0$, $j = 1,2,\dots,g$.
\end{enumerate}
 The following construction can be found in \cite{Ding2021}. Set\footnote{See Appendix~\ref{a:RS} for the definition of $R(z)$.}
 \begin{align}\label{eq:gp}
 \mathfrak g'(z) = \frac{Q_{g}(z)}{R(z)}, \qquad \text{where} \quad R(z)^2 = \prod_{j=1}^{g+1} (z- a_j)(z- b_j),
 \end{align}
 where $Q_{g}$ is a monic polynomial of degree $g$. The coefficients for $Q_g$, $Q_{g}(z) = \sum_\ell \alpha_\ell z^\ell$, are found by solving the linear system, which enforces (c)
 \begin{align*}
 \int_{ b_j}^{ a_{j+1}} \sum_{\ell=0}^{g-1} \alpha_k \frac{z^\ell}{R(z)} {\rm d}z = - \int_{ b_j}^{ a_{j+1}} \frac{z^{g}}{R(z)} {\rm d}z, \qquad j = 1,2,\dots,g.
 \end{align*}
This system is uniquely solvable due to abstract theory, see Appendix~\ref{a:RS}.

Then set
\begin{align}\label{eq:g}
 \mathfrak g(z) = \int_{a_1}^z \mathfrak g'(z') \dd z ', \qquad z \not \in [a_1, \infty),
\end{align}
where the path of integration is taken to be a straight line.

 The following properties then hold:
 \begin{enumerate}\itemsep=0pt
 \item[(i)] $\alpha_\ell \in \mathbb R$ for every $\ell$,
 \item[(ii)] $\mathfrak g^+(z) + \mathfrak g^-(z) = 0$, $z \in ( a_j, b_j),$
 \item[(iii)] $\mathfrak g^+(z) - \mathfrak g^-(z) = : \Delta_j \in \ii \mathbb R$ is constant on $[b_{j},a_{j+1}]$, and
 \item[(iv)] $\operatorname{Re} \mathfrak{g}(z)$ is strictly positive on any closed subset of $\mathbb R \setminus \bigcup_j [ a_j, b_j]$.
 \end{enumerate}
We will need to characterize the last property more quantitatively. Consider the function
\begin{align*}
 \phi(z) = z + \sqrt{z- 1} \sqrt{z + 1}.
\end{align*}
It follows that $\phi(z)$ is a conformal mapping from $\mathbb C \setminus [-1,1]$ to the exterior of the unit disk. It can also be shown that
\begin{align*}
 \log |\phi(z)| = \Re \int_{-1}^z \frac{\dd z'}{\sqrt{z'-1} \sqrt{z' + 1}}.
\end{align*}
Then for $a < b$, set $M(z;a,b) = \frac{2}{b-a} \bigl(z - \frac{b + a}{2} \bigr)$ and we establish the following in Appendix~\ref{s:proof-ggrowth}.
\begin{Proposition}\label{p:ggrowth}
For every $\epsilon > 0$ sufficiently small, there exist $\gamma, \gamma' > 0$ such that if $a_j - \epsilon \leq x \leq b_j + \epsilon$ and $|y| \leq \epsilon$, then
 \begin{align*}
\ee^{\Re \mathfrak g(z) } \geq |\phi(M(z;a_j,b_j))|^{\gamma} \qquad \text{and} \qquad \ee^{\Re \mathfrak g(z) } \leq |\phi(M(z;a_j,b_j))|^{\gamma'}.
 \end{align*}
\end{Proposition}

Define the vector
\smash{$
 \boldsymbol{\Delta} = (\Delta_j)_{j=1}^g$}.
We now define conformal maps in the neighborhood of each endpoints $a_j$, $b_j$. Specifically, define for $j = 1,2,\dots,g$,
\begin{align}\label{eq:varphi}
 \varphi(z;a_j) = \left(\int_{a_j}^z \mathfrak g'(z') \dd z'\right)^2, \qquad \varphi(z;b_j) = \left(\int_{b_j}^z \mathfrak g'(z') \dd z'\right)^2.
\end{align}
The following is immediate.
\begin{Lemma}
 For each $j =1,2,\dots,g+1$, $\varphi(z;a_j)$ $($resp. $\varphi(z;b_j))$ is a conformal mapping from a neighborhood of $a_j$ (resp. $b_j$) to a neighborhood of the origin. For $\epsilon$ sufficiently small, and $0 < s < \epsilon$, $\varphi(a_j + s ;a_j) < 0$ and $\varphi(b_j - s ;b_j) < 0$.
\end{Lemma}

\section{Extension operator}\label{a:extend}

In this appendix, we define an extension operator $^\exop$ that will take a sufficiently smooth function defined on the interval $[a,b]$ and extend it to the complex plane suitably. Let $\hat b(x)$ be a smooth function satisfying
\[\label{eq:b}
 \hat b(x) = \begin{cases} 0, & x < 1, \\
 1, & x > 2, \end{cases}
\]
that is monotonic on the interval $[1,2]$. Set
\begin{align}\label{eq:beps}
 b_\tau(x;c) = \hat b(|x-c|/\tau), \qquad 0 < \tau < (b-a)/8.
\end{align}
Then consider
$
 \tilde b(x) = b_\tau(x;a) b_{\tau}(x;b)$.
We use $b$ and $\tilde b$ to define an extension of a function $m\colon [a,b] \to \mathbb R$. We assume $m$ is $\kappa$ times differentiable with its $(\kappa+1)$-th derivative existing almost everywhere. Define the Taylor extension
\begin{align*}
 m(x,y) = \sum_{\ell=0}^\kappa \frac{m^{(\ell)}(x)}{\ell!} (\ii y)^\ell, \qquad x \in [a,b],
\end{align*}
and the Taylor approximations
{\begin{align*}
 m(z; c) = \sum_{\ell = 0}^\kappa \frac{{m}^{(\ell)}(c)}{\ell!} (z-c)^\ell, \qquad c = a, b.
\end{align*}}

For $0 < \theta_1 < \theta_2 < \pi/2$, construct an infinitely smooth function that for $\theta \in [-\pi/2, \pi/2]$ satisfies
\begin{align*}
 \tilde \beta(\theta) = \begin{cases} 0, & - \theta_1\leq \theta \leq \theta_1,\\
 1 , & |\theta| > \theta_2. \end{cases}
\end{align*}
Then extend $\tilde \beta(\theta)$ to the entirety of interval $[-\pi,\pi]$ using an even reflection about $\theta = \pi/2$. Then, one obtains a function on the unit circle in the obvious way
$
 \hat\beta(z) = \tilde \beta(\theta)$, $ z = \ee^{\ii \theta}$,
which then extends to a function on $\mathbb C \setminus \{0\}$
$
 \beta(x,y) = \hat\beta(z/|z|)$.
The extension is then defined as\footnote{Here $m(x,y)$ is only defined for $x \in [a,b]$, but because the $\tilde b(x) = 0$ for $|x- a| < \tau , |x - b| < \tau$, this naturally extends to $|z - c| < \tau$, $c = a,b$.}
\begin{align*}
 m^{\exop}(x,y) ={}& (1 - b_\tau(|y|;0))[\tilde b(x) m(x,y) \\
 &+ (1- b_\tau(x;a)) (1 - \beta(x - a,y)) m(x,y) + \beta(x-a,y) (1 - b_\tau(x;a)) m(z; a)\\
 &+ (1- b_\tau(x;b)) (1 - \beta(x - b,y)) m(x,y) + \beta(x-b,y) (1 - b_\tau(x;b)) m(z; b)] .
\end{align*}
We now verify the following properties of $m^\exop$:
\begin{enumerate}\itemsep=0pt
 \item[(1)] {$m^\exop(x,y)$ extends to be continuous for $a \leq x \leq b$, $y \in \mathbb R$,}
 \item[(2)] for $a \leq x \leq b$, $m^\exop(x,0) = m(x)$,
 \item[(3)] for $c = a, b$, $|z - c| < \tau$,
 \begin{align*}
 m^\exop(x,y) &= m(z;a), \qquad \theta_2 < |\arg(z - a)| < {\pi/2},\\
 m^\exop(x,y) &= m(z;b), \qquad \theta_2 < |\arg(b - z)| < {\pi/2},
 \end{align*}
 \item[(4)] for $a \leq x \leq b$ and $|y| \leq \tau$, $\|\dbar m^\exop(\cdot,y)\|_{L^\infty([a,b])} \leq C_\kappa \|m\|_{C^{\kappa,1}([a,b])} |y|^\kappa$, for some $C_\kappa > 0$.
\end{enumerate}

\subsection{Verification of the properties}

The following lemma will assist in the verification above properties.

\begin{Lemma}\label{l:taylor}
 Suppose $f \in C^{\kappa,1}(I)$, where $I$ is an open interval containing $x = 0$. Then for $x \in I$,
 \begin{align*}
 \left| \sum_{j=0}^\kappa \frac{f^{(j)}(x)}{j!} (\ii y)^j - \sum_{j=0}^\kappa \frac{f^{(j)}(0)}{j!} (x + \ii y)^j\right| \leq C \|f\|_{C^{\kappa,1}(I)} {\sum_{j=1}^{\kappa+1}} |x|^j |y|^{\kappa + 1 -j}
 \end{align*}
 for a constant $C$ that only depends on $I$ and $\kappa$.
\end{Lemma}
\begin{proof}
 We use Taylor's theorem with remainder to write
 \begin{align*}
 f^{(j)}(x) = \sum_{\ell = 0}^{\kappa - j} \frac{f^{(j+\ell)}(0)}{\ell!} x^\ell + R_j(x),
 \end{align*}
 where $|R_j(x)| \leq M_j \sup_{\xi \in I}\big|f^{({\kappa + 1})}(\xi)\big| |x|^{\kappa - j +1}$. Then, the binomial theorem gives
 \begin{align*}
 \sum_{j=0}^\kappa \frac{f^{(j)}(x)}{j!} (\ii y)^j = \sum_{j=0}^\kappa \frac{f^{(j)}(0)}{j!} ( x + \ii y)^j + \sum_{j=0}^\kappa \frac{R_j(x) (\ii y)^j}{j!},
 \end{align*}
 and the claim follows.
\end{proof}

\subsubsection{Property (1)}

The only points where continuity does not follow immediately is $z = a, b$. So, for~${|z - a| < \tau}$ we have
\begin{align*}
 m^\exop(x,y) = (1 - \beta(x - a,y)) m(x,y) + \beta(x-a,y) m(z; a).
\end{align*}
We then write, using Lemma~\ref{l:taylor},
\begin{align*}
 m^\exop(x,y) = m(z;a) + (1 - \beta(x - a,y)) P_{k+1}(x-a,y),
\end{align*}
where $P_{k+1}(x,y)$ is bounded by a homogeneous polynomial of degree $k+1$ in $x$ and $y$. With the definition $m^\exop(a,0) := m(a;a)$, continuity follows. Similar considerations hold near $z = b$.

\subsubsection{Property (2)}

For $a + 2\tau \leq x \leq b - 2\tau$, we have {$b_\tau(x;a) = b_\tau(x;b) = 1$} giving
$
 m^\exop(x,0) = m(x,0) = m(x)$.
Then for $a \leq x \leq a + 2\tau$, {we have $b_\tau(x;b) = 1$ and}
\begin{align*}
 m^\exop(x,0) = {b_\tau(x;a) m(x,0) + (1 - b_\tau(x;a)) m(x,0)}.
\end{align*}
A similar calculation follows near $x = b$.

\subsubsection{Property (3)}

By definition, for $|z - a| < \tau$ and $\theta_2 < |\arg(z - a)| < {\pi/2}$, we find that $m^\exop(x,y) = m(z;a)$. The same calculation follows near $x = b$.

\subsubsection{Property (4)}

We first note that, in general, for almost every $x$
\begin{align*}
 2 \dbar m(x,y) = \sum_{\ell=0}^\kappa \frac{m^{(\ell+1)}(x)}{\ell!} (\ii y)^\ell - \sum_{\ell=1}^\kappa \frac{m^{(\ell)}(x)}{(\ell-1)!} (\ii y)^{\ell-1}
 = \frac{m^{(\kappa+1)}(x)}{\kappa!} (\ii y)^\kappa.
\end{align*}
So, we immediately have the property we desire for $z$ such that $a +2 \tau \leq x \leq b - 2\tau$, when~${\tilde b(x) = 1}$. Then for $|z - a| \leq 2\tau$ and $- \theta_1 < \arg (z - a) < \theta_1$
\begin{align*}
 m^\exop(x,y) = m(x,y) (1- \tilde b(x)) + m(x,y) \tilde b(x).
\end{align*}
We again immediately have the desired property.

Now, consider $0 < |z - a| \leq \tau$ and $|\arg (z - a)| \geq \theta_1$. Then $\tilde b(x) = 0$ and we have
\begin{align*}
 m^\exop(x,y) = (1 - \beta(x-a,y)) m(x,y) + \beta(x - a,y) m(z;a).
\end{align*}
Since $\dbar m(z;a) = 0$, we find
\begin{align*}
 2 \dbar m^\exop(x,y) &= {(2 - 2 \beta(x-a,y))} \dbar m(x,y) + 2 \dbar\beta(x-a,y) \left[m(z;a) - m(x,y) \right].
\end{align*}
The first term admits a bound of the desired type using previous considerations. We use Lemma~\ref{l:taylor} to bound
\begin{align*}
 |m(z;a) - m(x,y)| \leq C \|m\|_{C^{\kappa,1}([a,b])} {\sum_{\ell=1}^{k+2}} |x-a|^{\ell} |y|^{k - \ell + 1}.
\end{align*}
And because $|\arg(z - a)| \geq \theta_1$, we find
\begin{align*}
 |m(z;a) - m(x,y)| \leq C' \|m\|_{C^{\kappa,1}([a,b])} |y|^{k+1}.
\end{align*}
It remains to bound $\dbar \beta$.

Then we compute
\begin{align*}
 \partial_x \beta(x,y) &= \partial_x \hat \beta \left( \frac{x + \ii y}{\sqrt{x^2 + y^2}} \right) = -\hat \beta'(z/|z|) \frac{\ii y}{\bar z |z|},\\
 \partial_y \beta(x,y) &= \partial_y \hat \beta \left( \frac{x + \ii y}{\sqrt{x^2 + y^2}} \right) = \hat \beta'(z/|z|) \frac{\ii x}{\bar z |z|}.
\end{align*}
For some $A > 1$, $A^{-1}|y| \leq |x| \leq A |y|$ whenever $|\arg(z - a)| \geq \theta_1$. This then implies that
\begin{align*}
 A^{-2} y^2 \leq x^2 \leq A^2 y^2,\qquad
 \bigl(A^{-2} + 1\bigr) y^2 \leq x^2 + y^2 \leq \bigl(A^2 + 1\bigr) y^2.
\end{align*}
So, within this sector, we have
\smash{$
 \big| \frac{\ii x}{\bar z |z|} \big| \leq C |y|^{-1}$}, \smash{$
 \big| \frac{\ii y}{\bar z |z|} \big| \leq C |y|^{-1}$},
for a constant $C >0$. Therefore, for a new constant $C'$, we have
\begin{align*}
 |\partial_x \beta(x,y)| \leq C' |y|^{-1},\qquad
 |\partial_y \beta(x,y)| \leq C' |y|^{-1}.
\end{align*}
 {The remaining case of $\tau \leq |z-a| \leq 2 \tau$ and $|\arg (z-a)| \geq \theta_1$ can be established analogously by using $\tilde b(x) = b_\tau(x;a)$ and therefore
\begin{align*}
 m^\sharp(x,y) = m(x,y) + (1 - b_\tau(x;a)) \beta(x-a,y) (m(z;a) - m(x,y)).
\end{align*}
This establishes property (4), with a $\tau$-independent constant for $|y| \leq \tau$.}

\section{Construction and estimation of local parametrices}\label{a:local}

 We first define the classical Bessel parametrix ($\sigma = 2\pi/3$)
\[
\vec P_{\Bes}(\xi;\alpha) = \begin{cases} \begin{bmatrix} \iia\bigl(2 \xi^{1/2}\bigr) & \frac{\ii}{\pi} \Ka\bigl(2\xi^{1/2}\bigr) \\
2 \pi \ii \xi^{1/2} \iia'\bigl(2 \xi^{1/2}\bigr) & -2 \xi^{1/2} \Ka'\bigl(2 \xi^{1/2}\bigr) \end{bmatrix},\qquad \arg \xi \in (- \sigma, \sigma),\vspace{4pt}\\
\begin{bmatrix} \frac{1}{2} \Ho\bigl(2 (-\xi)^{1/2}\bigr) & \frac{1}{2} \Ht \bigl(2 (-\xi)^{1/2}\bigr) \\ \pi \xi^{1/2} {\Ho}'\bigl(2(-\xi)^{1/2}\bigr) & \pi \xi^{1/2} {\Ht}'\bigl(2 (-\xi)^{1/2}\bigr) \end{bmatrix} \ee^{\frac{1}{2} \alpha \pi \ii \siggg},\qquad \arg \xi (\sigma, \pi),\vspace{4pt}\\
\begin{bmatrix} \frac{1}{2} \Ht\bigl(2 (-\xi)^{1/2}\bigr) & - \frac{1}{2} \Ho\bigl(2(-\xi)^{1/2}\bigr) \\ - \pi \xi^{1/2} {\Ht}'\bigl(2(-\xi)^{1/2}\bigr) & \pi \xi^{1/2} {\Ho}' \bigl(2(-\xi)^{1/2}\bigr) \end{bmatrix}\\
\qquad\times\ee^{-\frac{1}{2} \alpha \pi \ii \siggg},\qquad \arg \xi \in (-\pi,-\sigma).
\end{cases}
\]
Here $\iia$, $\Ka$, $\Ht$ and $\Ho$ are the modified Bessel and Hankel functions \cite{DLMF}. From \cite{KuijlaarsInterval}, we have the following.

The function $\vec P_{\Bes}$ satisfies the following jump conditions
\begin{gather*}
\vec P^+_{\Bes}(\xi;\alpha) = \vec P^-_{\Bes}(\xi;\alpha) \vec J_{\Bes}(\xi;\alpha), \qquad\xi \in \Gamma_1 \cup \Gamma_2 \cup \Gamma_3,\\
\vec J_{\Bes}(\xi;\alpha) = \begin{cases} \begin{bmatrix} 1 & 0 \\ \ee^{\alpha \pi \ii} & 1 \end{bmatrix},& \xi \in \Gamma_1,\vspace{4pt}\\\
\begin{bmatrix} 0 & 1 \\ -1 & 0 \end{bmatrix},& \xi \in \Gamma_2,\vspace{4pt}\\\
\begin{bmatrix} 1 & 0 \\ \ee^{-\alpha \pi \ii} & 1 \end{bmatrix},& \xi \in \Gamma_3,\end{cases}
\end{gather*}
where contours are oriented in the direction of {decreasing} modulus. The asymptotics at the origin are given by the following three cases. If $\alpha > 0$,
\begin{align*}
\vec P_{\Bes}(\xi;\alpha) &= \begin{cases} \OO \begin{bmatrix} |\xi|^{\alpha/2} & |\xi|^{-\alpha/2} \\
|\xi|^{\alpha/2} & |\xi|^{-\alpha/2} \end{bmatrix},& |\arg \xi| < 2 \pi /3,\\
\\
\OO \begin{bmatrix} |\xi|^{-\alpha/2} & |\xi|^{-\alpha/2} \\
|\xi|^{-\alpha/2} & |\xi|^{-\alpha/2} \end{bmatrix},& 2 \pi/3 < |\arg \xi| < \pi,
\end{cases} \text{as } \xi \to 0.
\end{align*}
If $\alpha = 0$, then
\begin{align*}
 \vec P_{\Bes}(\xi;\alpha) = \OO \begin{bmatrix} \log |\xi| & \log |\xi|\\
\log |\xi| & \log |\xi| \end{bmatrix},\qquad \text{as } \xi \to 0.
\end{align*}
And if $\alpha < 0$,
\begin{align*}
\vec P_{\Bes}(\xi;\alpha) &= \OO \begin{bmatrix} |\xi|^{\alpha/2} & |\xi|^{\alpha/2} \\
|\xi|^{\alpha/2} & |\xi|^{\alpha/2} \end{bmatrix},\qquad \text{as } \xi \to 0.
\end{align*}
We also remark that $\det \vec P_{\Bes} = 1$ (see \cite[10.28.2]{DLMF}) so that asymptotics for the inverse of~$\vec P_{\Bes}$ can be easily inferred. We also have
\begin{align*}
\vec P_{\Bes}(n^2 \xi;\alpha)={}& \begin{bmatrix} \frac{1}{2} \bigl(\frac{1}{\pi n}\bigr)^{1/2} \xi^{-1/4} \ee^{2 n \xi^{1/2}} & \frac{\ii}{2} \bigl(\frac{1}{\pi n}\bigr)^{1/2} \xi^{-1/4} \ee^{-2 n \xi^{1/2}}\\
\ii (\pi n)^{1/2}\xi^{1/4} \ee^{2 n \xi^{1/2}} & ({\pi n})^{1/2} \xi^{1/4} {\ee^{-2 n \xi^{1/2}}}
\end{bmatrix}\\
&\times\bigl(\vec I + \OO\bigl(n^{-1} |\xi|^{-1/2}\bigr)\bigr), \qquad n \to \infty.
\end{align*}
We rewrite this in a more convenient form
\begin{align}
\vec P_{\Bes}\bigl(n^2 \xi;\alpha\bigr) &= (\pi n)^{-\frac{1}{2} \siggg} \xi^{-\frac{1}{4} \siggg} \vec E_{\Bes}\bigl(n^2 \xi\bigr) \ee^{2 n \xi^{1/2} \siggg},\nonumber\\
\vec E_{\Bes}\bigl(n^2 \xi;\alpha\bigr) & =
 \vec E + \OO\bigl(n^{-1}|\xi|^{-1/2}\bigr), \qquad \vec E = \begin{bmatrix} \frac 1 2 & \frac \ii 2 \\ \ii & 1\end{bmatrix}.\label{eq:bessel-asy}
\end{align}
These asymptotics apply for all $\xi$ with $|\arg \xi| < \pi$. Furthermore, the asymptotics remain valid up to the boundary, $\arg \xi = \pm \pi$.

Now set, following \cite{Kuijlaars2003},
\begin{align*}
 W(z;a_j) &= {\ee^{-\ii \pi (\alpha_j + \beta_j)/2}}(z - b_j)_\ra^{\alpha_j/2} (z - a_j)_\ra^{\beta_j/2} f(z;a_j)^{-1/2},\\
 W(z;b_j) &= (z - b_j)^{\alpha_j/2} (z - a_j)^{\beta_j/2} f(z;b_j)^{-1/2}.
\end{align*}
 The last factor in each line can be seen to be analytic in a neighborhood of $z = a_j, b_j$, respectively, with the choice of the principal branch. For $c = a_j,b_j$, set
$\rho(z;c) = {\omega}_j(z)/f(z;c)$,
and we have for
\begin{gather}
 W^\pm(z;b_j) = \ee^{\pm \ii \pi \alpha_j/2} \rho(z;b_j)^{1/2}, \qquad b_j - \epsilon < z < b_j,\nonumber\\
 W^\pm(z;a_j) = \ee^{\mp \ii \pi \beta_j/2} \rho(z;a_j)^{1/2}, \qquad a_j < z < a_j + \epsilon,\label{eq:Wpm}
\end{gather}
so that for appropriate choices of $z$
\[
W^{+}(z;b_j) W^{-}(z;b_j) = \rho(z;b_j) , \qquad
W^{+}(z;a_j) W^{-}(z;a_j) = \rho(z;a_j).
\]

{\samepage Then for $c = a_j, b_j$, $\gamma = \beta_j, \alpha_j$, respectively, set
 \begin{align*}
 \vec P_n(z; c) = \vec A_n(z;c) \vec Q_{n}(z;c),\qquad
 \vec Q_{n}(z;c) = \vec P_{\Bes}\left(n^2 \frac{\varphi(z;c)}{4}; \gamma\right) W(z;c)^{-\siggg} \ee^{-n \mathfrak g(z) \siggg},
 \end{align*}
where $\vec A_n(z;c)$ is to be determined.}

 We recall that
 \begin{align*}
 \varphi(z;c) = \left(\int_c^z \mathfrak g'(z') \dd z' \right)^2.
 \end{align*}
 And for $\Im z > 0$,
 {\begin{align*}
\mathfrak g(z)& = \int_{a_1}^{b_j} \mathfrak g_+'(z') \dd z' + \int_{b_j}^z \mathfrak g'(z') \dd z' = \sum_{\ell=1}^j \int_{a_\ell}^{b_\ell} \mathfrak g_+'(z') \dd z' + \int_{b_j}^{z} \mathfrak g'(z') \dd z'\\
& = \frac{\Delta_j}{2} + \int_{b_j}^{z} \mathfrak g_+'(z') \dd z'.
 \end{align*}}
Similarly, for $\Im z < 0$, we have
\begin{align*}
\mathfrak g(z) = - \frac{\Delta_j}{2} + \int_{b_j}^{z} \mathfrak g'(z)\dd z'.
\end{align*}
Next, using the principal branch we consider $\varphi(z;b_j)^{1/2}$. As $\varphi(z;b_j)$ is injective for $|z - b_j| < \epsilon$, we conclude that $\varphi(z;b_j)^{1/2}$ can only fail to be analytic for $z \leq b_j$. Because $\varphi(z;b_j)$ is positive for $z > b_j$ the same is true of $\varphi(z;b_j)^{1/2}$. And then we consider
$
\int_{b_j}^z \mathfrak g'(z') \dd z'$, $ z > b_j$.
It can be shown that the monic polynomial $Q_g(z)$ in the \eqref{eq:gp} must have a root in each interval~${(b_j, a_{j+1})}$. From this, it follows that $\mathfrak g'(b_j + \epsilon') > 0$ for $\epsilon' > 0$ sufficiently small. Then
\[
 \varphi(z;b_j)^{1/2} = \int_{b_j}^z \mathfrak g'(z') \dd z'.
\]
And therefore
\begin{align*}
 \varphi(z;b_j)^{1/2} = \begin{cases}
 \mathfrak g(z) - \frac{\Delta_j}{2}, & |z - b_j| < \epsilon,
 \quad \Im z > 0,\\
 \mathfrak g(z) + \frac{\Delta_j}{2}, & |z - b_j| < \epsilon,\quad \Im z < 0.
 \end{cases}
\end{align*}
We repeat these calculations for $\varphi(z;a_j)$ and find that with the convention that $\Delta_0 = 0 = \Delta_{g+1}$
\begin{align*}
 \varphi(z;a_j)^{1/2} = \begin{cases}
 \mathfrak g(z) - \dfrac{\Delta_{j-1}}{2}, & |z - a_j| < \epsilon,\quad \Im z > 0,\\
 \mathfrak g(z) + \dfrac{\Delta_{j-1}}{2}, & |z - a_j| < \epsilon,\quad \Im z < 0.
 \end{cases}
\end{align*}
For $|z - a_j| = \epsilon$, we have $\Im z \neq 0$,
\begin{align*}
 \vec P_n(z;a_j) ={}& \vec A_n(z;a_j)\left(\frac{2}{\pi n}\right)^{\frac{1}{2} \siggg} \varphi(z;a_j)^{-\frac{1}{4} \siggg} \vec E_{\Bes}\left(n^2 \frac{\varphi(z;a_j)}{4}; \beta_j \right)\\
 & \times \ee^{- (\sign \Im z ) \frac{n}{2} \Delta_{j-1} \siggg} W(z;a_j)^{-\siggg}.
\end{align*}
And for $|z - b_j| = \epsilon$, we have $\Im z \neq 0$,
\begin{align*}
\begin{split}
 \vec P_n(z;b_j) ={}& \vec A_n(z;b_j)\left(\frac{2}{\pi n}\right)^{\frac{1}{2} \siggg} \varphi(z;b_j)^{-\frac{1}{4} \siggg} \vec E_{\Bes}\left(n^2 \frac{\varphi(z;b_j)}{4}; \alpha_j \right)\\
 & \times \ee^{- (\sign \Im z ) \frac{n}{2} \Delta_{j} \siggg} W(z;b_j)^{-\siggg}.
 \end{split}
\end{align*}
We set
 \begin{align*}
 \vec A_n(z;a_j) &= \vec G_n(z) W(z;a_j)^{\siggg} \ee^{(\sign \Im z ) \frac{n}{2} \Delta_{j-1} \siggg} \vec E^{-1} \left(\frac{2}{\pi n}\right)^{-\frac{1}{2} \siggg} \varphi(z;a_j)^{\frac{1}{4} \siggg},\\
 \vec A_n(z;b_j) &= \vec G_n(z) W(z;b_j)^{\siggg} \ee^{(\sign \Im z ) \frac{n}{2} \Delta_{j} \siggg} \vec E^{-1} \left(\frac{2}{\pi n}\right)^{-\frac{1}{2} \siggg} \varphi(z;b_j)^{\frac{1}{4} \siggg}.
 \end{align*}

In the following appendices, we discuss the behavior of the parametrix near $b_j$ leaving the calculations near $a_j$ to the reader.

\subsection[The jumps of vec Q\_n near b\_j]{The jumps of $\boldsymbol{\vec Q_n}$ near $\boldsymbol{b_j}$}
We then compute, for example, for $b_j - \epsilon < z < b_j$
\begin{align*}
 \vec Q_n^+(z;b_j) &= \left[\lim_{\epsilon \downarrow 0}\vec P_{\Bes}\left(n^2 \frac{\varphi(z - \ii \epsilon;b_j)}{4}; \alpha_j\right) \right] \begin{bmatrix}
 0 & 1 \\ -1 & 0 \end{bmatrix} \rho(z;b_j)^{-\siggg} W_-(z;b_j)^\siggg \ee^{n \mathfrak g^-(z) \siggg}\\
 & = \vec Q_n^-(z;b_j) \begin{bmatrix} 0 & \rho(z;b_j) \\ -1/ \rho(z;b_j) & 0 \end{bmatrix}.
\end{align*}
Then for $b_j < z < b_j + \epsilon$,
\begin{align*}
 \vec Q_n^+(z;b_j) &= \vec P_{\Bes}\left(n^2 \frac{\varphi(z;b_j)}{4}; \alpha_j\right) W_+(z;b_j)^{-\siggg} \ee^{-n \mathfrak g^+(z) \siggg}\\
 & =\vec P_{\Bes}\left(n^2 \frac{\varphi(z;b_j)}{4}; \alpha_j\right) W_-(z;b_j)^{-\siggg}\ee^{-n \mathfrak g^-(z) \siggg} \ee^{-n \Delta_j \siggg}\\
 & = \vec Q_n^-(z;b_j) \ee^{-n \Delta_j \siggg}.
\end{align*}
For $z \in \Sigma_{b_j,1}$,
\begin{align*}
\vec Q_n^+(z;b_j) = \vec Q_n^-(z;b_j)\begin{bmatrix} 1 & 0 \\ \ee^{\ii \alpha_j \pi}
 W(z;b_j)^{-2} \ee^{-{2}n \mathfrak g(z)} & 1 \end{bmatrix}.
\end{align*}
From \eqref{eq:Wpm}, it follows that
$
 W^+(z;b_j)^2 = \ee^{\ii \alpha_j \pi} \rho(z;b_j)$,
so that $\ee^{\ii \alpha_j \pi} W^+(z;b_j)^{-2}$ can be seen to be the analytic continuation of $1/\rho(z;b_j)$ to the upper-half plane
\begin{align*}
\vec Q_n^+(z;b_j) = \vec Q_n^-(z;b_j)\begin{bmatrix} 1 & 0 \\ \ee^{-{2}n \mathfrak g(z)}f_j^\exop(x,y)/\omega_j(z) & 1 \end{bmatrix}, \qquad z \in \Sigma_{b_j,1}.
\end{align*}
Then because $\ee^{-\ii \alpha_j \pi} W^-(z;b_j)^{-2}$ can be seen to be the analytic continuation of $1/\rho(z;b_j)$ to the lower-half plane
\begin{align*}
\vec Q_n^+(z;b_j) = \vec Q_n^-(z;b_j)\begin{bmatrix} 1 & 0 \\ \ee^{-{2}n \mathfrak g(z)} f_j^\exop(x,y)/\omega_j(z) & 1 \end{bmatrix}, \qquad z \in \Sigma_{b_j,3}.
\end{align*}

 \subsection[The jumps of vec A\_n near b\_j]{The jumps of vec $\boldsymbol{ A_n}$ near $\boldsymbol{b_j}$}
 We then need to analyze the analyticity and jump behavior of $\vec A_n$. We note that it is immediate that $\vec A_n$ is analytic off the real axis. For $b_j < z < b_j + \epsilon$,
 \begin{align*}
 \begin{split}
 \vec A_n^+(z;b_j) & = \vec G_n^+(z) W(z;b_j)^{\siggg} \ee^{\frac{n}{2} \Delta_{j} \siggg} \vec E^{-1} \left(\frac{2}{\pi n}\right)^{-\frac{1}{2} \siggg} \varphi(z;b_j)^{\frac{1}{4} \siggg}\\
 & = \vec G_n^-(z) W(z;b_j)^{\siggg} \ee^{-\frac{n}{2} \Delta_{j} \siggg} \vec E^{-1} \left(\frac{2}{\pi n}\right)^{-\frac{1}{2} \siggg} \varphi(z;b_j)^{\frac{1}{4} \siggg} = \vec A_n^-(z;b_j).
\end{split}
 \end{align*}
 Then for $b_j - \epsilon < z < b_j$ since $\varphi(z;b_j)$ is conformal and real-valued on the real axis,
 \begin{align*}
 \lim_{\epsilon \downarrow 0} \varphi(z \pm \ii \epsilon;b_j)^{1/4} = \ee^{\pm \ii \pi /4} |\varphi(z;b_j)|^{1/4}.
 \end{align*}
 Thus
 \begin{align*}
 \lim_{\epsilon \downarrow 0} \varphi(z + \ii \epsilon;b_j)^{1/4} =: \varphi(z;b_j)_+^{1/4} = \ii \varphi(z;b_j)_-^{1/4} := \lim_{\epsilon \downarrow 0} \varphi(z - \ii \epsilon;b_j)^{1/4}.
 \end{align*}
Then we compute
\begin{align*}
 \vec A_n^+(z;b_j) ={}& \vec G_n^-(z) \begin{bmatrix} 0 & \rho(z) \\ -1/\rho(z) & 0 \end{bmatrix} \rho(z;b_j)^{\siggg} W^-(z;b_j)^{-\siggg} \ee^{\frac{n}{2} \Delta_{j} \siggg} \\
 &\times \vec E^{-1} \left(\frac{2}{\pi n}\right)^{-\frac{1}{2} \siggg} \ii^\siggg \left(\varphi(z;b_j)_-^{\frac{1}{4}}\right)^\siggg.
\end{align*}
To simplify this expression, we first note that for $y \neq 0$,
\begin{align*}
\begin{bmatrix} 0 & y \\ -1/y & 0 \end{bmatrix} = \begin{bmatrix} 0 & 1\\ -1 & 0 \end{bmatrix} y^{-\sigma_3},
\end{align*}
and
\begin{align*}
 \begin{bmatrix} 0 & 1 \\ -1 & 0 \end{bmatrix} y^{\siggg} \vec E^{-1} \left( \frac{\ii}{\sqrt{2}} \right)^\siggg = y^{-\siggg}\begin{bmatrix} 0 & 1 \\ -1 & 0 \end{bmatrix} \begin{bmatrix} 1 & - \frac{\ii}{2} \\
 - \ii & \frac 1 2 \end{bmatrix} \left( \frac{\ii}{\sqrt{2}} \right)^\siggg = y^{-\siggg}\vec E^{-1} 2^{-\siggg/2},
\end{align*}
and therefore
\begin{gather*}
 \vec A_n^+(z;b_j) = \vec G_n^-(z) \left( \frac{\rho(z)}{\rho(z;b_j)} \right)^{\siggg} W^-(z;b_j)^\siggg \ee^{- \frac n 2 \Delta_j \siggg }\vec E^{-1} \left(\frac{2}{\pi n}\right)^{-\frac{1}{2} \siggg} \left(\varphi(z;b_j)_-^{\frac{1}{4}}\right)^\siggg\\
 \phantom{ \vec A_n^+(z;b_j) }{} = \vec A_n^-(z;b_j) \Psi^+(z)^{-\siggg} \vec E \left( \frac{\rho(z)}{\rho(z;b_j)} \right)^{\siggg} \vec E^{-1} \Psi^+(z)^\siggg,\\
 \Psi(z) = \left(\frac{2}{\pi n}\right)^{-\frac{1}{2} \siggg} \varphi(z;b_j)^{\frac{1}{4}\siggg}.
 \end{gather*}

\subsection[The jumps of vec P\_n near b\_j]{The jumps of vec $\boldsymbol{ P_n}$ near $\boldsymbol{ b_j}$}

For all but one contour near $b_j$, the jump of $\vec P_n$ is precisely the same as that of $\vec Q_n$ (and therefore the same as $\vec S_n$) because $\vec A_n$ is analytic. So, it remains to consider $b_j - \epsilon < z < b_j$
\begin{align*}
 \vec P^+_n(z;b_j) ={}& \vec A^+_n(z;b_j)\vec Q^+_n(z;b_j) = \vec A^+_n(z;b_j)\vec Q^-_n(z;b_j) \begin{bmatrix} 0 & \rho(z;b_j) \\ -1/ \rho(z;b_j) & 0 \end{bmatrix}\\
={}& \vec A_n^-(z;b_j)\Psi^+(z)^{-\siggg} \vec E \left( \frac{\rho(z)}{\rho(z;b_j)} \right)^{\siggg} \vec E^{-1} \Psi^+(z)^\siggg\vec Q^-_n(z;b_j) \rho(z;b_j)^{\siggg} \begin{bmatrix} 0 & 1 \\ -1 & 0 \end{bmatrix} \\
={}& \vec P_n^-(z;b_j) \vec Q_n^-(z;b_j)^{-1}\Psi^+(z)^{-\siggg} \vec E \left( \frac{\rho(z)}{\rho(z;b_j)} \right)^{\siggg} \\
 &\times \vec E^{-1} \Psi^+(z)^\siggg\vec Q^-_n(z;b_j)\rho(z;b_j)^{\siggg} \begin{bmatrix} 0 & 1 \\ -1 & 0 \end{bmatrix}.
\end{align*}

\subsection{Local estimates}

We have the following estimates which we state in a series of lemmas.
\begin{Lemma}\label{l:smallz}
 Suppose $|z - b_j| < \epsilon, z \in \Omega_{j,+} \cup \Omega_{j,-}$. If $|\alpha_j| > 0$, then
\begin{align*}
\|\vec Q_n(z;b_j) W(z;b_j)^\siggg \| = \OO\bigl( n^{|\alpha_j|} |z - b_j|^{-|\alpha_j|/2}\bigr) , \qquad n^2 (z-b_j) \to 0.
\end{align*}
If $\alpha_j =0$, then
\begin{align*}
\|\vec Q_n(z;b_j) W(z;b_j)^\siggg \| = \OO( \log n|z - b_j| ) , \qquad n^2 (z-b_j) \to 0.
\end{align*}
Also,
$
 \Re \mathfrak g(z) = \OO\bigl(|z - b_j|^{1/2}\bigr)$, $ z \to b_j$.
\end{Lemma}

\begin{Lemma}\label{l:A}
 The function
 \begin{align*}
 \vec A_n(z;c) \left( \frac{2}{\pi n} \right)^{\frac 1 2 \siggg}
 \end{align*}
 is uniformly bounded in both $z$ and $n$ for $|z - c| < \epsilon$, $c = a_j, b_j$. The same holds for its inverse.
\end{Lemma}
To prove this lemma we use the following lemma\footnote{The proof in \cite{Ding2021} was for $\alpha_j = \beta_j = 1/2$ but the same calculations hold in general.} of \cite[Lemma~A.2]{Ding2021} which shows that $\ee^{G(z) \siggg} W(z;b_j)^\siggg$ is bounded in a neighborhood of $z = b_j$.
 \begin{Lemma}\label{l:Gsing}
 If $k \geq 0$, for some $\epsilon > 0$, and for every $j =1,2,\dots,g+1$, we have
 \begin{align*}
 G(z) = -\frac {\alpha_j}{2} \log(z - b_j) -\frac {\beta_j}{2} \log(a_j -z) + R_{j}(z), \qquad \mathrm{dist}(z,[a_j, b_j]) \leq \epsilon,
 \end{align*}
 where $R_j(z)$ is a uniformly bounded function for $\mathrm{dist}(z,[ a_j, b_j]) \leq \epsilon$.
 \end{Lemma}
\begin{proof}[Proof of Lemma~\ref{l:A}]
By Lemma~\ref{l:Gsing}, it follows that $\ee^{G(z) \siggg} W(z;b_j)^\siggg$ is bounded in a neighborhood of $z = b_j$. And then one can show, as is classically used the derivation of the local solutions \cite{KuijlaarsInterval}, that if $\rho(z)$ is replaced with $\rho(z;b_j)$ in the definition of $G(z)$ for, say, $b_j - 2 \epsilon < z < b_j$, then $\vec A_n(z;b_j)$ is analytic for $|z - b_j| < \epsilon$ as the singular factors produce weaker-than-pole singularities. Thus, the actual $\vec A_n(z;b_j)$ is a bounded perturbation of this analytic approximation. The theta function terms that contribute to $ \vec Q_n$ are uniformly bounded, and the claim follows.
\end{proof}

\begin{Lemma}\label{l:midz}
 Suppose $|z - b_j| < \epsilon$, $z \in \Omega_{j,+} \cup \Omega_{j,-}$. For any $c < C = C(n)$, $C \ll n^2$, for sufficiently large $n$, there exists a constant $D$ such that if $c \leq \big|n^2 (b_j - z)\big| \leq C$, then
 \begin{align*}
 \|\Psi(z)^\siggg\| \leq D C^{1/4}, \qquad
 \|\vec A_n(z;b_j)\| &\leq D C^{1/4}, \qquad
 \|\vec Q_n(z;b_j) W(z;b_j)^{\siggg} \| \leq D C^{1/2}.
 \end{align*}
\end{Lemma}
\begin{proof}
 We have for $|z - b_j| \leq \epsilon$ and some $d \in \mathbb R$, ${D'} > 0$,
 \begin{gather*}
 |\varphi(z;b_j) - d (z-b_j)| \leq {D'} |z - b_j|^2, \\
 |n^2\varphi(z;b_j)| - \big|d n^2 (z-b_j)\big| \leq (n^4{D'} |z - b_j|^2) n^{-2} \leq {D'}C^2/n^2.
 \end{gather*}
 So
$
 |d|/C - {D'}C^2/n^2 \leq |n^2\varphi(z;b_j)| \leq |d| C + {D'}C^2/n^2$,
 and the bound on $\Psi(z)^\siggg$ follows for~$n$ sufficiently large. It follows from previous arguments that~${\vec A_n(z;b_j) \Psi(z)^{-\siggg}}$ is bounded for~${|z - b_j| < \epsilon}$ as $n \to \infty$. The claim then follows for $\vec A_n(z;b_j)$. Then, for $\vec Q_n(z;b_j) W(z;b_j)^\siggg $ factors of $W(z;b_j)^\siggg$ cancel. Then if $C = O(1)$, we have from Lemma~\ref{l:smallz}
$
 \ee^{n \mathfrak g(z) } = \OO(1)$,
 and~\smash{$\vec P_{\Bes}\bigl( n^2 \frac{\varphi(z;b_j)}{4}\bigr) = \OO(1)$} uniformly for $z$ in the range. Now, let $C' > 0$ be sufficiently large, but an absolute constant, so that the entry-wise bound from \eqref{eq:bessel-asy} holds for $|n^2 \xi| \geq C'$. And then for $n^2|z- b_j| \geq 2 C'$, there exists $n$ sufficiently large so that \smash{$\big| n^2 \frac{\varphi(z;b_j)}{4} \big| \geq C'$}, and for~${C' \leq \big|n^2 (z -b_j)\big| \leq C}$, we have
 \begin{align*}
 \|\Psi(z)^\siggg\vec Q_n(z;b_j) W(z;b_j)^{\siggg} \| &\leq D.
 \end{align*}
 The last claim follows from the estimates on $\Psi$.
\end{proof}

\begin{Lemma}\label{l:largez}
 Suppose $|z - b_j| < \epsilon$, $z \in \Omega_{j,+} \cup \Omega_{j,-}$. Then as $\big|n^2 (b_j - z)\big| \to \infty$
 \begin{align*}
&\|\vec A_n(z;b_j) \vec Q_n(z;b_j) W(z;b_j)^{\siggg}\| = \OO\bigl(|z - b_j|^{-1/4}\bigr),\\
&\big\| W(z;b_j)^{-\siggg} \vec Q_n(z;b_j)^{-1} \Psi(z)^{-\siggg} \vec E - \ee^{(\sign \Im z)\frac n 2 \Delta_{j} \siggg}\big\| = \OO\bigl(n^{-1} |z-b_j|^{-1/2}\bigr),\\
&\big\| \vec E^{-1} \Psi(z)^\siggg \vec A_n(z;b_j)^{-1} \big\| = \OO\bigl(|z-b_j|^{-1/4}\bigr).
\end{align*}
\end{Lemma}
\begin{proof}
 For the first claim, we can write
 \begin{align*}
 \vec A_n(z;b_j) \vec Q_n(z;b_j) W(z;b_j)^{\siggg} ={}& \vec G_n(z) W(z;b_j)^{\siggg} \ee^{(\sign \Im z)\frac{n}{2} \Delta_{j} \siggg} \vec E^{-1} \vec E_{\Bes}\\
 &\times\left(n^2 \frac{\varphi(z;b_j)}{4}; \alpha_j \right)\ee^{-(\sign \Im z)\frac{n}{2} \Delta_{j} \siggg}.
 \end{align*}
 The claim then follows from the asymptotics of $\vec E_{\Bes}$ and the fact that, from Lemma~\ref{l:Gsing}, $\vec G_n(z) W(z;b_j)^{\siggg}$ has quarter-root singularities.

 Then we write
 \begin{align*}
\vec E^{-1}\Psi(z)^{\siggg} \vec Q_n(z;b_j)W(z;b_j)^{\siggg}\ee^{(\sign \Im z) \frac n 2 \Delta_{j} \siggg}= \vec I + \OO\bigl(n^{-1} |z-b_j|^{-1/2}\bigr),
 \end{align*}
 from which the second claim follows. The last claim follows analogously.
\end{proof}

\subsection{Proof of Proposition~\ref{p:ggrowth}}\label{s:proof-ggrowth}

 We first consider $y > 0$ and $a_j \leq x \leq b_j$. We write
 \begin{align*}
 \mathfrak g'(z) = \frac{1}{\sqrt{z - a_j} \sqrt{z - b_j}} r(z),
 \end{align*}
 where $r(z)$ is analytic in a neighborhood of $[a_j,b_j]$ and is positive on $[a_j,b_j]$. We have
 \begin{align*}
 \Re \mathfrak g(z) = {-}\int_{0}^y \Im \mathfrak g'(x + \ii y') \dd y'.
 \end{align*}
 Similarly, {using $M$ in \eqref{eq:has_M}},
 \begin{align*}
 \log |{\phi}(M(z;a_j,b_j))| = {-}\int_{0}^y \Im \frac{\dd y'}{\sqrt{x + \ii y'-b_j} \sqrt{x +\ii y' -a_j}}.
 \end{align*}
 We write $r(z) = h(x,y) + \ii g(x,y)$, where $h$ and $g$ are real valued. From Taylor's theorem,
 \begin{align*}
 h(x + \ii y) = h(x) + \partial_y h (x,\xi_1) y,\qquad
 g(x + \ii y) = \partial_y g (x,\xi_2) y,
 \end{align*}
 for $0 < \xi_1, \xi_2 < y$, or, rather that
$
 |h(x + \ii y) - h(x)| \leq C |y|$, $
 |g(x + \ii y)| \leq C |y|$.
 Then we have that
 \begin{align*}
 \begin{split}
 \Im \mathfrak g'(x + \ii y) ={}& h(x + \ii y) \Im \frac{1}{\sqrt{x + \ii y-b_j} \sqrt{x +\ii y -a_j}} \\
 &+ g(x + \ii y) \Re \frac{1}{\sqrt{x + \ii y-b_j} \sqrt{x +\ii y -a_j}}.
 \end{split}
 \end{align*}
 We then claim that
 \begin{align*}
 \left |\Re \frac{\sqrt{y}}{\sqrt{x + \ii y-b_j} \sqrt{x +\ii y -a_j}} \right| \leq C'.
 \end{align*}
 Indeed,
 \begin{align*}
 \frac{\sqrt{y}}{|x + \ii y-b_j
|^{1/2} |x +\ii y -a_j|^{1/2}} & = \frac{1}{\bigl( (x- b_j)^2/y^2 + 1\bigr)^{1/4} \bigl( (x- a_j)^2 + y^2\bigr)^{1/4}}\\
\leq \frac{1}{\bigl( (x- a_j)^2 + y^2\bigr)^{1/4}}.
\end{align*}
And we find that
\begin{align*}
 \left |\Re \frac{\sqrt{y}}{\sqrt{x + \ii y-b_j} \sqrt{x +\ii y -a_j}} \right| &\leq \min\left\{\frac{1}{\bigl( (x- a_j)^2 + y^2\bigr)^{1/4}}, \frac{1}{\bigl( (x- b_j)^2 + y^2\bigr)^{1/4}}\right\} \\ &\leq \frac{2}{(b_j - a_j)^{1/2}}.
\end{align*}
This results in, for {$c' \geq h(x + \ii y) \geq c$,
\begin{gather*}
 c' \Im \frac{1}{\sqrt{x + \ii y-b_j} \sqrt{x +\ii y -a_j}} + C''\sqrt{y} \geq \Im \mathfrak g'(x + \ii y) {\dd y}\\
\phantom{ c' \Im \frac{1}{\sqrt{x + \ii y-b_j} \sqrt{x +\ii y -a_j}} + C''\sqrt{y} }{} \geq c \Im \frac{1}{\sqrt{x + \ii y-b_j} \sqrt{x +\ii y -a_j}} - C''\sqrt{y},\\
 c' \log |{\phi}(M(z;a_j,b_j))| + \frac{2C''}{3} y^{3/2} \geq \int_0^y \Im \mathfrak g'(x + \ii y) \dd y \\
 \phantom{ c' \log |{\phi}(M(z;a_j,b_j))| + \frac{2C''}{3} y^{3/2} }{}\geq c \log |{\phi}(M(z;a_j,b_j))| - \frac{2C''}{3} y^{3/2}.
\end{gather*}}
We then claim for $0 < y < \gamma$ and $a_j \leq x \leq b_j$, that there exists ${c''} > 0$ such that
\begin{align*}
 I(x,y) = \Im \frac{1}{\sqrt{x + \ii y-b_j} \sqrt{x +\ii y -a_j}} \geq {c''}.
\end{align*}
Indeed, if this were not the case, then there would exist a convergent sequence $(x_n,y_n)_{n \geq 1}$ such that $I(x_n,y_n) \to 0$. And we must have that $(x_n,y_n) \to (a_j,0)$ or $(x_n,y_n) \to(b_j,0)$. But neither can occur as $I(x_n,y_n) \to +\infty$. This implies
\begin{align*}
 \log |{\phi}(M(z;a_j,b_j))| \geq {c''} y,\qquad
 -y^{3/2} \geq -\frac{\sqrt{y}}{{c''}} \log |{\phi}(M(z;a_j,b_j))|.
\end{align*}
And therefore
\begin{align*}
 \left( c' +\sqrt{y}\frac{2C''}{3{c''}} \right) \log |{\phi}(M(z;a_j,b_j))| &\geq \int_0^y \Im \mathfrak g'(x + \ii y) {\dd y}\\ &\geq \left( c' -\sqrt{y}\frac{2C''}{3{c''}} \right) \log |{\phi}(M(z;a_j,b_j))|,
\end{align*}
from which the claim follows. Similar considerations work for $\Im y < 0$, establishing the claim for $a_j \leq x \leq b_j$. Continuity of the ratio $\Re \mathfrak g(z)/\log |{\phi}(M(z;a_j,b_j))|$ near $a_j$, $b_j$ allows the desired bound to hold in neighborhoods of $a_j$, $b_j$, completing the proof.

\section{Hyperelliptic Riemann surfaces}\label{a:RS}

\begin{figure}[th]\centering
\begin{tikzpicture}[scale=1]
\coordinate (a1) at (0,0);
\coordinate (b1) at (1,0);
\coordinate (a2) at (2.5,0);
\coordinate (b2) at (3.5,0);


 \draw[smooth, line cap = round, line width=1pt] plot[tension=0.65] coordinates{(6.5,-0.85) (6,-0.75) (4.5,-1.5) (3,-0.75) (1.5,-1.5) (0,0) (1.5,1.5) (3, 0.75) (4.5, 1.5) (6,0.75) (6.5,0.85)};

 \draw[smooth, line cap = round, line width=1pt] plot[tension=0.65] coordinates{(7.5,-0.85) (8,-0.75) (9.5,-1.5) (11,0) (9.5, 1.5) (8, 0.75) (7.5,0.85)};

 \node at (7,0) {$\cdots$};

\draw[smooth, line cap = round, line width=1pt] plot[tension=0.65] coordinates{(1,0) (1.5,0.25) (2,0)};
\draw[smooth, line cap = round, line width=1pt] plot[tension=0.65] coordinates{(0.875, 0.25) (1,0) (1.5, -0.25) (2,0) (2.125,0.25)};

\draw[smooth, line cap = round, line width=1pt] plot[tension=0.65] coordinates{(4,0) (4.5,0.25) (5,0)};
\draw[smooth, line cap = round, line width=1pt] plot[tension=0.65] coordinates{(3.875, 0.25) (4,0) (4.5, -0.25) (5,0) (5.125,0.25)};

\draw[smooth, line cap = round, line width=1pt] plot[tension=0.65] coordinates{(9,0) (9.5,0.25) (10,0)};
\draw[smooth, line cap = round, line width=1pt] plot[tension=0.65] coordinates{(8.875, 0.25) (9,0) (9.5, -0.25) (10,0) (10.125,0.25)};

\node[left] at (0,0) {\small $ a_1$};
\node at (0,0) {\color{NicePurple}\textbullet};
\node[left] at (1,0) {\small $ b_1$};
\node at (1,0) {\color{NicePurple}\textbullet};
\node[right] at (2,0) {\small $ a_2$};
\node at (2,0) {\color{NicePurple}\textbullet};

\node[left] at (4,0) {\small $ b_2$};
\node at (4,0) {\color{NicePurple}\textbullet};
\node[right] at (5,0) {\small $ a_3$};
\node at (5,0) {\color{NicePurple}\textbullet};

\node[left] at (9,0) {\small $ b_g$};
\node at (9,0) {\color{NicePurple}\textbullet};
\node[right] at (9.96,0) {\small $ a_{g+1}$};
\node at (10,0) {\color{NicePurple}\textbullet};

\node[right] at (11,0) {\small $ b_{g+1}$};
\node at (11,0) {\color{NicePurple}\textbullet};


\begin{scope}[decoration={
 	 markings,
	 mark=at position 0.56 with {\arrow[line width =1.2pt]{>}}
	 }
]
\draw[smooth cycle, line cap = round, line width=1pt, postaction=decorate, MidnightBlue] plot[tension=0.65] coordinates{(0.2,0) (1.5,-0.75) (2.8,0) (1.5, 0.75)};
\node[above, yshift=-2pt] at (1.5,0.75) {\color{MidnightBlue}\small ${\mathfrak{a}}_1$};

\draw[smooth cycle, line cap = round, line width=1pt, postaction=decorate, MidnightBlue] plot[tension=0.65] coordinates{(3+0.2,0) (3+1.5,-0.75) (3+2.8,0) (3+1.5, 0.75)};
\node[above, yshift=-2pt] at (3+1.5,0.75) {\color{MidnightBlue}\small ${\mathfrak{a}}_2$};

\draw[smooth cycle, line cap = round, line width=1pt, postaction=decorate, MidnightBlue] plot[tension=0.65] coordinates{(8+0.2,0) (8+1.5,-0.75) (8+2.8,0) (8+1.5, 0.75)};
\node[above, yshift=-2pt] at (8+1.5,0.75) {\color{MidnightBlue}\small ${\mathfrak{a}}_g$};

\end{scope}


\begin{scope}[decoration={
 	 markings,
	 mark=at position 0.56 with {\arrow[line width =1.2pt]{<}}
	 }
]
\draw[smooth, line cap = round, line width=1pt, postaction=decorate] plot[tension=0.65] coordinates{(1.5,-0.25) (1.25,-0.95) (1.5, -1.5)};
\draw[smooth, line cap = round, line width=1pt, postaction=decorate, dashed] plot[tension=0.65] coordinates{(1.5, -1.5) (1.65,-0.95) (1.5,-0.25)};
\node[xshift=-5pt, yshift=-2pt] at (1.25,-0.95) {\small ${\mathfrak{b}}_1$};

\draw[smooth, line cap = round, line width=1pt, postaction=decorate] plot[tension=0.65] coordinates{(3+1.5,-0.25) (3+1.25,-0.95) (3+1.5, -1.5)};
\draw[smooth, line cap = round, line width=1pt, postaction=decorate, dashed] plot[tension=0.65] coordinates{(3+1.5, -1.5) (3+1.65,-0.95) (3+1.5,-0.25)};
\node[xshift=-5pt, yshift=-2pt] at (3+1.25,-0.95) {\small ${\mathfrak{b}}_2$};

\draw[smooth, line cap = round, line width=1pt, postaction=decorate] plot[tension=0.65] coordinates{(8+1.5,-0.25) (8+1.25,-0.95) (8+1.5, -1.5)};
\draw[smooth, line cap = round, line width=1pt, postaction=decorate, dashed] plot[tension=0.65] coordinates{(8+1.5, -1.5) (8+1.65,-0.95) (8+1.5,-0.25)};
\node[xshift=-5pt, yshift=-2pt] at (8+1.25,-0.95) {\small ${\mathfrak{b}}_g$};

 \end{scope}

\end{tikzpicture}
\caption{An illustration of the Riemann surface $\Gamma$. }\label{fig_riemman}
\end{figure}
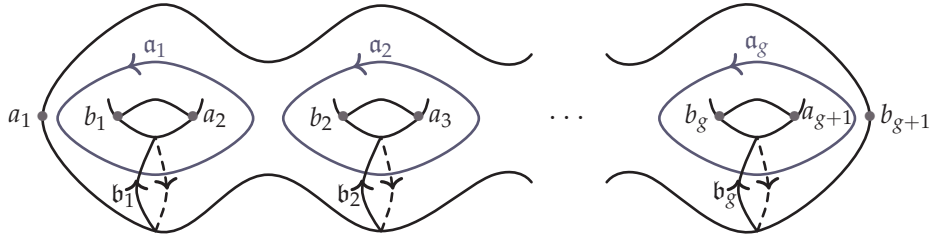

Classic references for what follows are \cite{Baik2007,Algebro, Deift1999}. Much of what follows here was closely adapted from that in \cite{Ding2021b}. Using the intervals $[ a_j, b_j]$, $1 \leq j \leq g+1$, we define a Riemann surface by the zero locus of
\begin{align*}
 w^2 - \prod_{j=1}^{g+1} (z - a_j)(z - b_j) =: w^2- P_{2g+2}(z)
\end{align*}
in $\mathbb C^2$. Then define the analytic function
\begin{align*}
 &R\colon\ \hat {\mathbb C} \to \mathbb C, \qquad R(z)^2 = P_{2g+2}(z), \qquad R(z)z^{-g-1} \to 1, \qquad \text{as} \quad z \to \infty,
\end{align*}
where \begin{align*}
 \hat {\mathbb C} = \mathbb C \setminus \bigcup_{j=1}^{g+1} [ a_j, b_j].
\end{align*}
A Riemann surface $\Gamma$ is built by joining two copies of $\hat {\mathbb C}$; see Figure~\ref{fig_riemman}. The Riemann surface has a canonical set of cycles, a homology basis. For a cartoon of these $\mathfrak a$-cycles and $\mathfrak b$-cycles see Figure~\ref{fig_riemman}. We have a natural projection operation $\pi\colon \Gamma \to \mathbb C$ defined by $\pi((z,w)) = z$ and its right-inverses $\pi_j^{-1}(z) = \bigl(z, (-1)^{j+1} R(z)\bigr)$, $j = 1,2$.

It is well-known that (see \cite{Algebro}, for example) a basis for holomorphic differentials on $\Gamma$ is given~by
\begin{align*}
 \dd \nu_j = \frac{z^{j-1}}{R(z)} \dd z, \qquad j = 1,2,\dots,g+1.
\end{align*}
Then define the $g\times g$ matrix of periods $\vec A = (A_{ij})$ by
\smash{$
 A_{ij} = \oint_{{\mathfrak a}_i} \dd \nu_j$}.
Note that if $\vec c = \smash{\begin{bmatrix} c_1 & c_2 & \cdots & c_g \end{bmatrix}^{\mathsf{T}}} = {\vec A}^{-1} \vec e_j$ for the standard basis vector $\vec e_j$, then
\begin{align*}
 \oint_{{\mathfrak a}_i} \sum_{k=1}^g c_k \dd \nu_k = \sum_{k=1}^g c_k A_{ik} = \vec e_i^{\mathsf{T}} \vec A \vec c = \vec e_i^{\mathsf{T}} \vec e_j = \delta_{ij}.
\end{align*}
A basis of normalized differentials is then given by
\begin{align*}
 \begin{bmatrix} \dd \omega_1 \\ \dd \omega_2 \\ \vdots \\ \dd \omega_g \end{bmatrix} = 2 \pi \ii{\vec A}^{-1} \begin{bmatrix} \dd \nu_1 \\ \dd \nu_2 \\ \vdots \\ \dd \nu_g \end{bmatrix}.
\end{align*}
These satisfy
$
 \oint_{{\mathfrak a}_i} \dd \omega_j = 2 \pi \ii \delta_{ij}$.
The fact that the matrix $\vec A$ is invertible follows from abstract theory \cite{Algebro}.

Define
\begin{align*}
 \vec u(z) = \left( \int_{a_1}^z \dd \omega_j \right)_{j=1}^g, \qquad z \not \in \mathbb R,
\end{align*}
where, for the sake of concreteness, the path of integration is taken to be a straight line connecting $a_1$ to $z$. We can treat $\vec u$ as both a function $\vec u(z)$ on $\hat {\mathbb C}$ and as a function $\vec u(P)$ on $\Gamma$. Then~$\vec u(P)$ is a single-valued analytic function on the Riemann surface $ \Gamma$ provided $\Gamma$ is cut along the cycles $\{{\mathfrak a}_1, \dots, {\mathfrak a}_g, {\mathfrak b}_1, \dots, {\mathfrak b}_g\}$, making it simply connected. Another important fact is for~${z \in \hat{\mathbb C}}$, $\vec u\bigl(\pi_1^{-1}(z)\bigr) = - \vec u\bigl(\pi_2^{-1}(z)\bigr)$.

The associated Riemann matrix of ${\mathfrak b}$ periods is given by
\[
\btau=\left( \tau_{i j}\right)=\left(\int_{{\mathfrak b}_{j}} \dd \omega_{i}\right)_{1 \leq i, j \leq g}.
\]
It follows, see \cite{Algebro}, that $\btau$ is symmetric, real and negative definite. The vector $\vec k$ of Riemann constants is defined componentwise via
\begin{align*}
 k_j = \frac{2 \pi \ii + \tau_{jj}}{2} - \frac{1}{2 \pi \ii} \sum_{\ell \neq j} \oint_{{\mathfrak a}_\ell} u_j \dd \omega_\ell, \qquad j =1,2,\dots,g.
\end{align*}
The associated theta function is defined by
\[
\theta(\vec z;\btau)=\sum_{m \in \mathbb{Z}^{g}} \exp \left( \frac 1 2 (m,\btau m) + (m,\vec z) \right), \qquad \vec z \in \mathbb{C}^{g},
\]
where $(\vec a, \vec b) = \vec a^{\mathsf{T}} \vec b$. We have
\begin{align*}
 \theta(\vec z + 2\pi \ii \vec e_j;\btau) &= \theta(\vec z;\btau),\qquad
 \theta(\vec z + \btau \vec e_j; \btau) = \exp \left( - \frac 1 2 \tau_{jj} - z_j\right) \theta(\vec z;\btau).
\end{align*}

A divisor $D = \sum_j n_j P_j$ is a formal sum of points $\{P_j\}$ on the Riemann surface $\Gamma$. The Abel map of a divisor is defined via
$
 \mathcal A (D) = \sum_j n_j \vec u(P_j)$.

 \newcommand{\bTheta}{\boldsymbol{\Theta}}
We construct an important function that will have piecewise constant jump conditions. Set
 \begin{align}\label{eq_vectortheta}
 \bTheta(z;\vec d;\vec v) = \bTheta(z) := \begin{bmatrix} \displaystyle \frac{\theta \left( \vec u(z) + \vec v - \vec d; \btau\right)}{\theta \left( \vec u(z) - \vec d; \btau\right)} & \displaystyle\frac{\theta \left( -\vec u(z) + \vec v - \vec d; \btau\right)}{\theta \left( -\vec u(z) - \vec d; \btau\right)} \end{bmatrix}, \qquad z \not\in \mathbb R.
 \end{align}
The first component function is nothing more than $\frac{\theta \left( \vec u(P) + \vec v - \vec d; \btau\right)}{\theta \left( \vec u(P) - \vec d; \btau\right)}$ restricted to the first sheet. The same is true for the second component function on the second sheet. The vector $\vec v$ is a free parameter. Then consider
\begin{align*}
 \vec u^+(z) + \vec u^-(z) = \left( 2 \sum_{k = 1}^{j-1} \int_{b_k}^{a_{k+1}} \dd \omega_\ell\right)_{\ell=1}^g = \left( \sum_{k = 1}^{j-1} \oint_{\mathfrak a_k} \dd \omega_\ell\right)_{\ell=1}^g = 2 \pi \ii \mathsf N, \qquad z \in [ a_j, b_j],
\end{align*}
for a vector $\mathsf N$ of zeros and ones. Then we compute
\begin{align*}
 \vec u^+(z) -\vec u^-(z) = \left( 2 \sum_{k = 1}^{j} \int_{a_k}^{b_{k}} \dd \omega_\ell\right)_{\ell=1}^g = \left( \oint_{\mathfrak b_{j}} \dd \omega_\ell\right)_{\ell=1}^g = \btau \vec e_{j}, \qquad z \in [ b_j, a_{j+1}].
\end{align*}
{Directly using this relation, it follows that, for $z \in [b_j, a_{j+1}]$},
\begin{align*}
 {
 \frac{\theta \left( \pm \vec u^+(z) + \vec v - \vec d; \btau\right)}{\theta \left( \pm \vec u^+(z) - \vec d; \btau\right)} = \frac{\theta \left( \pm \vec u^-(z) \pm \btau \vec e_j + \vec v - \vec d; \btau\right)}{\theta \left( \pm \vec u^-(z) \pm \btau \vec e_j- \vec d; \btau\right)} = \ee^{\mp v_j} \frac{\theta \left( \pm \vec u^-(z)+ \vec v - \vec d; \btau\right)}{\theta \left( \pm \vec u^-(z) - \vec d; \btau\right)}.}
\end{align*}
On the interval $(-\infty,a_1)$, $\vec u^+(z) = \vec u^-(z)$. And on $( b_{g+1},\infty)$,
\begin{align*}
 \vec u^+(z) - \vec u^-(z) = \left(\oint_\mathcal{C} \dd \omega_j\right)_{j=1}^g,
\end{align*}
where $\mathcal{C}$ is a clockwise-oriented simple contour that encircles $[ a_1, b_{g+1}]$. Because all the differentials $\dd \omega_j$ are of the form $P(z)/R(z)$ where $P$ is a degree $g-1$ polynomial and $R(z)/z^{g+1} \to 1$ as $z \to \infty$, it follows that $\oint_\mathcal{C} \dd \omega_j = 0$.
Thus, ignoring any poles $\bTheta$ may have, we find that $\bTheta$ satisfies the following jump conditions:
\begin{align*}
 \bTheta^+(z) = \begin{cases} \bTheta^-(z) \begin{bmatrix} 0 & 1 \\ 1 & 0 \end{bmatrix} ,& z \in ( a_j, b_j),\vspace{1.5pt}\\
 \bTheta^-(z) \begin{bmatrix} \ee^{-v_j} & 0 \\ 0 & \ee^{v_j} \end{bmatrix} ,& z \in ( b_j, a_{j+1}),\\
 \bTheta^-(z), & z \in (-\infty, a_1) \cup ( b_{g+1},\infty). \end{cases}
\end{align*}
Also, note that since $\vec u(\infty)$ is finite, $\bTheta$ is analytic at infinity.

Next, we must understand the poles of $\bTheta$. It is known that (see \cite{Algebro}, for example) if for $D = P_1 + \cdots + P_g$, $\theta( \vec u(P) - \mathcal A(D) - \vec k)$ is not identically zero,\footnote{This holds if $D$ is nonspecial.} then, counting multiplicities, $\theta( \vec u(P) - \mathcal A(D) - \vec k)$, has $g$ zeros on $\Sigma$. These zeros are then given by the points in the divisor~$D$. Consider the function
\[
\gamma(z)=\left[\prod_{j=1}^{g+1}\left(\frac{z- b_{j}}{z- a_{j}}\right)\right]^{1 / 4},
\]
analytic on $\mathbb{C} \setminus \bigcup_j [ a_j, b_j]$, with $\gamma(z) \sim 1, z \rightarrow \infty$. It can be shown that $\gamma - \gamma^{-1}$ has a single root~$z_j$ in $( b_j, a_{j+1})$ for $j =1,2,\dots,g$, while $\gamma + \gamma^{-1}$ does not vanish on $\mathbb C \setminus \bigcup_j [ a_j, b_j]$. So, define two divisors
\begin{align*}
 D_1 = \sum_{j=1}^g \pi_1^{-1}(z_j), \qquad D_2 = \sum_{j=1}^g \pi_2^{-1}(z_j).
\end{align*}
From \cite{DubrovinTheta} (see also \cite[Lemma 11.10]{TrogdonSOBook}), these divisors are nonspecial so that the $\theta$ functions we consider do not vanish identically.

Note that for $\vec d_1 := \mathcal A(D_1) + \vec k, $ the function $z \mapsto \theta(\vec u(z) - \vec d_1;\btau)$ has zeros at $z_j$, while the function $z \mapsto \theta(-\vec u(z) - \vec d_1;\btau)$ is non-vanishing. Similarly, for
$
 \vec d_2 := \mathcal A(D_2) + \vec k$,
the function $z \mapsto \theta(-\vec u(z) - \vec d_2;\btau)$ has zeros at $z_j$, while the function~${z \mapsto \theta(\vec u(z) - \vec d_2;\btau)}$ is non-vanishing.

Following \cite{Deift1997}, consider
\begin{align}\label{eq:L}
 \vec L(z;\vec v) = \begin{bmatrix} \left(\frac{\gamma(z) + \gamma(z)^{-1}}{2}\right) \bTheta_1(z;\vec d_2;\vec v) & \left(\frac{\gamma(z) - \gamma(z)^{-1}}{2\ii}\right) \bTheta_2(z;\vec d_2;\vec v) \\ \\
 \left(\frac{\gamma(z)^{-1} - \gamma(z)}{2\ii}\right) \bTheta_1(z;\vec d_1;\vec v) & \left(\frac{\gamma(z) + \gamma(z)^{-1}}{2}\right) \bTheta_2(z;\vec d_1;\vec v) \end{bmatrix},
\end{align}
which is analytic in $\mathbb C \setminus \bigcup_j [ a_j, b_j]$, with a limit as $z \to \infty$ and satisfies the jumps
\begin{align*}
 \vec L^+(z;\vec v) = \begin{cases} \vec L^-(z;\vec v) \begin{bmatrix} 0 & 1 \\ -1 & 0 \end{bmatrix}, & z \in ( a_j, b_j),\\ \\
 \vec L^-(z;\vec v) \begin{bmatrix} e^{-v_j} & 0 \\ 0 & e^{v_j} \end{bmatrix}, & z \in ( b_j, a_{j+1}),\\
 \vec L^-(z;\vec v), & z \in (-\infty, a_1) \cup ( b_{g+1},\infty). \end{cases}
\end{align*}
This follows because $\gamma^+(z) = \ii \gamma^-(z)$ for $z \in ( a_j, b_j)$ and therefore
\begin{align*}
 \gamma^+(z) + \bigl(\gamma(z)^{-1}\bigr)^+ &= \ii \bigl( \gamma^-(z) - \bigl(\gamma(z)^{-1}\bigr)^-\bigr),\\
 \gamma^+(z) - \bigl(\gamma(z)^{-1}\bigr)^+ &= \ii \bigl( \gamma^-(z) + \bigl(\gamma(z)^{-1}\bigr)^-\bigr).
\end{align*}
It is important to note that (\ref{eq:L}) was first used in \cite{Deift1997} and subsequently by many, see \cite{Chen2008,Deift1999b, Deift1999}, for example.

\section{Estimates of some exponential integrals}\label{a:estint}

In this appendix, we discuss how to estimate integrals of the form
\begin{gather*}
 \int_{E_\rho(a,b)} |z - a|^{-\gamma} |z - b|^{-\gamma} \big|\ee^{- m \mathfrak g(z)}\big| \dd A(z),\\
 E_\rho = \{z\mid |z + \sqrt{z+1}\sqrt{z -1}| < \rho \},\qquad \rho > 1,\\
 E_\rho(a,b) = M^{-1}(E_\rho), \qquad M^{-1}(z) = M^{-1}(z;a,b) = \frac{b - a}{2} x + \frac{b + a}{2},
\end{gather*}
as $m \to \infty$. From Proposition~\ref{p:ggrowth}, it follows that we can consider, for a new $m$,
\begin{align*}
I_m := \int_{E_\rho}& \big|z^2-1\big|^{-\gamma} |{\phi(z)}|^{-m} \dd A(z).
\end{align*}
We then set $z = 1/2\bigl( w + w^{-1}\bigr)$ for $1 < |w| < \rho$ and change variables.
{We have}
\begin{align*}
 \dd A(z) = \frac 1 4 \frac{\big|w^2 -1\big|^2}{|w|^4} \dd A(w).
\end{align*}
To finish the change of variables, we note that
\begin{align*}
 \varphi(z) = w, \qquad z^2 - 1 = \left[ \frac 1 2 \bigl( w - w^{-1} \bigr) \right]^2 = w^{-2} \left[ \frac 1 2 \bigl( w^2 - 1 \bigr) \right]^2.
\end{align*}
Therefore,
\begin{align*}
 I_m = \int_{1 < |w| < \rho} |w|^{2 \gamma -4 -m} 2^{{2}\gamma-2} \big|w^2 -1\big|^{2 -2\gamma } \dd A(w).
\end{align*}
For $\gamma \leq 1$, we have $I_m = \OO\bigl(m^{-1}\bigr)$, and therefore
\begin{align}\label{eq:exp_est}
 \int_{E_\rho(a,b)}& |z - a|^{-\gamma} |z - b|^{-\gamma} |\ee^{- m \mathfrak g(z)}| \dd A(z) = \OO\bigl(m^{-1}\bigr).
\end{align}

Now, suppose $0 \geq 2-2\gamma = -\sigma > -1$. We write $w = \rho \ee^{\ii \theta}$ and $\theta \in [-\pi/2,\pi/2]$, following \cite{Yattselev2023}, write
\begin{align*}
 |w - 1| = \sqrt{(\rho -1)^2 + 4 \rho \sin^2 ( \theta/2)} \geq C [ (\rho -1) + |\theta| ],
\end{align*}
for some $C > 0$. Consider
\begin{align*}
 \int_1^\rho \int_{-\pi/2}^{\pi/2} \rho^{-m} |w -1|^{-\sigma} \dd \rho \dd \theta \leq C^{-\sigma} \int_1^\rho \rho^{-m} \left(\int_{-\pi/2}^{\pi/2} \frac{\dd \theta}{((\rho -1) + |\theta|))^\sigma} \right) \dd \rho.
\end{align*}
For $\sigma < 1$, the $\theta$ integral produces a continuous function of $\rho$. Since it therefore must be bounded, we can conclude the integral is $O\bigl(m^{-1}\bigr)$. This argument can then be used for $w$ in the right and left half-planes separately, bounding $|w \pm 1|^{-\sigma} \leq 1$ for $\pm \operatorname{Re} w \geq 0$ to also conclude~\eqref{eq:exp_est} for~${\gamma < 3/2}$.

\subsection*{Acknowledgments}

This material is based upon work supported by NSF DMS-2306438. Any opinions, findings, and conclusions or recommendations expressed in this material are those of the author and do not necessarily reflect the views of the National Science Foundation. The author is very grateful to the anonymous referees whose comments not only improved this manuscript but caught inaccuracies.

\pdfbookmark[1]{References}{ref}
\LastPageEnding

\end{document}